\documentclass[10pt,twoside]{article}

\usepackage{amssymb,amsmath,makeidx}
\usepackage[T1]{fontenc}
\usepackage{latexsym}
\usepackage{exscale}
\usepackage{pb-diagram}
\usepackage[dvips]{graphicx}
\usepackage{xr}
\usepackage{hyperref,url}

\usepackage{xcolor,pgf,tikz}

\newcommand{\al}{\alpha}
\newcommand{\alvec}{{\vec{\al}}}

\newcommand{\alpr}{{\alpha^\prime}}

\newcommand{\alnod}{{\al_0}}

\newcommand{\aln}{{\al_n}}

\newcommand{\alnmin}{{\al_{n-1}}}
\newcommand{\almmin}{{\al_{m-1}}}

\newcommand{\alm}{{\al_m}}

\newcommand{\alplus}{{\al^+}}

\newcommand{\be}{\beta}

\newcommand{\bevec}{{\vec{\be}}}

\newcommand{\ga}{\gamma}

\newcommand{\Ga}{\Gamma}
\newcommand{\de}{\delta}

\newcommand{\De}{\Delta}
\newcommand{\Depr}{{\Delta^\prime}}
\newcommand{\Depre}{{\Delta^\prime_1}}
\newcommand{\Deprz}{{\Delta^\prime_2}}
\newcommand{\Depri}{{\Delta^\prime_i}}

\newcommand{\Deprk}{{\Delta^\prime_k}}
\newcommand{\Deprl}{{\Delta^\prime_l}}

\newcommand{\eps}{\varepsilon}

\newcommand{\epsn}{\varepsilon_0}
\newcommand{\epsom}{\varepsilon_{\Omega+1}}

\newcommand{\la}{\lambda}

\newcommand{\om}{\omega}

\newcommand{\Om}{\Omega}

\newcommand{\Omi}{{\Omega_i}}
\newcommand{\Omj}{{\Omega_j}}
\newcommand{\Omk}{{\Omega_k}}

\newcommand{\Omm}{{\Omega_m}}
\newcommand{\Omme}{\Om_{m+1}}
\newcommand{\Ommz}{\Om_{m+2}}
\newcommand{\Omkz}{\Om_{k+2}}
\newcommand{\Omke}{\Om_{k+1}}
\newcommand{\Omle}{\Om_{l+1}}
\newcommand{\Omiz}{\Om_{i+2}}
\newcommand{\Omie}{\Om_{i+1}}
\newcommand{\Omjz}{\Om_{j+2}}
\newcommand{\Omje}{\Om_{j+1}}

\newcommand{\rhopr}{{\rho^\prime}}

\newcommand{\si}{\sigma}
\newcommand{\tht}{\vartheta}
\newcommand{\Xipr}{{\Xi^\prime}}
\newcommand{\xipr}{{\xi^\prime}}

\newcommand{\ze}{\zeta}
\newcommand{\zepr}{{\zeta^\prime}}

\newcommand{\etapr}{{\eta^\prime}}

\newcommand{\etanod}{{\eta_0}}

\renewcommand{\phi}{\varphi}

\newcommand{\N}{{\mathbb N}}

\newcommand{\Hz}{{\mathbb P}}

\newcommand{\Lz}{{\mathbb L}}
\newcommand{\Ez}{{\mathbb E}}
\newcommand{\Ezone}{{\mathbb E}_1}
\newcommand{\On}{{\mathrm{Ord}}}

\newcommand{\CNF}{{\mathrm{\scriptscriptstyle{CNF}}}}
\newcommand{\ANF}{{\mathrm{\scriptscriptstyle{ANF}}}}

\newcommand{\NF}{{\mathrm{\scriptscriptstyle{NF}}}}

\newcommand{\Lim}{\mathrm{Lim}}

\newcommand{\Image}{\mathrm{Im}}

\newcommand{\logend}{{\mathrm{logend}}}
\newcommand{\sumend}{{\mathrm{end}}}

\newcommand{\klex}{<_\mathrm{\scriptscriptstyle{lex}}}

\newcommand{\thtm}{\tht_m}
\newcommand{\thtme}{\tht_{m+1}}
\newcommand{\thtn}{\tht_n}
\newcommand{\thtnod}{\tht_0}
\newcommand{\thte}{\tht_1}
\newcommand{\thtz}{\tht_2}
\newcommand{\thti}{\tht_i}
\newcommand{\thtie}{\tht_{i+1}}
\newcommand{\thtiz}{\tht_{i+2}}
\newcommand{\thtj}{\tht_j}
\newcommand{\thtje}{\tht_{j+1}}
\newcommand{\thtk}{\tht_k}
\newcommand{\thtke}{\tht_{k+1}}

\newcommand{\thtt}{\tht^\tau}

\newcommand{\T}{{\operatorname{T}}}

\newcommand{\Tm}{{\operatorname{T}_m}}

\newcommand{\Tt}{{\operatorname{T}^\tau}}

\newcommand{\Ts}{{\operatorname{T}^\si}}

\newcommand{\stark}{{\ast_k}}

\newcommand{\start}{{\ast^\tau}}

\newcommand{\starnod}{{\ast_0}}

\newcommand{\thetam}{{\theta_m}}
\newcommand{\thetame}{{\theta_{m+1}}}

\newcommand{\htt}{{\operatorname{ht}_\tau}}
\newcommand{\htomi}{{\operatorname{ht}_{\Omi}}}

\newcommand{\htomje}{{\operatorname{ht}_{\Omje}}}

\newcommand{\htomme}{{\operatorname{ht}_{\Omme}}}

\newcommand{\mc}{{\operatorname{mc}}}

\newcommand{\pioneonecanod}{\Pi^1_1{\operatorname{-CA}_0}}

\newcommand{\id}{\operatorname{id}}

\newcommand{\qed}{\mbox{ }\hfill $\Box$\vspace{2ex}}
\newcommand{\imp}{\Rightarrow}
\newcommand{\aeq}{\Leftrightarrow}
\newcommand{\andsp}{\:\&\:}
\newcommand{\veesp}{\:\vee\:}

\newcommand{\set}[2]{\{ #1 \:|\: #2\}}

\newcommand{\singleton}[1]{\{ #1 \}}

\newlength{\hilflh}
\newcommand{\hilfminus}[1]{
  \settowidth{\hilflh}{$#1-$}\mbox{$#1-\hspace{-0.5\hilflh}
  \makebox[0pt]{\raisebox{0.24\hilflh}{$#1\cdot$}}\hspace{0.5\hilflh}$}}
\newcommand{\minusp}{\mathbin{\mathchoice {\hilfminus{\displaystyle}}
  {\hilfminus{\textstyle}}{\hilfminus{\scriptstyle}}
  {\hilfminus{\scriptscriptstyle}}}}

\newtheorem{theo}{Theorem}[section]
\newtheorem{cor}[theo]{Corollary}
\newtheorem{lem}[theo]{Lemma}
\newtheorem{defi}[theo]{Definition}
\newtheorem{prop}[theo]{Proposition}

\newtheorem{rmk}[theo]{Remark}
\newtheorem{claim}[theo]{Claim}

\newcommand{\oneinf}{1^\infty}

\newcommand{\tauinf}{\tau^\infty}

\newcommand{\dom}{\mathrm{dom}}

\newcommand{\kpr}{{k^\prime}}

\def\vec#1{\mathchoice{\mbox{\boldmath$\displaystyle#1$}}
{\mbox{\boldmath$\textstyle#1$}}
{\mbox{\boldmath$\scriptstyle#1$}}
{\mbox{\boldmath$\scriptscriptstyle#1$}}}

\newcommand{\Romannumeral}[1]{\uppercase\expandafter{\romannumeral #1\relax}}

\newcommand{\Si}{\Sigma}
\newcommand{\Sipr}{\Sigma^\prime}
\newcommand{\Sinod}{\Sigma_0}
\newcommand{\stari}{{\star_{i}}}
\newcommand{\starie}{{\star_{i+1}}}
\newcommand{\starj}{{\star_{j}}}
\newcommand{\starje}{{\star_{j+1}}}

\newcommand{\chiomie}{\chi^{\Omie}}
\newcommand{\chiomiz}{\chi^{\Omiz}}
\newcommand{\Tomi}{\T^{\Omi}}
\newcommand{\Tomie}{\T^{\Omie}}
\newcommand{\Tomj}{\T^{\Omj}}
\newcommand{\chiomje}{\chi^{\Omje}}
\newcommand{\chiomjz}{\chi^{\Omjz}}

\newcommand{\chiomke}{\chi^{\Omke}}
\newcommand{\chiomkz}{\chi^{\Omkz}}
\newcommand{\chiomle}{\chi^{\Omle}}
\newcommand{\Tomke}{\T^{\Omke}}
\newcommand{\Tomk}{\T^{\Omk}}
\newcommand{\Tomm}{\T^{\Omm}}

\newcommand{\domf}{\mathrm{d}}
\newcommand{\domfinv}{\mathrm{d}^{-1}}

\newcommand{\Piset}{{\operatorname{P}_i}}

\newcommand{\ual}{\underline{\al}}

\newcommand{\ube}{\underline{\be}}

\newcommand{\uGa}{\underline{\Ga}}

\newcommand{\uDe}{\underline{\De}}
\newcommand{\uDee}{\underline{\De_1}}
\newcommand{\uDepre}{\underline{\Depre}}

\newcommand{\uSi}{\underline{\Si}}

\newcommand{\Ttcirc}{{\mathring{\T}^\tau}}

\def\rddots#1{\mbox{ }\cdot^{\cdot^{\cdot^{#1}}}}

\def\norm#1{|\!|{#1}|\!|}

\newcommand{\Omom}{{\Om_\om}}

\newcommand{\Cqi}{{\bar{\operatorname{C}}_i}}

\newcommand{\Cqj}{{\bar{\operatorname{C}}_j}}
\newcommand{\Cqk}{{\bar{\operatorname{C}}_k}}
\newcommand{\Cqm}{{\bar{\operatorname{C}}_m}}

\newcommand{\thtq}{{\bar{\tht}}}
\newcommand{\thtqnod}{\thtq_0}
\newcommand{\thtqe}{\thtq_1}
\newcommand{\thtqz}{\thtq_2}
\newcommand{\thtqt}{{\bar{\tht}^\tau}}

\newcommand{\thtqi}{{\bar{\tht}_i}}
\newcommand{\thtqie}{{\bar{\tht}_{i+1}}}
\newcommand{\thtqj}{{\bar{\tht}_j}}
\newcommand{\thtqje}{{\bar{\tht}_{j+1}}}
\newcommand{\thtqk}{{\bar{\tht}_k}}
\newcommand{\thtqke}{{\bar{\tht}_{k+1}}}
\newcommand{\thtqm}{{\bar{\tht}_m}}

\newcommand{\thtqme}{{\bar{\tht}_{m+1}}}

\newcommand{\thtqn}{{\bar{\tht}_n}}

\newcommand{\domthtqi}{{\operatorname{dom}({\thtqi})}}
\newcommand{\domthtqj}{{\operatorname{dom}({\thtqj})}}

\newcommand{\domthtqk}{{\operatorname{dom}({\thtqk})}}

\newcommand{\domthtqm}{{\operatorname{dom}({\thtqm})}}
\newcommand{\domthtqme}{{\operatorname{dom}({\thtqme})}}

\newcommand{\Tq}{{\bar{\operatorname{T}}}}
\newcommand{\Tqm}{{\Tq_m}}

\newcommand{\Tqt}{{\Tq^\tau}}

\newcommand{\Tqomj}{\Tq^{\Omj}}
\newcommand{\Tqomi}{\Tq^{\Omi}}
\newcommand{\Tqomie}{\Tq^{\Omie}}
\newcommand{\Tqomke}{\Tq^{\Omke}}
\newcommand{\Tqomm}{\Tq^{\Om_m}}
\newcommand{\Tqomd}{\Tq^{\Om_d}}

\newcommand{\Kt}{\operatorname{K}^\tau}
\newcommand{\tstar}{{\tau\ast}}

\newcommand{\Ktstar}{\operatorname{K}^\tstar}

\newcommand{\ktstar}{\operatorname{k}^\tstar}

\newcommand{\starti}{{\ast^\tau_i}}

\newcommand{\ft}{\operatorname{f}^\tau}

\newcommand{\gt}{\operatorname{g}^\tau}
\newcommand{\gOmie}{\operatorname{g}^{\Omie}}

\newcommand{\Tqtcirc}{{{\mathring{\bar{\T}}}^\tau}}

\newcommand{\plainit}{{\operatorname{it}}}
\newcommand{\itm}{{\plainit_m}}
\newcommand{\itme}{{\plainit_{m+1}}}
\newcommand{\itk}{{\plainit_k}}
\newcommand{\itke}{{\plainit_{k+1}}}
\newcommand{\iti}{{\plainit_i}}
\newcommand{\itie}{{\plainit_{i+1}}}

\newcommand{\plainrt}{{\operatorname{rt}}}
\newcommand{\rtm}{{\plainrt_m}}
\newcommand{\rtme}{{\plainrt_{m+1}}}
\newcommand{\rti}{{\plainrt_i}}
\newcommand{\rtie}{{\plainrt_{i+1}}}
\newcommand{\rtj}{{\plainrt_j}}
\newcommand{\rtjmin}{{\plainrt_{j-1}}}
\newcommand{\rtje}{{\plainrt_{j+1}}}
\newcommand{\rtk}{{\plainrt_k}}
\newcommand{\rtke}{{\plainrt_{k+1}}}
\newcommand{\rtl}{{\plainrt_l}}

\begin{document}

\title{Fundamental sequences based on localization 
\footnote{Extended version of \cite{W24} containing complete proofs and generalization to simultaneously defined collapsing functions.\\
To appear in: Annals of Pure and Applied Logic (2025), https://doi.org/10.1016/j.apal.2025.103658}}

\author{Gunnar Wilken\\
Okinawa Institute of Science and Technology\\
904-0495 Okinawa, Japan\\
{\tt wilken@oist.jp}
}

\maketitle

\begin{abstract}
Building on Buchholz' assignment for ordinals below Bachmann-Howard ordinal, see \cite{B03},
we introduce systems of fundamental sequences for two kinds of relativized $\vartheta$-function-based notation systems of strength
$\pioneonecanod$ and prove Bachmann property for these systems, which is essential for monotonicity properties of subrecursive hierarchies 
defined on the basis of fundamental sequences. The central notion of our construction is the notion of \emph{localization}, which was introduced in 
\cite{W07a}. 

The first kind of \emph{stepwise defined} $\tht$-functions over ordinal addition as basic function fits the framework of the ordinal arithmetical toolkit 
developed in \cite{W07a}, whereas the second kind of $\vartheta$-functions is defined \emph{simultaneously} and will allow for further generalization
to larger proof-theoretic ordinals, see \cite{WW11}.

The systems of fundamental sequences given here enable the investigation of fundamental sequences and independence phenomena also in the 
context of patterns of resemblance, an approach to ordinal notations that is both semantic and combinatorial and was first introduced by Carlson in  \cite{C01} and further analyzed in \cite{W06,W07b,W07c,CWa}. 

Our exposition is put into the context of the abstract approach to fundamental sequences developed by Buchholz, Cichon, and Weiermann in \cite{BCW94}.
The results of this paper will be applied to the theory of Goodstein sequences, extending results of \cite{FDW24}.
\end{abstract}



\section{Introduction}

In this article we are going to define systems of fundamental sequences that satisfy Bachmann property for two kinds of ordinal notation systems
which we will briefly review in the preliminaries section. The fundamental sequences chosen here coincide with those Buchholz gave in \cite{B03}
for the initial segment up to Bachmann-Howard ordinal. Here, however, we work with syntactically different terms that reach Takeuti ordinal, the 
proof-theoretic ordinal of subsystem $\pioneonecanod$ of second order number theory, and come with a straightforward relativization as introduced in \cite{W07a} and \cite{WW11}.

For independent work also building on Buchholz' fundamental sequences in \cite{B03}, but in the direction of Goodstein sequences of strength up to 
Bachmann-Howard ordinal, see \cite{FDW24}. Results of this paper will also be beneficial for an elegant definition of (maximal) Goodstein sequences of 
a strength at least matching $\pioneonecanod$, which is subject to work in progress with the authors of \cite{FDW24}, who observed a miniaturized 
analogue of the mechanism of localization that was introduced in \cite{W07a} and used here to determine fundamental sequences that satisfy 
Bachmann property.

We are going to work with notations used in the analysis of patterns of resemblance of orders $1$ and $2$. Patterns of resemblance were first introduced by Carlson, 
see \cite{C01}, and constitute an approach to ordinal notations that is both semantic and combinatorial. 
We thus further extend the ordinal arithmetical results elaborated as a byproduct when analyzing patterns of resemblance that were laid out in \cite{W07a}, 
applied in \cite{W07b} and  \cite{W07c}, slightly extended in \cite{CWa}, and put to full usage in \cite{W18,W21}.\footnote{For an overview and gentle introduction 
to the analysis of patterns, see \cite{W06} and \cite{W20}.} 
Notations derived from injective $\vartheta$-functions as in \cite{W07a} provide unique terms for ordinals, which is helpful when ordinal arithmetical 
analysis is the main focus of usage of such notations. 
Regarding patterns of resemblance the present article also sets the stage for a theory of \emph{pattern related fundamental sequences} that gives rise 
to independence phenomena related to patterns.

While notations of \cite{W07a}, which are defined in a stepwise manner, are sufficient for ordinal segments characterizing Takeuti ordinal, 
it is desirable to also consider a notation system that still provides unique terms for ordinals like notations from \cite{W07a} but allows for natural 
generalizability. Notations of such kind involve simultaneously defined collapsing functions. These were investigated in \cite{WW11} after an adaptation
and restriction of notations used by Wilfried Buchholz in \cite{B81}, and an order isomorphism with stepwise defined notations from \cite{W07a} was 
given in Definition 3.9 of \cite{WW11}, thus providing translations between the two notation systems. 

We therefore extend our results from \cite{W24}, the proofs of which are given here, to the more generalizable notations considered in \cite{WW11}. 
Simultaneously defined collapsing functions, denoted as $\thtq$-functions, have more segmented domains, which motivated the earlier usage of 
stepwise defined collapsing functions, denoted as $\tht$-functions.
The order isomorphism given in \cite{WW11} will be restated in concise form in the preliminaries section and in essence provides for each $m<\om$ a
strictly increasing mapping $\rtm$ of the fragmented domain of (simultaneously defined) function $\thtqm$ onto the domain of (stepwise defined) function 
$\thtm$, which is a segment of ordinals denoted by $\theta_{m+1}$, cf.\ Corollary 3.7 of \cite{WW11}. $\rtm$ is therefore the Mostowski collapse of 
$\domthtqm$. For details of the stepwise definition process of $\tht$-functions, on which we do not rely in this article, we refer the reader to \cite{W07a}.

\section{Preliminaries}

For basics on the arithmetic of ordinal numbers we refer the reader to Pohlers' book \cite{P09}.
We denote the class of \emph{additive principal numbers}, i.e.\ ordinals greater than $0$ that are closed under ordinal addition, by $\Hz$. Such ordinals are 
characterized as the image of $\om$-exponentiation, where $\om$ denotes the least infinite ordinal. Likewise, we 
denote the class $\Lim(\Hz)$ of limits of additive principal numbers by $\Lz$, and the class of epsilon numbers by $\Ez$ (epsilon numbers are ordinals
closed under $\om$-exponentiation).

Any ordinal $\al\not\in\Hz\cup\{0\}$ is called \emph{additively decomposable}, and any $\al$ can be written in \emph{additive normal form}
$\al=_\ANF\al_1+\ldots+\al_n$ where $\al_1,\ldots,\al_n$ is a weakly decreasing sequence of additive principal numbers. We will often use the notation
$\al=_\NF\be+\ga$ for additively decomposable $\al$, which means that $\be,\ga>0$, $\ga$ is additively indecomposable (additive principal), 
and $\be$ is minimal such that $\al=\be+\ga$, hence $\be=\al_1+\ldots+\al_{n-1}$ and $\ga=\aln$. We will also use the notation $\sumend(\al)$ for $\ga$, 
and clearly $\sumend(\al)=\al$ if $\al\in\Hz\cup\{0\}$.
We write $\mc(\al):=\al_1$ for the maximal (additive) component of $\al$ and set $\mc(0):=0$.
For ordinals $\al,\be$ such that $\al\le\be$ we write $-\al+\be$ for the unique $\ga$ such that $\al+\ga=\be$. And for $k<\om$ we write 
$\al\minusp k$ for the least ordinal $\alpr$ such that $\al$ can be written as $\al=\alpr+l$ where $l\le k$.

For $\al\not\in\Ez\cup\{0\}$ the \emph{Cantor normal form} representation of $\al$ is useful and indicated as 
$\al=_\CNF\om^{\al_1}+\ldots+\om^{\al_n}$ 
where $\al>\al_1\ge\ldots\ge\al_n$. For completeness, in the case $\al=0$ we have $n=0$, and clearly $\al=\om^\al$ for $\al\in\Ez$.
We will use the notation $\logend(\al)$ for the ordinal $\al_n$. Clearly, $\logend(0)=0$ and $\logend(\al)=\al$ for $\al\in\Ez$.

\subsection{\boldmath Stepwise defined notation systems $\Tt$}\label{Ttsubsec}

For $\tau\in\Ezone$, i.e.\ $\tau$ is either equal to $1$ or any epsilon number, 
we defined the relativized notation system $\Tt$ in Subsection 2.1.2 of \cite{W21}.
$\Tt$ is built up over a sequence $\tau=\Om_0=\thtnod(0),\Om_1=\thte(0), \Om_2=\tht_2(0), \ldots$,
where $\Om_1,\Om_2,\ldots$ is a strictly increasing sequence of regular ordinals such that $\Om_1>\tau$.
The canonical choice is to assume that $\tau$ is countable (recursive) and that $\Omi=\aleph_i$ for $0<i<\om$.

Terms in systems $\Tt$ are composed of parameters below $\tau$, ordinal summation, and $\thtj$-functions
for $j<\om$, where $\thtnod$ is also written as $\thtt$. 
$\thtj$-functions uniquely denote ordinals closed under ordinal addition in the segment $[\Omj,\Omje)$.
The only restriction for the application of a $\thtj$-function
to a term in the system is that the argument be strictly below $\Om_{j+2}$.
The operation $\cdot^\starj$, cf.\ Definition 2.24 of \cite{W21} 
searches a $\Tt$-term for its $\thtj$-subterm of largest ordinal value, but under the restriction to treat $\thti$-subterms 
for $i<j$ as atomic. 
If such largest $\thtj$-subterm does not exist, $\cdot^\starj$ returns $0$, rather than any value from the interval $(0,\Omj)$.
We also write ${}^\start$ instead of ${}^\starnod$ if $\Om_0=\tau$.

By slight abuse of notation we can consider notation systems $\Tomie$ to be systems relativized 
to the initial segment $\Omie$ of ordinals and built up over 
$\Omie=\thtie(0), \Om_{i+2}=\tht_{i+2}(0),\ldots$, i.e.\ without renaming the indices of $\tht$-functions.
Results on $\Tt$ from \cite{W21} then directly carry over to such systems $\Tomie$ and the corresponding $\tht$-functions.
Note that $\Tm$ in Definition 3.22 of \cite{W07a} is equal to $\Tomm$ as introduced here (up to renaming indices of $\tht$-functions),
since it is defined to be the closure of parameters from below $\Omm$ under $+$ and functions $\thtk$ for $k\ge m$ and arguments 
below $\Omkz$.
Fixing domains $\thetame$ of $\thtm$-functions we have the following theorem. 
\begin{theo}[cf.\ 3.23 of \cite{W07a}] For $m<\om$ we have 
\[\thetam:=\Tm\cap\Omme=\sup_{n\ge m}\thtm(\cdots(\thtn(0))\cdots).\]
\end{theo}

For convenience, we cite a well-known useful proposition regarding $\tht$-functions  
that we can use as an alternative definition of $\tht$-functions in order to shortcut the definition process
that was fully carried out in \cite{W07a} (see also Subsection 2.1 of \cite{W21}).
\begin{prop}\label{thetaprop} For $\al\in\Tt\cap\Om_{j+2}$ we have
\[\thtj(\al)=\min\{\theta\ge\Omj\mid\al^\starj<\theta\mbox{ and }
  \forall\:\be\in\Tomj\cap\al\left(\be^\starj<\theta\to\thtj(\be)<\theta\right)\}.\]
\end{prop}

The following well-known comparison lemma of $\tht$-terms is rooted in the pivotal property $\be^\star<\tht(\be)$ 
and was stated and proven in Lemma 3.30 of \cite{W07a} (see also Lemma 2.26 of \cite{W21}).
\begin{prop}\label{thetacompprop} For $\al,\be,\ga\in\Tt\cap\Om_{j+2}$ we have $\be^\starj<\thtj(\be)$ and
\[\thtj(\al)<\thtj(\ga)\quad\aeq\quad\left(\al<\ga\;\wedge\;\al^\starj<\thtj(\ga)\right)\;\vee\;\thtj(\al)\le\ga^\starj.\]
\end{prop}
Note that in the above two  propositions we could as well regard $\al$, respectively $\al$, $\be$, and $\ga$, 
as terms of $\T^{\Omi}$ for any $i\le j$.

\emph{Convention:} We will make use of a convention also used throughout \cite{W07a,W21}, namely that when writing the argument of a 
$\thtj$-function, $j<\om$, as $\De+\eta$, more generally, using a sum of a capital Greek letter and a lower case Greek letter, 
we automatically mean that $\De$ is a multiple of $\Omje$ (possibly $0$), while $\eta<\Omje$.

Informally, for an ordinal $\al=\thtt(\De+\eta)$ the term $\De$ indicates the \emph{fixed point level} of $\al$, in arithmetized form. 
This mechanism creates the strength of ordinal notations using functions like $\tht$, as e.g.\ opposed to the various versions of 
Veblen functions $\phi$ of increasing arity. An easy observation into this direction is given by the following lemma.

\begin{lem}[cf.\ 4.3 of \cite{W07a}, 2.29 of \cite{W21}]\label{epscharlem} For $\al=\thtj(\De+\eta)\in\Tt$ we have $\al\in\Ez^{>\Omj}$ 
if and only if $\De>0$.
\end{lem}

Note also that the function $\si\mapsto\thtj(\De+\si)$ is strictly increasing for $\si<\Omje$, hence 
$\eta\le\sup_{\si<\eta}\thtj(\De+\si)$, and let the condition $F_j(\De,\eta)$ be defined as
\[F_j(\De,\eta)\quad:\aeq\quad\eta={\sup_{\si<\eta}}^+\thtj(\De+\si),\]
where the usage of proper supremum excludes the unintended equality
$1=\sup_{\si<1}\thtnod(\si)=\thtnod(0)$ (for $\tau=1$). 
\footnote{An acknowledgement to Samuel Alexander for pointing out this flaw in earlier publications like \cite{W21}.}
An easy application of the above proposition is the following. 
\begin{prop}\label{continuityprop} For $\al=\thtj(\De+\eta)$ such that
$\eta\in\Lim$ and $\neg F_j(\De,\eta)$ we have the continuous, strictly monotonically increasing approximation
\begin{equation}\label{continuityeq} 
\al=\sup_{\si<\eta}\thtj(\De+\si).
\end{equation}
\end{prop}
Another very useful proposition in this context is Lemma 4.4 of \cite{W07a} (Lemma 2.30 of \cite{W21}):
\begin{prop}\label{fixpcondprop}
For $\al=\thtj(\De+\eta)$ we have
\begin{equation}\label{fixpcondequiv}
F_j(\De,\eta)\quad\aeq\quad\eta\mbox{ is of a form }\eta=\thtj(\Ga+\rho)\mbox{ such that }\Ga>\De\mbox{ and }\eta>\De^\starj.
\end{equation}
\end{prop}

\subsection{\boldmath Simultaneously defined notation systems $\Tqt$}\label{Tqtsubsec}

This subsection as well as Subsection \ref{orderisosubsec} are preliminaries only for Section \ref{fundseqTtbarsec} and therefore 
can be skipped at first reading.
For the reader's convenience we begin with the formal definition of simultaneously defined $\vartheta$-functions, denoted here by indexed 
$\thtq$ to avoid clash of notations.
We keep conventions regarding $\Om_0=\tau$ and the sequence $(\Omi)_{0<i<\om}$ of regular ordinals strictly above $\tau$, 
and set $\Omom:=\sup_{n<\om}\Om_n$. 

\begin{defi}[Def. 2.1 of \cite{WW11}] Simultaneously for all $i<\om$ we define sets of ordinals $\Cqi(\al,\be)$ where $\be\le\Om_{i+1}$ and,
whenever $\al\in\Cqi(\al,\Om_{i+1})\cap\Omom$, ordinals $\thtqi(\al)$ by recursion on $\al\le\Omom$.
For each $\be\le\Om_{i+1}$ the set $\Cqi(\al,\be)$ is defined inductively by
\begin{enumerate}
\item $\Om_i\cup\be\subseteq\Cqi(\al,\be)$
\item $\xi,\eta\in\Cqi(\al,\be)\;\imp\;\xi+\eta\in\Cqi(\al,\be)$
\item $j\in[i,\om)\andsp\xi\in\Cqj(\xi,\Om_{j+1})\cap\Cqi(\al,\be)\cap\al\;\imp\;\thtqj(\xi)\in\Cqi(\al,\be).$
\end{enumerate}
For $\al\in\Cqi(\al,\Om_{i+1})$ we define
\[\thtqi(\al):=\min\left(\set{\xi<\Om_{i+1}}{\al\in\Cqi(\al,\xi)\andsp
                         \Cqi(\al,\xi)\cap\Om_{i+1}\subseteq\xi}\cup\singleton{\Om_{i+1}}\right).\]
We then define
\[\domthtqi:=\set{\al<\Omom}{\al\in\Cqi(\al,\Om_{i+1})}.\]
We will also write $\thtqt$ instead of $\thtq_0$ in order to make the setting of relativization at the
$0$-th level visible.
\end{defi}

Clearly, the $\Cqi$-sets defined above have the straightforward monotonicity and continuity properties in both arguments and the index $i$ 
(cf.\ Lemma 2.3 of \cite{WW11}). 
As in the case of $\vartheta$-functions we have $\thtqi(0)=\Omi$ for all $i<\om$, and by the usual regularity argument,
for all $\al\in\domthtqi$ we have $\thtqi(\al)<\Om_{i+1}$, i.e.\ the functions $\thtqi$ are collapsing functions (cf.\ Lemma 2.4 of \cite{WW11}). 
As a corollary, cf.\ Corollary 2.5 of \cite{WW11}, we obtain that the sets $\Cqi$ are also monotonically increasing in their index $i$, 
in particular \[\domthtqi\subseteq\domthtqj\] for $i\le j$,
that for all $\al\in\domthtqi$ we have \[\thtqi(\al)=\Cqi(\al,\thtqi(\al))\cap\Om_{i+1}\in\Hz\cap[\Om_i,\Om_{i+1}),\]
that images of functions $\thtqi$ and $\thtqj$ for $i\not=j$ are disjoint, and that functions $\thtqi$ are injective. $\Cqi$-sets are closed under additive
decomposition and satisfy the property 
\[\ga\in\domthtqj\andsp j\ge i\andsp\be\le\thtqj(\ga)\in\Cqi(\al,\be)\quad\imp\quad\ga\in\Cqi(\al,\be)\cap\al,\]
which now follows by induction on the build-up of $\Cqi(\al,\be)$, see Lemma 2.6 of \cite{WW11}. This allows one to prove the following theorem.

\begin{theo}[2.8 of \cite{WW11}]\label{segmenttheo} For every $\al\le\Omom$ and every $\be\le\Om_{i+1}$ we have 
\[\Cqi(\al,\be)\cap\Om_{i+1}\in\On.\]
For every $\be\in\Om_{i+1}-\Cqi(\al,\be)$ we have $\Cqi(\al,\be)\cap\Om_{i+1}=\be$.
\end{theo}

\begin{defi}[cf.\ 2.9 of \cite{WW11}] For any $m<\om$ we define \[\Tqm:=\Cqm(\Omom,0).\]
We set $\Tqt:=\Tq_0$ in order to explicitly show the setting of relativization.
Again, by slight abuse of notation we can identify notation systems $\Tqomm$ and $\Tqm$ to be systems relativized 
to the initial segment $\Om_m$ of ordinals and built up over 
$\Om_m=\thtqm(0), \Om_{m+1}=\thtq_{m+1}(0),\ldots$, i.e.\ without renaming the indices of $\thtq$-functions in  $\Tqomm$.
\end{defi}

\begin{cor}[2.10 of \cite{WW11}] For every $m<\om$ we have $\Tqomm\cap\Omme=\Tqm\cap\Omme\in\On$.
\end{cor}

The following definition, cf.\ \cite{B81}, is helpful for comparison of $\thtq$-terms and in order to characterize domains $\dom(\thtqi)$. 
\begin{defi}[4.1 of \cite{WW11}] The sets of subterms $\Kt_i(\al)$ and $\Ktstar_i(\al)$ for $\al\in\Tqt$ and $i<\om$ are
defined by the following clauses.
\begin{enumerate}
\item $\Kt_i(\al):=\Ktstar_i(\al):=\emptyset$ for $\al<\tau$.
\item $\Kt_i(\al):=\Kt_i(\xi)\cup\Kt_i(\eta)$ and $\Ktstar_i(\al):=\Ktstar_i(\xi)\cup\Ktstar_i(\eta)$ for $\al=_\NF\xi+\eta$.
\item $\Kt_i(\al):=\left\{\begin{array}{cl}
                  \emptyset&\mbox{if }j<i\\
                  \singleton{\al}&\mbox{if }j=i\\
                  \Kt_i(\xi)&\mbox{if }j>i
                  \end{array}\right.
         \quad\mbox{ and }\quad
         \Ktstar_i(\al):=\left\{\begin{array}{cl}
                  \emptyset&\mbox{if }j<i\\
                  \singleton{\xi}\cup\Ktstar_i(\xi)&\mbox{if }j\ge i
                  \end{array}\right.
         \quad\mbox{ for } \al=\thtqj(\xi)$.
\end{enumerate}
We also set \[\al^\starti:=\max\left(\Kt_i(\al)\cup\singleton{0}\right) \quad\mbox{ and }\quad \ktstar_i(\al):=\max\left(\Ktstar_i(\al)\cup\singleton{0}\right).\]
Whenever confusion concerning the setting $\tau$ of relativization is unlikely we simply write
$\al^\stari$ instead of $\al^{\starti}$. Instead of $\al^{\starnod}$ we also write $\al^\start$.
\end{defi}

An induction on the build-up of $\Tqt$ now establishes the following lemma, which entails the subsequent characterization of $\domthtqm$
as a simple condition on subterms and the familiar lemma for comparison of $\thtqi$-terms. 

\begin{lem}[4.2 of \cite{WW11}]\label{tqcqkcharlem}
Let $\ga\in\Tqt$ and $k<\om$. For $\al\le\Omom$ and $\be<\Omke$ we have
\[\ga\in\Cqk(\al,\be)\aeq\Kt_k(\ga)\subseteq\Cqk(\al,\be)\andsp\Ktstar_{k+1}(\ga)<\al.\]
In particular, for every $\ga\in\Tqt\cap\domthtqk$ we have \[\ga^\stark<\thtqk(\ga).\]
\end{lem}

\begin{lem}[4.3 of \cite{WW11}]\label{tqdomthtqmlem}
$\al\in\domthtqm\aeq\Ktstar_{m+1}(\al)<\al$ for $\al\in\Tqt$.
\end{lem}

\begin{lem}[2.11 of \cite{WW11}] $\Tqm$ is inductively characterized as the closure of $\Om_m$ under $+$ and
$(i,\al)\longmapsto\thtqi(\al)$ where $i\in[m,\om)$ and $\al\in\domthtqi$, i.e.\ $\Ktstar_{i+1}(\al)<\al$.
\end{lem}

\begin{prop}[4.4 of \cite{WW11}] For $\al,\be\in\Tqt\cap\domthtqi$ we have
\[\thtqi(\al)<\thtqi(\be)\aeq\left(\al<\be\andsp\al^\stari<\thtqi(\be)\right)\veesp\thtqi(\al)\le\be^\stari.\]
\end{prop}

\emph{Convention:} For any $i<\om$ and $\De\in\domthtqi$ such that $\Omie\mid\De$ we have $\De+\eta\in\domthtqi$ for all $\eta<\Omie$, and we adopt 
the convention introduced for $\tht$-functions that writing $\thtqi(\De+\eta)$ implies that $\eta<\Omie\mid\De$ \emph{and} $\De\in\domthtqi$. 
Note that for any suitable term of the form $\De+\eta$ where $\eta<\Omie\mid\De$ we have $\De\in\domthtqi$ 
if and only if $\De+\eta\in\domthtqi$, cf. Lemma 2.16 of \cite{WW11}.
Note further that we no longer have $\De<\Omiz$, as $\thtq$-functions are not stepwise collapsing but rather have domains unbounded in $\Omom$.
We will therefore sometimes write $\thtqi$-terms in the form $\thtqi(\Xi+\De+\eta)$ where $\eta<\Omie\mid\De<\Omiz\mid\Xi$ to indicate the initial sum
of terms $\ge\Omiz$ of the argument of a $\thtqi$-function. Clearly, any of $\eta,\De,\Xi$ in this denotation can be $0$. 

For $\De\in\domthtqi$ the mapping $\eta\mapsto\thtqi(\De+\eta)$ is strictly increasing, and Lemma \ref{epscharlem} carries over,
see part c) of Lemma 2.12 of \cite{WW11}. 
As before we use the abbreviation $F_j(\De,\eta)$ for the fixed point condition, modified in the context of $\thtq$-functions to 
\[F_j(\De,\eta)\quad:\aeq\quad\De\in\domthtqj\quad\mbox{ and }\quad\eta={\sup_{\si<\eta}}^+\thtqj(\De+\si),\]
and for completeness we state the analogues of Propositions \ref{thetaprop} and \ref{fixpcondprop}.
\begin{prop}\label{thetaqprop} For $\al\in\Tqt\cap\domthtqj$ we have
\[\thtqj(\al)=\min\{\theta\ge\Omj\mid\al^\starj<\theta\mbox{ and }
  \forall\:\be\in\Tqomj\cap\domthtqj\cap\al\left(\be^\starj<\theta\to\thtqj(\be)<\theta\right)\}.\]
\end{prop}

The proof of the following proposition is a straightforward adaptation of the proof given for Lemma 4.4 of \cite{W07a}.
\begin{prop}\label{fixpcondqprop}
For $\al=\thtqj(\De+\eta)$ we have
\begin{equation}
F_j(\De,\eta)\quad\aeq\quad\eta\mbox{ is of a form }\eta=\thtqj(\Ga+\rho)\mbox{ such that }\Ga>\De\mbox{ and }\eta>\De^\starj.
\end{equation}
\end{prop}

The height of terms in $\Tqt$ is now defined as in \cite{W07a} for terms in $\Tt$.
Since confusion is unlikely we use the same notation for the respective height functions in $\Tt$ and $\Tqt$.

\begin{defi}[4.5 of \cite{WW11}]\label{qhtdef} We define a function $\htt:\Tqt\to\om$ as follows:
\[\htt(\al):=\left\{\begin{array}{ll}
   m+1 & {\small \mbox{ if } m=\max\set{k}{\mbox{there is a subterm of $\al$ of shape $\thtqk(\eta)$}}}\\
   0   & {\small \mbox{ if such $m$ does not exist.}}
   \end{array}\right.\]
\end{defi}

The above definition applies to $\Tqomie$ as well, considering parameters $<\Omie$ not as $\thtq$-terms. 
$\htomi$ is weakly increasing on $\Tqomi\cap\Omie$, and for $\al\in\Tqomi\cap\Omie$ we have 
\[\htomi(\al)=\left\{\begin{array}{ll}
      0 & \mbox{ if }\al<\Omi\\
      \min\set{n\ge i+1}{\al<\thtqi(\thtqn(0))} & \mbox{ otherwise,}\end{array}\right.\] 
cf.\ Lemma 4.6 of \cite{WW11}.

\subsection{Localization}\label{localizationsec}
The notion of \emph{localization} was introduced in Section 4 of \cite{W07a}, a summary is also given in Subsection 2.2 of \cite{W21}.
We are going to define this notion in the context of stepwise defined notations. 
We will state a version of the definition of localization (cf.\ Definition 4.6 of \cite{W07a} and Definition 2.32 of \cite{W21}) for segments 
$[\Omi,\Omie)$, $i<\om$, called $\Omi$-localization, which is completely analogous to $\tau$-localization for the lowest level
$[\Om_0,\Om_1)$ in the setting $\Om_0=\tau$. The corresponding notion of localization for systems $\Tqt$ is completely analogous,
as pointed out in Section 5 of \cite{WW11}, as are all ordinal arithmetical tools introduced in \cite{W07a} and Section 5 of \cite{CWa}.

By $\Piset(\al)$ we denote the set of all $\thti$-subterms of a term $\al\in\Tt$ that are not in the 
scope of a $\thtj$-subterm of $\al$ for any $j<i$, cf.\ 2.23-2.27 of \cite{W21}. Thus, for any $\thti$-term $\al=\thti(\xi)$, where
$\xi\in\Tt\cap\Omiz$, we have $\xi^\stari=\max(\Piset(\xi))=\max(\Piset(\al)\setminus\{\al\})$. Note that the function ${}^\stari$ is 
constantly $0$ for arguments strictly below $\Omi$ and weakly increasing on $\Tomi\cap[\Omi,\Omie)$ as $\ze^\stari$ returns the 
largest additive principal number $\le\ze$ for such $\ze$.

For $\al=\thtt(\De+\eta)\in\Tt$ we set 
\[\alplus:=\thtt(\De+\eta+1).\]
The crucial property of the interval $(\al,\alplus)$ is given by the following 
\begin{lem}\label{alalpluslem}[2.31 of \cite{W21}, 4.5 of \cite{W07a}]
For $\al=\thti(\De+\eta)\in\Tt$ and $\thti(\Ga+\rho)\in(\al,\alplus)$ we have $\Ga<\De$.
\end{lem}

\begin{defi}\label{localizationdef}
Let $\al=\thti(\De+\eta)\in\Tt$ where $i<\om$.
We define a finite sequence of ordinals as follows:
Set $\al_0:=\Omi$. Suppose $\al_n$ is already defined and $\al_n<\al$. Let
$\al_{n+1}:=\thti(\xi)\in\Piset(\al)\setminus(\al_n+1)$ where $\xi$ is maximal.
This yields a finite sequence $\Omi=\al_0<\ldots<\al_n=\al$ for some $n<\om$ which we call the
\boldmath{\bf $\Omi$-localization}\unboldmath\ of $\al$.\index{localization!$\Omi$-localization of $\al$}  
\end{defi} 

\emph{Example:} As a simple example consider the ordinal $\al$ denoted by $\thtnod(\thte(0)+\tau)$ where $\tau:=\thtnod(\thte(\thtz(0)))$
is the Bachmann-Howard ordinal. The ($\Om_0$-)localization of $\al$ is $(\tau,\al)$. $\al$ is the least epsilon number strictly above $\tau$.
In order to describe where the successor-epsilon number $\al$, the fixed point level of which is $\Om_1=\thte(0)$, is located, it is essential 
to specify the greatest ordinal below $\al$ of higher fixed point level, namely $\tau$, the fixed point level of which is $\thte(\Om_2)=\epsom$.

Note again that we could as well regard $\al$ in the above definition (and in the lemmas below) as a term of $\T^{\Omj}$ for any
$j\le i$. For convenience we also include the corresponding reformulation of lemmas stating key
properties and structural importance of localization.
\begin{lem}[cf.\ 2.33 of \cite{W21} and 4.7, 4.8, 4.9 of \cite{W07a}]\label{localipic}
Let $\al=\thti(\De+\eta)\in\Tt$, $\al>\Omi$, and denote the $\Omi$-localization of $\al$ 
by $(\Omi=\alnod,\ldots,\aln=\al)$ where $\al_j=\thti(\De_j+\eta_j)$ for $j=1,\ldots,n$. Then 
\begin{enumerate} 
\item For $j<n$ and any $\be=\thti(\Ga+\rho)\in(\al_j,\al_{j+1})$ we have 
$\Ga+\rho<\De_{j+1}+\eta_{j+1}$.
\item $(\De_j)_{1\le j\le n}$ forms a strictly descending sequence of multiples of $\Omie$.
\item For $j<n$ the sequence $(\al_0,\ldots,\al_j)$ is the $\Omi$-localization of $\al_j$.
\end{enumerate} 
We have the following {\bf guiding picture} for localizations:
\[\Omi<\al_1<\ldots<\al_n=\al<\alplus=\al_n^+<\ldots<\al_1^+.\]
\end{lem}

Note that for $\al$ as in the above lemma, the $\Omi$-localization of $\alplus$ is $(\alnod,\ldots,\almmin,\alplus)$,
provided that $\al>\Omi$, while $\Om_i^+=\thti(1)$ has $\Omi$-localization $(\Omi,\Omi\cdot\om)$.

\begin{lem}[2.35 of \cite{W21}, 5.2 of \cite{CWa}]\label{loclexordlem} 
Let $\al,\be\in\Tt\cap(\Omi,\Omie)\cap\Hz$ with $\al<\be$, $i<\om$. 
For the $\Omi$-localizations 
$\Omi=\al_0,\ldots,\al_n=\al$ and $\Omi=\be_0,\ldots,\be_m=\be$ we have
\[\alvec :=(\al_1,\ldots,\al_n)\klex(\be_1,\ldots,\be_m)=:\bevec .\]
\end{lem}

\begin{lem}[2.36 of \cite{W21}, 5.3 of \cite{CWa}]\label{subtermloclem}
Let $\al=\thti(\xi)\in\Tt$ with $\Omi$-localization $\alvec :=(\al_0,\ldots,\al_n)$
and $\be\in\Piset(\xi)$.
If there is $\al_j\in\Piset(\be)$ where $0\le j\le n$ then $(\al_0,\ldots,\al_j)$ is an initial
sequence of the $\Omi$-localization of any $\ga\in\Tt\cap[\be,\al]$.
\end{lem}

\begin{rmk}\label{noteworthyrmk} Noteworthy consequences of Lemma \ref{loclexordlem} are the following:
\begin{enumerate}
\item For any $\be\in(\al,\alplus)\cap\Hz$ where $\al,\be$ are $\thti$-terms with $\Omi$-localizations $\alvec,\bevec$, it follows that 
$\alvec$ is a proper initial sequence of $\bevec$, cf.\ also Lemma \ref{alalpluslem}.
\item For $\al=\thti(\xi)$ with $\Omi$-localization
$\Omi=\alnod,\ldots,\aln=\al$ the sequence $\alnod,\ldots,\alnmin$ is an initial sequence
of the $\Omi$-localization of $\xi^\stari$, provided that $\xi^\stari\ge\Omi$, since in that case 
$\alnmin\le\xi^\stari<\al$ as $\xi^\stari$ then indeed is the largest proper $\thti$-subterm of $\al$. 
Moreover, starting from the $\Omi$-localization of $\xi^\stari$, $\alnod,\ldots,\alnmin$ is obtained by erasing those ordinals that 
have fixed point level less than or equal to the fixed point level of $\al$.
\end{enumerate}
\end{rmk}

\subsection{The Order Isomorphism of \boldmath$\Tqt$ and $\Tt$\unboldmath}\label{orderisosubsec}

Having clarified the notion of localization, we are in the position to state the order isomorphism given in Section 3 of \cite{WW11} in a 
concise form, cf.\ Definitions 3.2 and 3.3, Lemma 3.5, and Corollary 3.7 of \cite{WW11}, as a preparation for Section \ref{fundseqTtbarsec}.  

\begin{defi}\label{itmdefi}
Simultaneously for $m<\om$ and recursively along the subterm relation we define \[\itm:\thetame\to\domthtqm.\] 
Let $\al=\De+\eta<\thetame$ where $\eta<\Omme\mid\De$.
    \[\itm(\al):=\Xipr+\Depr+\eta,\]
where $\Xipr$ and $\Depr$ are defined as follows:
\begin{enumerate}
\item If $\De=0$ we set $\Xipr:=0$ and $\Depr:=0$.

\item In case of $\De>0$ let $\De=_\ANF\De_1+\ldots+\De_k$. If $\De_1=\Omme$ we put $\Xi:=0$, otherwise
let $\Xi$ be such that $\thtme(\Xi+\rho)$, where $\rho<\Ommz\mid\Xi$, is the first element of the $\Omme$-localization of $\De_1$. We then define
\[\Xipr:=\itme(\Xi).\]
Letting $\xi_i$ be such that $\De_i=\thtme(\xi_i)$ for $i=1,\ldots,k$ we set
 \[\Depri:=\thtqme(\itme(\xi_i))\]
and define
    \[\Depr:=\left\{\begin{array}{ll}
       \Depre+\ldots+\Deprk&\mbox{if }\;\Xi=0\mbox{ or }\thtme(\Xi)<\De_1\\
       \Deprz+\ldots+\Deprk&\mbox{otherwise.}
    \end{array}\right.\]
\end{enumerate}
\end{defi}
\emph{Example:} We have $\thtnod(\thte(\thtz(0)+1))=\thtqnod(\thtqz(0)+\thtqe(\thtqz(0)+1))$.

\begin{defi}\label{rtmdefi}
Recursively along the subterm relation and for a given term descending in $m$ we define \[\rtm:\domthtqm\to\thetame.\]
Let $\al=\Xi+\De+\eta\in\domthtqm$ where $\eta<\Omme\mid\De<\Ommz\mid\Xi$. We first define the auxiliary term $\Depr$ as follows.
\begin{enumerate}
\item If $\De=0$ we put $\Depr:=0$.
\item Otherwise let $\De=_\ANF\De_1+\ldots+\De_k$, pick $\xi_i$ for $i=1,\ldots,k$ such that $\De_i=\thtqme(\xi_i)$, put
      \[\Depri:=\thtme(\rtme(\xi_i)),\]
and set \[\Depr:=\Depre+\ldots+\Deprk.\]
\end{enumerate}
Noting that $\Xi\in\domthtqme$, we then define
      \[\rtm(\al):=\left\{\begin{array}{ll}
        \Depr+\eta&\mbox{if }\;\Xi=0\\
        \thtme(\rtme(\Xi))+\Depr+\eta&\mbox{otherwise.}
      \end{array}\right.\]
\end{defi}

\begin{rmk} Note that for $\al=\Xi\ge\Ommz$ we need a subsidiary recursion in $\htomme(\al)-m$.
\end{rmk}

\begin{theo}[cf.\ 3.7 of \cite{WW11}] \label{isomtheo} For $m<\om$ the domain transformation functions $\itm$ and $\rtm$ are strictly increasing and
are inverses of one another. We have $\rtm\circ\itm=\id\restriction{\thetame}$ and $\itm\circ\rtm=\id\restriction{\domthtqm}$,  in particular
          \[\domthtqm=\Image(\itm).\]
We have $\thtm=\thtqm\circ\itm$ and $\thtqm=\thtm\circ\rtm$.
\end{theo}

\begin{defi}[3.9 of \cite{WW11}]\label{isomdefi}
Let $\ft:\Tt\to\Tqt$ denote the order-isomorphism of $\Tt$ onto $\Tqt$, that is
\[\ft(\al):=\left\{\begin{array}{ll}
\al&\mbox{if }\;\al<\tau\\
\ft(\xi)+\ft(\eta)&\mbox{if }\;\al=_\NF\xi+\eta>\tau\\
\thtqm(\itm(\xi))&\mbox{if }\;\al=\thtm(\xi).
\end{array}\right.
\]
The inverse mapping $\gt:=(\ft)^{-1}$ is characterized similarly using the mapping $\rtm$ instead of $\itm$.
\end{defi}

\begin{lem}\label{fixpcondtransformlem}
For any term of a form $\Xi+\De+\eta\in\domthtqi$ such that $\eta<\Omie\mid\De<\Omiz\mid\Xi$, 
by definition $\rti(\Xi+\De+\eta)$ is of a form $\Ga+\eta$ such that
\[\Ga=\rti(\Xi+\De)=\left\{\begin{array}{ll}
  \Depr&\mbox{ if }\Xi=0\\
  \thtie(\rtie(\Xi))+\Depr&\mbox{ otherwise,}\end{array}\right.\]
with $\Depr$ defined as in Definition \ref{rtmdefi}. We have \[F_i(\Xi+\De,\eta)\quad\aeq\quad F_i(\Ga,\eta).\]
\end{lem}
{\bf Proof.} Note first that $\eta=\gt(\eta)$. In order to prove the equivalence we apply Propositions \ref{fixpcondqprop} and \ref{fixpcondprop}, respectively. 
We have \[\eta>(\Xi+\De)^\stari\quad\aeq\quad\eta>\Ga^\stari,\]
since domain transformation $\rti$ does not change the $\mbox{}^\stari$-values. Furthermore, $\eta$ is of a shape $\eta=\thtqi(\Si+\rho)$ such that $\Si>\Xi+\De$ 
if and only if  $\gt(\eta)$ is of a shape $\gt(\eta)=\thti(\Sipr+\rho)$ such that $\Sipr>\Ga$, where $\Sipr=\rti(\Si)$, 
since $\rti$ is strictly increasing.
\qed

\begin{lem}\label{localizationtransformlem}
Let $\al=\thti(\De+\eta)\in\Tt$ with $\Omi$-localization $(\al_0,\ldots,\aln)$. Then the $\Omi$-localization of $\ft(\al)$ is \[(\ft(\al_0),\ldots,\ft(\aln)).\]
\end{lem}
{\bf Proof.} This follows directly from the fact that the mappings $\itm$ do not change the ordinal value of $\thtm$-subterms but merely transform their
arguments again using the strictly increasing mapping $\itm$. Therefore the \emph{set} of \emph{values} of additive principal subterms in the interval 
$[\Omi,\Omie)$ of $\al$ and $\ft(\al)$, respectively, is identical (with possibly differing number of duplicate occurrences). The mapping via $\itm$ of the 
arguments of these $\thtm$-subterms to the corresponding $\thtqm$-subterms of $\ft(\al)$ is strictly increasing, which implies the claim.
\qed

\subsection{Systems of fundamental sequences}
We now fix concrete instances of well-known notions of systems of fundamental sequences and Bachmann systems and refer
the reader to the seminal paper by Buchholz, Cichon, and Weiermann \cite{BCW94} for more background and an abstract 
approach.

\begin{defi}\label{fundseqdefi} We define the notion of \emph{(Cantorian) system of fundamental sequences} on a
set $S$ of ordinals. A mapping $\cdot\{\cdot\}: S\times\N\to \sup S$ is called a 
\emph{system of fundamental sequences}
if the following conditions hold:
\begin{enumerate}
\item $0\{n\}=0$ if $0\in S$.
\item $\ze\{n\}=\zepr$ if $\ze=\zepr+1$ and $\ze\in S$.
\item For $\la\in S\cap\Lim$, the sequence $(\la\{n\})_{n\in\N}\subseteq \sup S$ is strictly increasing 
such that $\la=\sup_{n\in\N}\la\{n\}$.
\end{enumerate}
A system $\cdot\{\cdot\}: S\times\N\to \sup S$ of fundamental sequences is called \emph{Cantorian} 
if the following conditions hold:
\begin{enumerate}
\item[4.] For $\al=_\NF\xi+\eta\in S$ we have $\eta\in S$ and $\al\{n\}=\xi+\eta\{n\}$ for $n<\om$. 
\item[5.] For $\al=\om^\ze\in S$ such that $\ze=\zepr+1$ we have $\al\{n\}=\om^\zepr\cdot(n+1)$.
\end{enumerate}
A system $\cdot\{\cdot\}: S\times\N\to \sup S$ of fundamental sequences is said to have the 
\emph{Bachmann property} (see Schmidt \cite{Schm76}), if for all $\al,\be\in S$ and $n<\om$ we have $\al\{n\}\le\be\{0\}$ 
whenever $\al\{n\}<\be<\al$. Such systems are called \emph{Bachmann systems}.
\end{defi}

For convenience we introduce the notion of \emph{base system} for the set of parameters. The intention is that we want to allow
$\al[0]=0$ only for $\al=0,1,\Omi$ where $i\in(0,\om)$.

\begin{defi}\label{basesysdefi}
For $\tau\in\Ezone\cap\aleph_1$ ($\tau$ needs to be countable) suppose $\cdot\{\cdot\}:(\tau+1)\times\N\to\tau$ 
is a Cantorian system of fundamental sequences with the following additional property: 
\begin{enumerate}
\item[6.] For $\al\in (\tau+1)\cap\Hz^{>1}$, i.e.\ for additive principal limits $\al\le\tau$, we have $\al\{0\}\in\Hz$.
\end{enumerate}
We call such a system a \emph{base system}.
\end{defi}

\subsection{Organization of the article}

In Section \ref{fundseqTtsec} we will first extend a mapping as in the above Definition \ref{basesysdefi} to a mapping 
$\cdot\{\cdot\}:\Ttcirc\times\N\to\Ttcirc$ where $\Ttcirc$ is the maximal subset of $\Tt$ such that $(\Ttcirc,\cdot\{\cdot\})$ is a Cantorian system of 
fundamental sequences. It will be shown that the set \[\Ttcirc:=\{\al\in\Tt \mid \domf(\al)=0\},\]
where the \emph{domain indicator function} $\domf$ is specified in Definition \ref{dominddefi}, yields a system of fundamental sequences with the 
desired properties. Note that for simplicity we have suppressed a superscript notation $\domf^\tau$ and written simply $\domf$
instead, assuming that the corresponding $\tau$ is understood from the context. In Section \ref{BachmannTtsec} we prove Bachmann property 
for $\Ttcirc$, and in Section \ref{normsec} we will introduce norms and define a Hardy hierarchy for $\Ttcirc$. 
This covers the results that were given without proofs in \cite{W24} due to page limit.

In Section \ref{fundseqTtbarsec} we will carry out the corresponding procedure for $\Tqt$ by providing a system of fundamental sequences, 
proving Bachmann property, and introducing norms and Hardy hierarchy.

\section{Fundamental sequences for \boldmath $\Tt$\unboldmath}\label{fundseqTtsec}

Before we can define the domain indicator function $\domf$, we need to define a crucial function already introduced as essential tool in the analysis
of pure patterns of order $2$ in \cite{CWc}. It can be seen as a characteristic function for uncountable ``moduli'' in ordinal terms and enables the 
definition of a function $\domf$ in Definition \ref{dominddefi} such that the cofinality of any limit ordinal $\al\in\Tt$ is $\aleph_{\domf(\al)}$ (for countable
$\tau$ and $\Omi=\aleph_i$ for $0<i<\om$).
\begin{defi}[cf.\ Definition 3.3 of \cite{W21}]\label{chidefi}
We define a characteristic function $\chiomie:\Tomie\to\{0,1\}$, where $i<\om$, by recursion on the build-up of $\Tomie$:
\begin{enumerate}
\item $\chiomie(\al):=\left\{\begin{array}{cl} 
                           0&\mbox{ if } \al<\Omie\\
                           1&\mbox{ if } \al=\Omie,
           \end{array}\right.$ 
\item $\chiomie(\al):=\chiomie(\eta)$ if $\al=_\NF\xi+\eta$,
\item $\chiomie(\al):=\left\{\begin{array}{cl}
          \chiomie(\De) & \mbox{ if } \eta\not\in\Lim\mbox{ or }F_j(\De,\eta)\\[2mm]
          \chiomie(\eta)& \mbox{ otherwise,}
      \end{array}\right.$\\[2mm]
if $\al=\thtj(\De+\eta)>\Omie$ and hence $j\ge i+1$.
\end{enumerate}
\end{defi}
\begin{rmk}\label{bracketsupportrmk} 
Recall the first paragraphs of Subsection \ref{Ttsubsec} and note that the mapping $\chiomie$ remains well defined if 
restricted to notations from $\T^{\Omj}$ for $j\le i$. 
For $j=0$ and suitable $\Om_0=\tau\in\Ezone$ this holds with $\T^0=\Tt$.  
We will tacitly make use of the canonical embeddings from $\T^\Omi$ into $\T^\Omj$ for $i<j$.
\end{rmk}

\begin{defi}\label{dominddefi}
We define a domain indicator function $\domf:\Tt\to\om$ recursively in the term build-up.
\begin{enumerate}
\item $\domf(\al):=0$ if $\al<\tau$,
\item $\domf(\al):=\domf(\eta)$ if $\al=_\NF\xi+\eta>\tau$,
\item for $\al=\thti(\De+\eta)$,
\begin{enumerate}
\item[3.1.] $\domf(\al):=\domf(\eta)$, if $\eta\in\Lim$ and $F_i(\De,\eta)$ does not hold,
\item[3.2.] if $\eta\not\in\Lim$ or $F_i(\De,\eta)$:
\begin{enumerate}
\item[3.2.1.] $\domf(\al):=\left\{\begin{array}{cl} i & \mbox{ if }\eta=0\\0 & \mbox{ otherwise}
       \end{array}\right\}$ in case of $\De=0$,\\[2mm] 
\item[3.2.2.] $\domf(\al):=0$ in case of $\chiomie(\De)=1$,
\item[3.2.3.] $\domf(\al):=\domf(\De)$ otherwise.
\end{enumerate}
\end{enumerate}
\end{enumerate}
\end{defi}

The following lemma shows the partitioning of $\Tt$ into terms of equal cofinality, using the just introduced auxiliary functions 
$\chi$ and $\domf$.
\begin{lem}\label{partitioninglem}
$\Tt$ is \emph{partitioned} into the union of disjoint sets
\[\Tt=\Ttcirc\;\dot{\cup}\;\sum_{i<\om}\{\al\in\Tt\mid\chiomie(\al)=1\}.\]
\end{lem}
{\bf Proof.} 
A straightforward induction on the build-up of terms shows for $i<k<\om$ that 
\[\chiomie(\al)+\chiomke(\al)<2\]
for all $\al\in\Tomie$,
with the canonical embedding $\Tomie\subseteq\Tomke$. 
Defining \[M_0:=\{\al\in\Tt\mid\chiomie(\al)=0\mbox{ for all }i<\om\}\quad\mbox{and }\quad 
                M_{i+1}:=\{\al\in\Tt\mid\chiomie(\al)=1\}\mbox{ for }i<\om,\]
we see that the sets $(M_i)_{i<\om}$ are pairwise disjoint.

In order to prove the lemma, we are going to show the more informative claim that  
\begin{equation}
\domfinv(0)=M_0
\quad\mbox{ and }\quad
\domfinv(i+1)=M_{i+1}\mbox{ for }i<\om.
\end{equation} 
Since $\domf$ is a well-defined function on the entire domain $\Tt$, we then obtain the desired partitioning result.
We proceed by induction on the build-up of terms in $\Tt$ along the definition of $\domf$.
If $\al\le\tau$, we have $\domf(\al)=0$ and $\chiomie(\al)=0$ for all $i<\om$, hence $\al\in M_0$.
If $\al=_\NF\xi+\eta>\tau$, we have $\domf(\al)=\domf(\eta)$, and $\chiomie(\al)=\chiomie(\eta)$ for any $i<\om$, so that
the claim follows from the i.h. Suppose $\al=\thti(\De+\eta)>\tau$. 
\\[2mm]
{\bf Case 1:} $\eta\in\Lim$ and $F_i(\De,\eta)$ does not hold. By definition we have both $\domf(\al)=\domf(\eta)$ and
$\chiomje(\al)=\chiomje(\eta)$ for all $j<\om$, as is easily checked considering the cases $i>j$ and $i\le j$. The i.h.\ thus
applies to $\eta$, and the claim follows. 
\\[2mm]
{\bf Case 2:} $\eta\not\in\Lim$ or $F_i(\De,\eta)$. 
\\[2mm]
{\bf Subcase 2.1:} $\De=0$. 
\\[2mm]
{\bf 2.1.1:} $\eta=0$. We then have $\al=\Omi$, and clearly $\Omi\in\domfinv(i)$. We already dealt with the case $i=0$ in this
context, and for $j$ such that $i=j+1$ we have $\chiomje(\Omi)=1$, that is, $\al\in M_i$.  
\\[2mm]
{\bf 2.1.2:} $\eta>0$. So, $\al=\thti(\eta)\in\domfinv(0)=\Ttcirc$, and for $j<i$ we have $\chiomje(\al)=\chiomje(0)=0$, while 
$\chiomje(\al)=0$ for $i\le j$ since $\al<\Omje$. This shows that $\al\in M_0$.
\\[2mm]
{\bf Subcase 2.2:} $\chiomie(\De)=1$. Here we again have $\al\in\domfinv(0)=\Ttcirc$. For any $j\ge i$ we have $\chiomje(\al)=0$
since $\al<\Omje$, and for $j<i$ we have $\chiomje(\al)=\chiomje(\De)=0$ since $j\not= i$ using the disjointness of the $M$-sets.
\\[2mm]
{\bf Subcase 2.3:} Otherwise, i.e.\ $\De>0$ such that $\chiomie(\De)=0$. Then by definition $\domf(\al)=\domf(\De)$, as well as
$\chiomje(\al)=\chiomje(\De)$ for all $j<\om$, checking the cases $i>j$, $i=j$, and $i<j$. Now the i.h.\ applies to $\De$, 
the claim follows.
\qed

We are now prepared to define for each $\al\in\Tt$, $\tau<\al\in\Lim$, a strictly increasing sequence $\al[\cdot]$ with supremum $\al$,
given a base system for all parameters $\le\tau$ according to Definition \ref{basesysdefi}. The extended system $\cdot[\cdot]$ 
inherits the properties $0[n]=0$ and $(\al+1)[\ze]=\al$ for all $\al\in\Tt$ as well as the property that for $\al\in\Hz^{>1}$
we have $\al[0]\in\Hz$, \emph{provided that} $\al\not=\Omi$ for all $i\in(0,\om)$. 

\begin{defi}\label{bsystemdefi}
For $\tau\in\Ezone\cap\aleph_1$ let $\cdot\{\cdot\}:(\tau+1)\times\N\to\tau$ be a base system according to Definition \ref{basesysdefi}.  
Fix the canonical assignment $\Om_0:=\tau$ and $\Omie:=\aleph_{i+1}$ for $i<\om$.
Let $\al\in\Tt$. By recursion on the build-up of $\al$ we define the function $\al[\cdot]:\aleph_d\to\T^{\Om_d}$ where $d:=\domf(\al)$.
Let $\ze$ range over $\aleph_d$.
\begin{enumerate}
\item $\al[\ze]:=\al\{\ze\}$ if $\al\le\tau$.
\item $\al[\ze]:=\xi+\eta[\ze]$ if $\al=_\NF\xi+\eta>\tau$.
\item For $\al=\thti(\De+\eta)$ where $i<\om$, $\tht_0=\thtt$, note that $d\le i$, and denote the $\Omi$-localization of 
$\al$ by $\Omi=\al_0,\ldots,\alm=\al$. We define a support term $\ual$ by
\[\ual:=\left\{\begin{array}{cl}
            \almmin & \mbox{ if either } F_i(\De,\eta)\mbox{, or: } \eta=0\mbox{ and }\De[0]^{\star_i}<\almmin=\De^{\star_i}
            \mbox{ where }m>1\\[2mm]
            \thti(\De+\etapr) & \mbox{ if } \eta=\etapr+1\\[2mm]
            0 & \mbox{ otherwise.}
     \end{array}\right.\] 
For $\al>\tau$ the definition then proceeds as follows.
\begin{enumerate}
\item[3.1.] If $\eta\in\Lim$ and $\neg F_i(\De,\eta)$, that is, $\eta\in\Lim\cap\sup_{\si<\eta}\thti(\De+\si)$, we have $d=d(\eta)$ 
and define 
   \[\al[\ze]:=\thti(\De+\eta[\ze]).\]
\item[3.2.] If otherwise $\eta\not\in\Lim$ or $F_i(\De,\eta)$, we distinguish between the following 3 subcases.
\begin{enumerate}
\item[3.2.1.] If $\De=0$, define \[\al[\ze]:=\left\{\begin{array}{cl} \ual\cdot(\ze+1) & \mbox{ if }\eta>0\mbox{ (and hence $d=0$)}\\
            \ze & \mbox{ otherwise.}
       \end{array}\right.\]
\item[3.2.2.] $\chiomie(\De)=1$. This implies that $d=0$, and we define recursively in $n<\om$
\[\al[0]:=\thti(\De[\ual])\quad\mbox{ and }\quad\al[n+1]:=\thti(\De[\al[n]]).\]
\item[3.2.3.] Otherwise. Then $d=d(\De)$ and \[\al[\ze]:=\thti(\De[\ze]+\ual).\]
\end{enumerate}
\end{enumerate} 
\end{enumerate}
Recalling that \[\Ttcirc=\domfinv(0),\] we call the system $(\Ttcirc,\cdot\{\cdot\})$ 
(more sloppily also simply the entire mapping $\cdot\{\cdot\}$),  where the mapping $\cdot\{\cdot\}$ is simply the restriction
of $\cdot[\cdot]$ to $\Ttcirc$, 
a \emph{Buchholz system over $\tau$}. Note that this system is determined uniquely modulo the choice of 
$\cdot\{\cdot\}:(\tau+1)\times\N\to\tau$, which in turn is trivially determined if $\tau=1$.   
\end{defi}

\begin{rmk}\label{alcdotrmk}
\begin{enumerate}
\item Any $\al\in\Tt$ gives rise to a unique function \[\al[\cdot]:\aleph_{\domf(\al)}\to1+\al\] 
such that $\sup_{\ze<\aleph_{\domf(\al)}}\al[\ze]=\sup\al$ and which is strictly increasing for $\al\in\Lim$ and constant otherwise, 
as will be confirmed shortly. 
\item\label{star_makes_sense_rmk} Recall the first paragraphs of Subsection \ref{Ttsubsec} and Remark \ref{bracketsupportrmk}, and note that 
restricting parameters $\ze$ to $\Tomk$ for some $k<\domf(\al)$ yields 
a function $\al[\cdot]:\Tomk\cap\aleph_d\to\Tomk$, so that the application of $\mbox{ }^\stari$ for $i\ge k$ becomes meaningful. 
\end{enumerate}
\end{rmk}

\begin{rmk}\label{examplesrmk}
We list a few instructive examples of fundamental sequences, writing $\Om$ for $\Om_1$ and 
$\mbox{ }^\star$ for $\mbox{ }^\starnod$.
\begin{enumerate}
\item $\al:=\phi_\om(0)=\thtnod(\Om\cdot\om)=\thtnod(\thte(1))$, for which $\al\{n\}=\phi_{n+1}(0)$.
\item $\be:=\thtnod(\Om^\om)=\thtnod(\thte(\thte(1)))$, the so-called small Veblen number, for which
$\be\{n\}=\phi(1,\vec{0}_{n+2})$ where $\vec{0}_k$ is a $0$-vector of length $k$, cf.\ Lemma 4.11 and subsequent remark of \cite{W06}. 
\footnote{Note that $\tht$-functions in \cite{W06} are defined 
over $+$ and $\om$-exponentiation as basic functions, which however does not change the equality at this level.}
\item $\ga:=\phi_{\Ga_0}(1)=\thtnod(\thte(\Ga_0))=\thtnod(\thte(\thtnod(\thte(\thte(0)))))$, for which 
$\ga\{n\}=\thtnod(\thte(\Ga_0\{n\})+\Ga_0)$ with $\Ga_0\{0\}=\thtnod(\thte(0))=\epsn$ and 
$\Ga_0\{n+1\}=\thtnod(\thte(\Ga_0\{n\}))$, that is, $\Ga_0\{n\}=(\thtnod\circ\thte)^{(n+1)}(0)$.
\item $\de:=\phi_{\eps_{\Ga_0+1}}(0)=\thtnod(\thte(\thtnod(\Om_1+\Ga_0)))$, for which 
$\de\{n\}=\thtnod(\De\{n\})$, as $\uDe=0$ because $\De\{0\}^\star=\Ga_0\cdot\om>\Ga_0$,
where $\De:=\thte(\thtnod(\Om_1+\Ga_0))$.
$\De\{n\}=\thte(\eps_{\Ga_0+1}\{n\})$, where $\eps_{\Ga_0+1}\{0\}=\Ga_0\cdot\om$ and
$\eps_{\Ga_0+1}\{n\}$ consists of a tower of height $n+1$ for $n>0$:
\[\eps_{\Ga_0+1}\{n\}=\Ga_0^{\rddots{\Ga_0^\om}}.\] 
\item $\eta:=\thtnod(\thte(\tht_2(\Om_1)))$. Then $\eta\{0\}$ is the Bachmann-Howard ordinal, and 
in general we have \[\eta\{n\}=(\thtnod\circ\thte\circ\tht_2)^{(n+1)}(0).\] Note that we have
$\tht_2(\Om_1)=\Om_2\cdot\Om_1$ and $\chi^{\Om_1}(\De)=1$ for $\De:=\thte(\tht_2(\Om_1))=\phi_{\Om_1}(1)$,
and for $\ze<\Om_1$ we have $\De[\ze]=\thte(\tht_2(\ze))=\phi_{\om^\ze}(\Om_1+1)$.
\end{enumerate}
\end{rmk}

The following remark contains observations regarding the support term $\ual$. Before discussing this in technical detail, we should
mention that localizing an ordinal $\al$ in terms of fixed point levels can be refined to the notion of \emph{fine-localization}, which was 
introduced in Section 5 of \cite{CWa} and takes the level of \emph{limit point thinning} into account. In the context of Definition \ref{bsystemdefi},
in particular in view of clause 3.1, we only need to refine from taking the predecessor $\almmin$ to taking $\ual$, which is a simple case of
fine-localization. Other cases of fine-localization are \emph{automatically} taken care of by clause 3.1.

\begin{rmk}\label{defaultforlargecofinalityrmk}
\begin{enumerate}
\item If $\eta=\etapr+1$ for some $\eta$ we have $\al=\ual^+$, so the next ordinal of fixed point level $\De$ is built up above
$\ual$. For example, for $\al=\thti(\Omie+1)$ the ordinal $\al[n]$ is \emph{not} $(\thti(\Omie))\cdot (n+1)$, but rather $\thti(\thti(\Omie))=\eps_{\Om_i+1}\cdot\om$ 
for $n=0$ and $\thti(\al[k])$ for $n=k+1$, thus building up iterated applications of $\thti$ starting from the support term $\ual=\thti(\Omie)=\eps_{\Om_i+1}$, approximating
the ordinal $\al=\eps_{\Om_i+2}$.

\item In case of $F_i(\De,\eta)$ (when considering an ordinal $\al=\thti(\De+\eta)$) we have $\eta=(\De+\eta)^{\star_i}=\almmin=\ual$ with $m>1$  
by (\ref{fixpcondequiv}) of Proposition \ref{fixpcondprop}, so that the next ordinal of fixed point level $\De$ is built up above the predecessor $\eta$ of $\al$ 
in its localization, which has fixed point level strictly greater than $\De$.

\item As mentioned above, in clause 3.1 of Definition \ref{bsystemdefi} of $\al[\ze]$ the support issue is automatically taken care of by $\eta[\ze]$ as we do
not have $\eta=\almmin$, so that $\si\mapsto\thti(\De+\si)$ is continuous at $\eta$, cf.\  (\ref{continuityeq}) of Proposition \ref{continuityprop}.

\item In clause 3.2.1 we either have $\eta>0$ to the effect that $\al=\ual\cdot\om$ is a successor-additive principal number that must 
be approximated by finite multiples of $\ual$. If $\eta=0$, we have $\al=\Omi$ with $i>0$, so that $\al$ is 
of uncountable cofinality and hence is approximated by parameters $\ze<\Omi$.

\item In clause 3.2.2 the ordinal $\al$ is approximated by approximation of its fixed point level $\De$ via term nesting starting 
from the support term $\ual$, which as just described takes care of the possible contribution of $\eta$, taking into account 
the collapsing that takes place when $\thti$ is applied to $\De+\eta$, where $\De$ is of cofinality $\Omie$ before collapsing.

\item In clause 3.2.3, if $\eta=0$, recall that as a consequence of Lemma \ref{loclexordlem} we have $\almmin\le\De^\stari$. 
We need the support term to be $\ual=\almmin$ if $\De[0]^{\star_i}<\almmin=\De^\stari$, since otherwise we would have $\al[\ze]<\almmin$ for all $\ze$, and the 
sequence for $\al$ would fall short. 
The need for the support term is seen by examples like $\al=\thte(\tht_2(\nu))$ where 
$\nu=\thte(\tht_2(\tht_2(\Om_1)))$, the $\Om_1$-localization of which is $(\Om_1,\nu,\al)$,
and for which we have $\thte(\tht_2(\nu[\ze]))<\nu$ for every $\ze<\Om_1$. 

If on the other hand $\eta=0=\ual$ in clause 3.2.3, it will follow from parts \ref{bracketlocalizationclaim} 
and \ref{stardomclaim} of Lemma \ref{bracketsmainlem}  that the situation $\De[0]^\stari<\almmin<\De^\stari$ can not occur, 
in other words: if $\eta=0$, $m>1$, and
$\De[0]^\stari<\De^\stari$, it follows that $\almmin\le\De[0]^\stari$. 

In summary, if we need a support term $\ual>0$, $\almmin$ suffices and is equal to $\De^\stari$.

\item Note that this \emph{minimal} choice of support term $\ual$ to guarantee convergence of the sequence to 
$\al$ is crucial for the assignment of fundamental sequences to satisfy Bachmann property. A simple example illustrates this. Consider
in clause 3.2.3 for $\al=\thti(\De)$ (assuming that $\eta=0$) the alternative definition \[\al\langle\ze\rangle:=\thti(\De[\ze]+\De^\stari),\] 
which still converges to $\al$.
Now consider, writing $\om$ instead of $\thtnod(\thtnod(0))$, $\al:=\thtnod(\thte(\om))=\phi_{\om^\om}(0)$. We then have  
$\al[n]=\thtnod(\thte(n+1))=\phi_{\om^{n+1}}(0)$ and $\al\langle n\rangle=\thtnod(\thte(n+1)+\om)=\phi_{\om^{n+1}}(\om)$. For
$\be:=\al[n+1]=\phi_{\om^{n+2}}(0)$ we have $\be[0]=\phi_{\om^{n+1}}(0)$ and $\be\langle 0\rangle=\phi_{\om^{n+1}}(1)$, and we observe that
$\al[n],\al\langle n\rangle<\be<\al$ and $\al[n]=\be[0]<\be\langle0\rangle<\al\langle n\rangle$, demonstrating that the alternative definition using the 
$\mbox{ }^\star$-operation instead of localization (only where a support term is necessary) does not satisfy Bachmann property.
\end{enumerate}
\end{rmk}

Due to frequent occurrence in upcoming proofs we consider a few special cases when dealing with ordinals $\al$ and $\be$ such that $\al[\ze]<\be<\al$.
\begin{lem}\label{auxlem} Let $\al,\be\in\Tt$ be of the form $\al=\thti(\De+\eta)$ and $\be=\thti(\Ga+\rho)$ such that $\al[\ze]<\be<\al$ for some
$\ze<\aleph_{\domf(\al)}$.
\begin{enumerate}
\item If either $\eta\not\in\Lim$ or $F_i(\De,\eta)$ holds, and $\Ga=\De$, then $\be\le\ual$.
\item Suppose that either $\rho\not\in\Lim$ or $F_i(\Ga,\rho)$ holds, and that we are in one of the following two scenarios:
\begin{enumerate}
\item $\eta\in\Lim$, $F_i(\De,\eta)$ does not hold, and $\Ga=\De$, or
\item either $\eta\not\in\Lim$ or $F_i(\De,\eta)$ holds, $\chiomie(\De)=0$, and $\Ga=\De[\ze]$. 
\end{enumerate} 
Then $\al[\ze]\le\ube$.
\end{enumerate}
\end{lem}
{\bf Proof.} For part 1, note that the assumptions imply $\rho<\eta$. We obtain
\[\be=\thti(\De+\rho)\left\{\begin{array}{ll}
\le\thti(\De+\etapr) & \mbox{ if }\eta=\etapr+1\\[2mm]
<\sup_{\si<\eta}\thti(\De+\si)=\eta & \mbox{ if }F_i(\De,\eta)\mbox{ holds}\end{array}\right\}=\ual.\] 
\\[2mm]
Regarding part 2, in scenario (a) we have $\domf(\eta)=\domf(\al)$ and $\eta[\ze]<\rho$, while in scenario (b) we have $\domf(\De)=\domf(\al)$ and $\ual<\rho$.
Thus 
\[\al[\ze]=\left\{\begin{array}{ll}
\thti(\De+\eta[\ze]) & \mbox{ in scenario (a) }\\[2mm]
\thti(\De[\ze]+\ual) & \mbox{ in scenario (b) }\end{array}\right\} 
\left\{\begin{array}{ll}
\le\thti(\Ga+\rhopr) & \mbox{ if } \rho=\rhopr+1\\[2mm]
<\sup_{\si<\rho}\thti(\Ga+\si)=\rho & \mbox{ if }F_i(\Ga,\rho)\mbox{ holds}\end{array}\right\}=\ube,\] 
which concludes the proof. \qed

The following main lemma provides several properties shown by simultaneous induction crucial to verify that $\Ttcirc$ is a Cantorian system of fundamental 
sequences enjoying Bachmann property. 

\begin{lem}\label{bracketsmainlem} 
In the setting of Definition \ref{bsystemdefi}, let $\al\in\Tt\setminus\tau$ be a limit number, and let $\ze$ range over $\aleph_d$ 
where $d:={\domf(\al)}$.
\begin{enumerate}
\item\label{ualcontrolclaim} For $\al=\thtj(\De+\eta)$ such that $\ual>0$ and $\chiomje(\De)=0$ we have \[\De[\ze]^\starj<\ual.\]
\item\label{fundseqclaim} The mapping $\ze\mapsto\al[\ze]$ is strictly increasing with proper supremum $\al$. 
\item\label{starcontrolclaim} The mapping $\ze\mapsto\al[\ze]^{\stark}$ is weakly increasing with upper bound $\al^{\stark}$
 for $d\le k$. 
\item\label{bracketparestimclaim} If $d=i+1>0$, we have 
\[\ze^\stark\le\al[\ze]^{\stark}\le\max\{\al^{\stark},\ze^{\stark}\}\quad\mbox{ and }\quad\al^{\stark}\le\max\{\al[\ze]^{\stark},1\}\] 
for any $k\le i$.
\item\label{bracketlocalizationclaim} For $\al=\thtj(\De+\eta)>\Omj$, denote the $\Omj$-localization of $\al$ by 
$\Omj=\alnod,\ldots,\alm=\al$. For such $\al$ we have \[\almmin\le\al[0].\]
\item\label{stardomclaim} If $\al=\thtj(\De+\eta)$ and $i$ is such that $\domf\le i\le j$ and $\al[0]^\stari<\al^\stari$
where $\al^\stari>\Omi$, then we have $\domf(\al^\stari)=\domf$, and one of the following two cases applies.
\begin{enumerate}
\item $\al^\stari$ is of the form $\nu\cdot\om$ for some additive principal $\nu\ge\Omi$, where we have $d=0$,
$\al[n]^\stari=\nu$, and $\al^\stari[n]=\nu\cdot(n+1)$ for $n<\om$. 
\item Otherwise. Then $\al^\stari[\ze]\le\al[\ze]^\stari$ for all $\ze<\aleph_\domf$, 
and there exists $\ze_0<\aleph_\domf$ such that $\al^\stari[\ze]=\al[\ze]^\stari$ for all $\ze\in(\ze_0,\aleph_\domf)$.
\end{enumerate}
\item\label{sandwichclaim} For any $\be$ such that $\al[\ze]\le\be<\al$ we have \[\be^\stark\ge\al[\ze]^\stark\] 
for all $k$ if $\ze=0$, and 
for all $\ze<\aleph_\domf$ and all $k$ such that $k+1\ge\domf$.
\end{enumerate}
\end{lem}
{\bf Proof.} We argue by induction along the build-up of terms of $\Tt$. Before proving the lemma for the main cases where
$\al$ is of a form $\thtj(\De+\eta)>\Omj$, we discuss the other cases for $\al$, provide general argumentation regarding 
part \ref{sandwichclaim}, which is based on subterm relationship, and show part \ref{ualcontrolclaim} in general on the basis of the i.h.

For $\al=\tau$ claim \ref{fundseqclaim} holds by assumption, and the other assertions are trivially fulfilled. 
For $\al=\Omie=\thtie(0)$, so that $j=d=i+1$, we have $\al[\ze]=\ze$, and recalling part \ref{star_makes_sense_rmk} of 
Remark \ref{alcdotrmk}, all claims follow immediately, noticing that for part
\ref{sandwichclaim} we either trivially have $\ze=0$, or we have $k\ge i$ by assumption. 
The case $\al=_\NF\xi+\eta$ also directly follows from
the i.h., apart from part \ref{sandwichclaim}, where we first note that the additive normal form of $\be$ must start with $\xi$, as 
$\al[\ze]=\xi+\eta[\ze]$.

Note that \emph{in general for part \ref{sandwichclaim}} the case $\al[\ze]=\be$ is trivial, so that we may suppose $\al[\ze]<\be$. 
The case $\sumend(\al)[\ze]<\Omk$ is trivial as well, so we may assume that $\sumend(\al)[\ze]\ge\Omk$, and since the situation 
$\sumend(\al)=\nu\cdot\om$ for some $\nu$, where $\sumend(\al)[n]=\nu\cdot(n+1)$, is trivial as it implies that the 
additive normal form of $\be$ also begins with $\al[n]$, we may assume that $\sumend(\al)[\ze]\in\Hz\setminus\Omk$. 
Without loss of generality we may therefore assume that in the case where $\al=_\NF\xi+\eta$, i.e.\ $\sumend(\al)=\eta<\al$, 
$\be=_\NF\xi+\rho$ for some additive principal $\rho<\eta$, and that $\be\in\Hz$ in the cases where $\al\in\Hz$, 
i.e.\ $\sumend(\al)=\al$. 

Returning to the verification of part \ref{sandwichclaim} in the particular case $\al=_\NF\xi+\eta$, we have $\eta[\ze]<\rho<\eta$
and $k+1\ge\domf=\domf(\eta)$ at least if $\ze>0$, so that the i.h.\ for $\eta$ applies to yield $\rho^\stark\ge\eta[\ze]^\stark$, 
which implies the desired $\be^\stark\ge\al[\ze]^\stark$.

In order to prove all parts of the lemma it now remains to consider terms of the form $\al=\thtj(\De+\eta)>\Omj$, 
and we denote the $\Omj$-localization of $\al$ by $\Omj=\alnod,\ldots,\alm=\al$. Note that $d\le j$.

{\bf Part \ref{ualcontrolclaim}} is seen to hold by inspecting the definition of $\ual$, using (\ref{fixpcondequiv}) of Proposition \ref{fixpcondprop}
and part \ref{starcontrolclaim} of the i.h.\ for $\De$, noting that $\domf(\De)\le j$ as $\chiomje(\De)=0$.
This is immediate if $\eta=\etapr+1$ for some $\etapr$ since then $\ual=\thtj(\De+\etapr)$ and $\De[\ze]^\starj\le\De^\starj<\ual$.
Similarly, if $F_j(\De,\eta)$, we have $\De[\ze]^\starj\le\De^\starj<\eta=\almmin=\ual$.
If $\eta=0$ and $\ual>0$, we have $m>1$ and $\De[0]^\starj<\De^\starj=\almmin=\ual\in\Ez$, so that by part \ref{stardomclaim} of 
the i.h.\ for $\De$ we have $\domf(\De^\starj)=\domf(\De)$,  $\De^\starj[\ze]\le\De[\ze]^\starj$ for all $\ze$, and there exists 
$\ze_0<\aleph_{\domf(\De)}$ such that 
$\De^\starj[\ze]=\De[\ze]^\starj$ for $\ze\in(\ze_0,\aleph_{\domf(\De)})$. Using the i.h.\ for $\De$ to see that $\De[\ze]^\starj$ is
weakly increasing in $\ze$ with upper bound $\De^\starj$ and the i.h.\ for $\De^\starj$ to see that  $\De^\starj[\ze]$ is strictly
increasing in $\ze$ with supremum $\De^\starj$, we obtain $\De[\ze]^\starj<\De^\starj=\ual$ for all $\ze<\domf(\De)$.

For the other parts of the lemma, we proceed by case distinction along the definition of $\al[\ze]$.
\\[2mm]
\noindent{\bf Case 1:} $\eta$ is a limit ordinal for which $F_j(\De,\eta)$ does not hold. 
We then have $d=\domf(\eta)$ and $\al[\ze]=\thtj(\De+\eta[\ze])$.

Here, {\bf parts \ref{fundseqclaim}, \ref{starcontrolclaim}, and \ref{bracketparestimclaim}} directly follow from the i.h.\ since 
the function $\si\mapsto\thtj(\De+\si)$ is strictly increasing and continuous at $\eta$, cf.\  (\ref{continuityeq}) of Proposition \ref{continuityprop}. 

For {\bf part \ref{bracketlocalizationclaim}} assume the non-trivial situation $\almmin>\Omj$. If $\almmin$ is a subterm
of $\De$ we immediately obtain $\almmin<\al[0]$. Otherwise $\almmin$ is a subterm of $\eta$. Since $\neg F_j(\De,\eta)$,
it must be a proper subterm, and hence by Lemma \ref{subtermloclem} $\almmin<\eta$ is an element of the $\Omj$-localization of $\mc(\eta)$, 
so that the i.h.\ for $\eta$ yields $\almmin\le\eta[0]<\al[0]$. 

We turn to the verification of {\bf part \ref{stardomclaim}}. Assume that $\domf\le i\le j$ for some $i$, 
and $\Omi,\al[0]^\stari<\al^\stari$. 
We first consider the case $i=j$. We then have $\al^\stari=\al$. 
Since by assumption of Case 1, $\al$ is a limit of additive principal numbers, we have $\al[\ze]^\stari=\al[\ze]=\al^\stari[\ze]$.  
Assume now that $i<j$.
Then we have $\al[0]^\stari=\max\{\De^\stari,\eta[0]^\stari\}<\eta^\stari=\al^\stari$, 
hence $\Omi,\De^\stari,\eta[0]^\stari<\eta^\stari$, so that the i.h.\ for $\sumend(\eta)$ applies. 
We may therefore assume that $\eta\in\Hz$. 
Thus $\domf(\al)=\domf(\eta)=\domf(\eta^\stari)=\domf(\al^\stari)$, and either $\eta^\stari$ and hence also $\al^\stari$ is of the
form $\nu\cdot\om$ for some $\nu\in\Hz\setminus\Omi$, so that we obtain $\al^\stari[n]=\eta^\stari[n]=\nu\cdot(n+1)$ and
$\al[n]^\stari=\max\{\De^\stari,\eta[n]^\stari\}=\nu$ for all $n<\om$, or we obtain eventual equality of $\al^\stari[\ze]$ and
$\al[\ze]^\stari$ for large enough $\ze$ from this property for $\eta^\stari[\ze]$ and $\eta[\ze]^\stari$, using that 
$\ze\mapsto\al[\ze]^\stari$ is weakly increasing with upper bound $\al^\stari$, i.e.\ claim \ref{starcontrolclaim} for $\al$,
and that by claim \ref{fundseqclaim} for $\eta^\stari$, we eventually have $\eta^\stari[\ze]>\De^\stari$ for large $\ze$.

To conclude the treatment of Case 1 we show {\bf part \ref{sandwichclaim}}. Assume that for $\be$ of the form $\thtj(\Ga+\rho)$, 
cf.\ the remarks at the beginning of the proof, either $k+1\ge\domf$ and $\al[\ze]<\be<\al$ for some $\ze$,
or $\ze=0$, $k<\om$, and $\al[0]<\be<\al$. 
If $k\ge j$, it follows immediately that $\be^\stark\ge\al[\ze]^\stark$. Suppose therefore that $k<j$.
We argue by subsidiary induction on the build-up of $\be$. The criterion for comparison of $\tht$-terms, Proposition \ref{thetacompprop}, yields
\begin{equation}\label{sandwichcompone} 
\be<\al\quad\aeq\quad\left(\Ga+\rho<\De+\eta\;\wedge\;(\Ga+\rho)^{\starj}<\al\right)\;\;\vee\;\;\be\le(\De+\eta)^\starj.
\end{equation}
Note that according to part \ref{bracketlocalizationclaim}  we have $\al_{m-1}\le\al[\ze]$, so that 
$\be\in(\al_{m-1},\al)$, hence according to part 1 of Lemma \ref{localipic} we can not have $\Ga>\De$.
By definition we have $\al[\ze]=\thtj(\De+\eta[\ze])$ with $\domf(\al)=\domf(\eta)$, and according to the equivalence for 
comparison of $\tht$-terms 
\[\al[\ze]<\be\quad\aeq\quad\left(\De+\eta[\ze]<\Ga+\rho\;\wedge\;(\De+\eta[\ze])^\starj<\be\right)\;\;\vee\;\;\al[\ze]\le(\Ga+\rho)^\starj.
\]
If $\al[\ze]\le(\Ga+\rho)^\starj$, then by the s.i.h.\ applied to $(\Ga+\rho)^\starj$ we obtain 
$\be^\stark\ge((\Ga+\rho)^\starj)^\stark\ge\al[\ze]^\stark$ and are done. Now assume that $(\Ga+\rho)^\starj<\al[\ze]$, 
so we must have $\De+\eta[\ze]<\Ga+\rho$ and $(\De+\eta[\ze])^\starj<\be$. 
But this implies that $\De=\Ga$, as we already showed that $\De<\Ga$ is not possible.
This further entails that $\be^\stark\ge\Ga^\stark=\De^\stark$ as well as $\eta[\ze]<\rho<\eta$, 
whence by the i.h.\ for $\eta$ and $\rho$ also 
$\be^\stark\ge\rho^\stark\ge\eta[\ze]^\stark$.
We thus obtain $\be^\stark\ge\al[\ze]^\stark$.\\[2mm]
\indent For the remaining cases we assume that Case 1 does not apply, so $\eta$ is either $0$, a successor $\eta=\etapr+1$, 
or a limit for which $F_j(\De,\eta)$ holds.
\\[2mm]
{\bf Case 2:} $\De=0$. Then we have $d=0$, $\eta>0$, and $\almmin\le\al[n]=\ual\cdot(n+1)$, 
as we already dealt with the case $\al=\Omj$. 
We know that $\al=\ual\cdot\om$. The lemma obviously holds 
in this case, which covers all successor-additive-principal numbers, see Lemma \ref{epscharlem}.
\\[2mm]
\noindent{\bf Case 3:} $\chiomje(\De)=1$. Then it follows that $d=0$, while $\domf(\De)=j+1$,
and using in this case the more familiar
variable $n<\om$ instead of $\ze<\om$ we have
\[\al[0]=\thtj(\De[\ual])\quad\mbox{ and }\quad\al[n +1]=\thtj(\De[\al[n ]]).\]
Note that if $\eta=0$, due to part \ref{bracketparestimclaim} of the i.h.\ for $\De$ the condition for $\ual>0$ is
not satisfied, so that $\eta=0$ implies $\ual=0$ in this case.
We begin with showing {\bf part 2}.
Using again the i.h.\ for $\De$, in particular part \ref{bracketparestimclaim}, we see that 
$\ual\le\De[\ual]^\starj<\thtj(\De[\ual])=\al[0]$, 
and also $\De[\ual]^\starj\le\De[\al[0]]^\starj<\al[1]$ and
$\De[\ual]<\De[\al[0]]$, hence $\al[0]<\al[1]$, and 
by induction on $n$ we obtain $\De[\al[n]]^\starj\le\De[\al[n+1]]^\starj<\al[n+2]$ 
and $\De[\al[n]]<\De[\al[n+1]]$, hence $\al[n+1]<\al[n+2]$.

We also have $\De[\ual]<\De$ and $\De[\ual]^\starj<\al$, hence $\al[0]<\al$, and inductively in $n$
in order to see that $\al[n+1]<\al$ we have
$\De[\al[n]]<\De$ and $\De[\al[n]]^\starj<\al$ by the i.h., 
thus $\al[n]<\al$ for all $n<\om$.

We have now seen that $(\al[n])_{n<\om}$ is a strictly increasing sequence with strict upper bound 
$\al$, and it remains to be shown that $\al$ actually is the supremum of this sequence. 
We therefore define
\[\xi:=\sup_{n<\om}\al[n]\]
and use Proposition \ref{thetaprop} to show that $\al\le\xi$:

1. $(\De+\eta)^\starj<\xi$, since 
$\De^\starj\le\De[\ual]^\starj+1\le\al[0]<\xi$ and
$\eta^\starj\le\ual<\al[0]<\xi$. Note that this also implies that $\almmin\le\al[0]$, i.e.\ {\bf part \ref{bracketlocalizationclaim}}.

2. Let $\Ga+\rho<\De+\eta$ such that $(\Ga+\rho)^\starj<\xi$ be given. We need to verify that 
$\thtj(\Ga+\rho)<\xi$. If $\Ga=\De$, Lemma \ref{auxlem} yields $\thtj(\Ga+\rho)\le\ual$, and clearly $\ual<\xi$.

Now consider the case $\Ga+\rho<\De$. We need to find an $n<\om$ such that $\Ga+\rho<\De[\al[n]]$.
Assume to the contrary that $\De[\al[n]]\le\Ga<\De$ for every $n<\om$. 
Then parts \ref{bracketparestimclaim} and \ref{sandwichclaim} of the i.h.\ for $\De$ yield $\Ga^\starj\ge\al[n]^\starj=\al[n]$ 
for every $n<\om$, 
which however contradicts the assumption $(\Ga+\rho)^\starj<\xi$. Thus, there exists an $n<\om$
such that $\Ga+\rho<\De[\al[n]]$ which can as well be chosen so that $(\Ga+\rho)^\starj<\al[n+1]$,
resulting in $\thtj(\Ga+\rho)<\al[n+1]<\xi$. This concludes the proof of part \ref{fundseqclaim} in this case.
Note that {\bf part \ref{starcontrolclaim}} regarding $\mbox{}^\stark$-terms for $k<j$ is easily seen from the definition of the
terms $\al[n ]$. As to {\bf part \ref{bracketparestimclaim}}, note that the assumption $d=i+1$ for some $i$ is not fulfilled.

{\bf Part \ref{stardomclaim}} trivially holds in the base case $i=j$. If $i<j$, we have $\al[0]^\stari=\De[\ual]^\stari$, and
according to part \ref{bracketparestimclaim} of the i.h.\ for $\De$ we have 
$\De^\stari\le\max\{\De[\ual]^\stari,1\}<\eta^\stari=\al^\stari$ due to the assumptions $\al^\stari>\Omi\ge 1$ and 
$\al^\stari>\al[0]^\stari$. Since $\ual^\stari\le\De[\ual]^\stari$, the assumptions imply that $\ual^\stari<\eta^\stari$, which
is impossible as $\eta^\stari>1$. Thus, if $i<j$ the assumptions of part \ref{stardomclaim} do not hold.

{\bf Part \ref{sandwichclaim}} trivially holds for $k\ge j$, so we assume that $k<j$.
Since $\domf=0$, we write $n$ instead of $\ze$.
Arguing again by side induction on the build-up of $\be$, for $\be$ of a form $\thtj(\Ga+\rho)$ we have
\[\al[n]<\be\quad\aeq\quad\left(\De[\alpr]<\Ga+\rho\;\wedge\;\De[\alpr]^\starj<\be\right)
                                           \;\;\vee\;\;\al[n]\le(\Ga+\rho)^\starj,\]
where $\alpr:=\ual$ if $n=0$ and $\alpr:=\al[n-1]$ if $n>0$.                                        
In the case $\al[n]\le(\Ga+\rho)^\starj$ we apply the s.i.h\ as in Case 1, so we may assume that 
\[(\Ga+\rho)^\starj<\al[n]\mbox{, hence }
\De[\alpr]<\Ga+\rho\mbox{ and }\De[\alpr]^\starj<\be.\]
Note that in the case $\eta=0$ we have $\ual=0$, since according to part \ref{bracketparestimclaim} we can not have
$\De[0]^{\starj}<\De^\starj=\almmin$ with $m>1$. 
In case of $\eta=\etapr+1$ for some $\etapr$ we have $\ual=\thtj(\De+\etapr)$, and if $F_j(\De,\eta)$ holds, we have
$\ual=\eta=\almmin=(\De+\eta)^\starj$ with $m>1$. In both cases where $\ual>0$ we have $\be\in(\ual,\ual^+)$ using Lemma \ref{localipic} 
to see that $\al\le\ual^+$, and the $\Omj$-localization of $\ual$ is an initial segment of the $\Omj$-localization of $\be$, cf.\ Remark \ref{noteworthyrmk}, 
hence $\be^\stark\ge\ual^\stark$, as $\ual$ then is a subterm of $\be$. 
We can not have $\Ga=\De$, since according to Lemma \ref{auxlem} this would imply $\be\le\ual$, contradicting the assumption as $\ual<\al[0]$.
Since as in Case 1, $\almmin\le\al[0]<\be<\al$ according to the assumption and part \ref{bracketlocalizationclaim},  
by part 1 of Lemma \ref{localipic} we can not have $\Ga>\De$ either, we must therefore have $\Ga<\De$.
Then $\De[0]\le\De[\alpr]<\Ga+\rho<\De$, and the i.h.\ for $\De$ yields $(\Ga+\rho)^\stark\ge\De[0]^\stark$.
We thus obtain $\be^\stark\ge\max\{\De[0]^\stark,\ual^\stark\}=\al[n]^\stark$, as follows inductively from 
part \ref{bracketparestimclaim}, since $k<j$.
\\[2mm]
{\bf Case 4:} Otherwise.
In this final case we have $d=\domf(\al)=\domf(\De)$ and $\al[\ze]=\thtj(\De[\ze]+\ual)$. 
We start out with showing {\bf part 2}.
By the i.h.\ the function $\ze\mapsto\De[\ze]+\ual$ is strictly increasing with
supremum $\De$, as $\De>\De[\ze]\ge\Omje>\ual$ since $\chiomje(\De)=0$ 
ensures that $\De[\ze]$ is a nonzero multiple of $\Omje$, and the function 
$\ze\mapsto(\De[\ze]+\ual)^{\star_j}$ is weakly increasing with upper bound $(\De+\ual)^{\star_j}$,
which in turn is strictly less than $\al$. This implies that the mapping
$\ze\mapsto\thtj(\De[\ze]+\ual)=\al[\ze]$ is strictly increasing with upper bound $\al$.
Setting \[\xi:=\sup_{\ze<\aleph_d}\al[\ze],\] we obtain $\xi\le\al$, and in order to see that also $\al\le\xi$ 
we apply Proposition \ref{thetaprop}: 

1. $(\De+\eta)^\starj<\xi$ is shown as follows. If $F_j(\De,\eta)$ holds, invoking Proposition \ref{fixpcondprop} and 
Remark \ref{noteworthyrmk} we have $\eta>\De^\starj$, $\eta=\almmin=\ual<\al[0]<\xi$. 
If $\eta=\etapr+1$ for some $\etapr$, we have $\ual=\thtj(\De+\etapr)$
and $(\De+\eta)^\starj\le\ual<\al[0]<\xi$, where $\almmin\le(\De+\eta)^\starj$ in case of $m>1$. 
We are left with cases where $\eta=0$. 
The case $\De^\starj\le\Omj$ is trivial, so assume that $\De^\starj>\Omj$.
If $\De[0]^\starj<\almmin=\De^\starj$ then we have $m>1$ and $\xi>\ual=\De^\starj$. Otherwise we have $\ual=0$.
If in this situation $\De^\starj=\almmin$, we must have $\De^\starj\le\De[0]^\starj$, and if this latter inequality holds,
we are done since $\De[0]^\starj<\al[0]<\xi$.
Suppose finally, that $\almmin,\De[0]^\starj<\De^\starj$. 
According to Lemma \ref{localipic}, $\De^\starj$ must be of a form
$\thtj(\Ga+\rho)$ such that $\Ga+\rho<\De$ and hence $\Ga+\rho<\De[\ze]$ for large enough $\ze$. 
Since $(\Ga+\rho)^\starj<\De^\starj$, we find $\ze$ large enough to both satisfy 
$\Ga+\rho<\De[\ze]$ and $(\Ga+\rho)^\starj<\De^\starj[\ze]$, so that 
by part \ref{stardomclaim} of the i.h.\ for $\De$ we also have $(\Ga+\rho)^\starj\le\De[\ze]^\starj$, and thus
$\De^\starj=\thtj(\Ga+\rho)<\thtj(\De[\ze])=\al[\ze]<\xi$, as needed.
Note further that, according to part 2 of Remark \ref{noteworthyrmk}, the $\Omj$-localization of
$\De^\starj$ has - in the situation we are in currently - proper initial sequence $(\alnod,\ldots,\almmin)$,
so that  by the i.h.\ for $\De^\starj$ (part \ref{bracketlocalizationclaim}) and $\De$ (part \ref{stardomclaim}) we have
$\almmin\le\De^\starj[0]\le\De[0]^\starj<\al[0]$.

Looking back we see that the above verification of  $(\De+\eta)^\starj<\xi$ also showed that $\almmin<\al[0]$
holds generally in Case 4, which settles {\bf part \ref{bracketlocalizationclaim}}.  

2. Now let $\Ga+\rho$ be given with the properties $\Ga+\rho<\De+\eta$ and $(\Ga+\rho)^{\star_j}<\xi$. We have
to show that $\thtj(\Ga+\rho)<\xi$ holds as well.  
We find a large enough $\ze<\aleph_d$ such that $(\Ga+\rho)^{\star_j}<\thtj(\De[\ze]+\ual)$ and in case of 
$\Ga<\De$ also $\Ga+\rho<\De[\ze]$. In this latter case we directly obtain $\thtj(\Ga+\rho)<\thtj(\De[\ze]+\ual)<\xi$,
while otherwise $\Ga=\De$, so that $\thtj(\Ga+\rho)\le\ual$ by Lemma \ref{auxlem} and clearly $\ual<\xi$.

We have now seen that $\ze\mapsto\thtj(\De[\ze]+\ual)=\al[\ze]$ is strictly increasing with supremum $\al$, which
concludes the proof of part \ref{fundseqclaim}.
Regarding the remaining claims, first suppose $k$ satisfies $d\le k<j$ in order to show {\bf part \ref{starcontrolclaim}}.
We then have 
$\al[\ze]^{\star_k}=(\De[\ze]+\ual)^{\star_k}$, and since by the i.h.\ the mapping
$\ze\mapsto\De[\ze]^{\star_k}$ is weakly increasing with upper bound $\De^{\star_k}$,
also $\ze\mapsto\al[\ze]^{\star_k}$ is weakly increasing with upper bound $\al^{\star_k}$.

Let us turn to {\bf part \ref{bracketparestimclaim}}, and assume that $k$ satisfies $0\le k<d=i+1$ for some $i<\om$. 
By the i.h.\ we then have 
\[\ze^\stark\le\De[\ze]^\stark\le\al[\ze]^{\star_k}=(\De[\ze]+\ual)^{\star_k}\le\max\{\De^{\star_k},\ze^{\star_k},\eta^{\star_k}\}
        =\max\{\al^{\star_k},\ze^{\star_k}\}.\]
Note that if $k=0$, $\tau=1$, $\eta=\etapr+1$, and $\al[\ze]^{\star_0}=0$ the third inequality can be strict, 
since $\eta^{\star_0}=1$ in such case. 
On the other hand, using the i.h.\ for $\De$ we have 
\[\al^\stark=(\De+\eta)^\stark\le\max\{\De[\ze]^\stark,1,\eta^\stark\}\le\max\{\De[\ze]^\stark,\ual^\stark,1\}=\max\{\al[\ze]^\stark,1\}.\] 

As to {\bf part \ref{stardomclaim}}, assume that $\domf\le i\le j$ for some $i$ and $\Omi, \al[0]^\stari<\al^\stari$.
Here we have $\domf=\domf(\De)$, $\De>0$, and $\chiomje(\De)=0$. 
In the base case $i=j$, we have $\al^\stari=\al$, $\al^\stari[\ze]=\al[\ze]=\al[\ze]^\stari$, and hence the claim holds.
Assume now that $i<j$. We then have
$\al[\ze]^\stari=\max\{\De[\ze]^\stari,\ual^\stari\}$ and $\al^\stari=\max\{\De^\stari,\eta^\stari\}$. 
We consider the defining cases for $\ual$: 
\begin{enumerate}
\item $F_j(\De,\eta)$. Then we have $\ual=\eta>\De^\starj$, and it follows that $\eta^\stari,\De[0]^\stari<\De^\stari=\al^\stari$. 
The i.h.\ for
$\De$ now applies, and we obtain $\domf(\al^\stari)=\domf(\De^\stari)=\domf(\De)=\domf$ and
$\al^\stari[\ze]=\De^\stari[\ze]$. Now, if $\al^\stari=\nu\cdot\om$ for some $\nu\ge\Omi$, we obtain $\al[n]^\stari=\De[n]^\stari=\nu$
and $\al^\stari[n]=\De^\stari[n]=\nu\cdot(n+1)$. For $\al^\stari$ and hence $\De^\stari$ not a successor-principal number, since 
$\al[\ze]^\stari\le\al^\stari$ is weakly increasing by part \ref{starcontrolclaim}, the eventual equality of $\al^\stari[\ze]$ and
$\al[\ze]^\stari$ follows from the i.h.\ for $\De$ and $\De^\stari$ (part \ref{fundseqclaim}), and we have 
$\al^\stari[\ze]=\De^\stari[\ze]\le\De[\ze]^\stari\le\al[\ze]^\stari$ for all $\ze$. 
\item $\eta=\etapr+1$ for some $\etapr$. Then $\ual=\thtj(\De+\etapr)$, and the assumptions $\Omi,\al[0]^\stari<\al^\stari$
can not be satisfied.
\item $\eta=0$, $m>1$, and $\De[0]^\starj<\De^\starj=\almmin=\ual$. Then 
$(\De^\starj)^\stari,\De[0]^\stari<\De^\stari=\al^\stari$, the i.h.\ for $\De$ applies, and we argue as in the case where $F_j(\De,\eta)$.
\item $\eta=\ual=0$, hence $\al[\ze]^\stari=\De[\ze]^\stari$
and $\al^\stari=\De^\stari$, so that the claim directly follows from the i.h.\ for $\De$.
\end{enumerate}
This concludes the proof of claim \ref{stardomclaim}.

Finally, we show {\bf part \ref{sandwichclaim}}. Assume that for $\be$ of the form $\thtj(\Ga+\rho)$,
either $k+1\ge\domf$ and $\al[\ze]<\be<\al$ for some $\ze$, or $\ze=0$, $k<\om$, and  $\al[0]<\be<\al$. 
Since for $k\ge j$ the claim is trivial, we assume that $k<j$.
We argue by subsidiary induction on the build-up of $\be$. 
Recall that we have $\domf(\De)=\domf$, $\al[\ze]=\thtj(\De[\ze]+\ual)$, and that
\[\al[\ze]<\be\quad\aeq\quad\left(\De[\ze]+\ual<\Ga+\rho\;\wedge\;(\De[\ze]+\ual)^\starj<\be\right)
                                           \;\;\vee\;\;\al[\ze]\le(\Ga+\rho)^\starj.\]
In the case $\al[\ze]\le(\Ga+\rho)^\starj$ the side i.h.\ applies as in Case 1, so we may assume that 
\[(\Ga+\rho)^\starj<\al[\ze]\mbox{, hence }
\De[\ze]+\ual<\Ga+\rho\mbox{ and }(\De[\ze]+\ual)^\starj<\be.\]
Recall the definition of $\ual$: 
\[\ual=\left\{\begin{array}{cl}
      \almmin  & \mbox{ if either } F_{j+1}(\De,\eta)\mbox{, or }\eta=0 \mbox{ and } \De[0]^{\starj}<\De^\starj=\almmin\mbox{ where }m>1
                            \\[2mm]
      \thtj(\De+\etapr) & \mbox{ if } \eta=\etapr+1 \mbox{ for some }\etapr\\[2mm]
      0 & \mbox{ otherwise.}
      \end{array}\right.\]
      
For the following argumentation, recall again properties and guiding picture of localization, cf.\ Lemma \ref{localipic}. 
We discuss the clauses of the definition of $\ual$:
\begin{enumerate}
\item
If $\eta=\etapr+1$ for some $\etapr$, then $\be\in(\ual,\ual^+=\al)$, and the $\Omj$-localization of $\ual$ is an initial
segment of the $\Omj$-localization of $\be$, cf.\ Remark \ref{noteworthyrmk}, hence $\be^\stark\ge\ual^\stark$, as $\ual$ then is a subterm of $\be$.
\item
If $\eta$ satisfies $F_{j+1}(\De,\eta)$, we have $\ual=\eta=\almmin$, thus again $\be\in(\ual,\ual^+)$ and  
$\be^\stark\ge\ual^\stark$. 
\item
$\eta=0$. Then $\ual=0$, unless $\De[0]^\starj<\almmin=\De^\starj$ where $m>1$, so that $\ual=\almmin$, in which
case again $\be\in(\ual,\ual^+)$ and thus $\be^\stark\ge\ual^\stark$ follows.
\end{enumerate}
Now that we have verified that $\be^\stark\ge\ual^\stark$, 
we also observe that we can not have $\Ga=\De$, since invoking Lemma \ref{auxlem} this would imply $\be\le\ual<\al[\ze]$.
Again, part \ref{bracketlocalizationclaim}, the assumption of part 7, and part 1 of Lemma \ref{localipic} imply that we can not have $\Ga>\De$ either, 
so we must have $\Ga<\De$.
Then $\De[\ze]\le\Ga<\De$, and the i.h.\ for $\De$ yields $\Ga^\stark\ge\De[\ze]^\stark$.
We finally obtain $\be^\stark\ge\max\{\De[\ze]^\stark,\ual^\stark\}=\al[\ze]^\stark$.
\qed

\begin{cor}\label{cantorianfundseqcor}
Any Buchholz system $\Ttcirc$ is a Cantorian system of fundamental sequences  
provided that its restriction to $\tau+1$ is a Cantorian system of fundamental sequences.\qed
\end{cor}

We now take a closer look at $\Omi$-localizations of elements of fundamental sequences. This will be crucial in proving 
Bachmann property in the next subsection.

\begin{lem}\label{localizationlem}
Let $\al=\thti(\De+\eta)>\Omi$ with $\Omi$-localization $(\Omi=\alnod,\ldots,\alm=\al)$ and $\ze<\aleph_{\domf(\al)}$.
We display the $\Omi$-localization of $\al[\ze]$, distinguishing between the cases of definition. 
\begin{enumerate}
\item $\eta\in\Lim$ such that $F_i(\De,\eta)$ does not hold. Then the $\Omi$-localization of $\al[\ze]$ is 
$(\alnod,\ldots,\almmin,\al[\ze])$, unless $\ze=0$ and $\al=\thti(\Omk)$ with $0<k\le i$, in which case the $\Omi$-localization of
$\al[0]$ is $(\Omi)$.
\item $\eta\not\in\Lim$ or $F_i(\De,\eta)$ holds.
\begin{enumerate}
\item $\De=0$. Here $\al[\ze]\in\Hz$ if and only if $\ze=0$. The $\Omi$-localization of $\al[0]$ then is $(\alnod,\ldots,\almmin)$,
unless $\eta=\etapr+1$ for some $\etapr>0$, in which case it is $(\alnod,\ldots,\almmin,\ual)$ where $\ual=\thti(\etapr)$.
\item $\chiomie(\De)=1$. The $\Omi$-localization of $\al[\ze]$ then is $(\alnod,\ldots,\almmin,\ual,\al[\ze])$ if $\eta=\etapr+1$
for some $\etapr$ where $\ual=\thti(\De+\etapr)$, otherwise it is $(\alnod,\ldots,\almmin,\al[\ze])$, unless $\ze=0$, $\De=\Omie$, 
and $\eta=0$, in which case it simply is $(\Omi)$.
\item $\De>0$ and $\chiomie(\De)=0$. Then the $\Omi$-localization of $\al[\ze]$ is $(\alnod,\ldots,\almmin,\al[\ze])$, unless
$\eta=\etapr+1$ for some $\etapr$, in which case it is $(\alnod,\ldots,\almmin,\ual,\al[\ze])$ where $\ual=\thti(\De+\etapr)$.
\end{enumerate}
\end{enumerate}
\end{lem}
{\bf Proof.} We discuss cases for $\al$. Note, first of all, that $\almmin\le\al[0]\le\al[\ze]<\al$ by parts \ref{fundseqclaim} and 
\ref{bracketlocalizationclaim} of Lemma \ref{bracketsmainlem}, so that $(\alnod,\ldots,\almmin)$ is an initial segment of the 
$\Omi$-localization of $\al[\ze]$ according to Lemma \ref{loclexordlem} on lexicographic monotonicity of localizations.
\\[2mm]
{\bf Case 1:} $\eta\in\Lim\cap\sup_{\si<\eta}\thti(\De+\si)$. Then $\domf(\al)=\domf(\eta)$ and $\al[\ze]=\thti(\De+\eta[\ze])$.
In this case the fixed point level of $\al[\ze]$ is equal to the fixed point level of $\al$, namely $\De$. This implies that
$\al[\ze]=\almmin$ if and only if $m=1$, $\De=0$, and $\eta[\ze]=0$, which in turn is only possible if $\ze=0$ and $\eta=\Omk$ 
for some $k\le i$, as $\eta<\Omie$.

Assume now that $\almmin<\al[\ze]$ and suppose $\be=\thti(\Ga+\rho)$ is the predecessor of $\al[\ze]$ in its $\Omi$-localization. 
Assume $\almmin<\be$.
By Lemma \ref{localipic} part 1 this implies that $\Ga\le\De$, but according to part 2 we must have $\Ga>\De$. 
So we must have $\be=\almmin$.
\\[2mm]
{\bf Case 2:} $\eta\not\in\Lim$ or $F_i(\De,\eta)$ holds. 
\\[2mm]
{\bf Subcase 2.1:} $\De=0$. Then we have $\domf(\al)=0$, $\al=\ual\cdot\om=\ual^+$, and $\al[n]=\ual\cdot(n+1)$ for $n<\om$.
If $\ual=\Omi$, we have $m=1$, and the $\Omi$-localization of $\al[0]=\ual$ is $(\Omi)$. 
If $F_i(0,\eta)$ holds, we have $m>1$ and $\al[0]=\ual=\almmin=\eta$ with $\Omi$-localization $(\alnod,\ldots,\almmin)$.
And if $\eta=\etapr+1$ for some $\etapr>0$, we have $\al[0]=\ual=\thti(\etapr)>\Omi$ with $\Omi$-localization 
$(\alnod,\ldots,\almmin,\ual)$. Note that $\al[n+1]=\ual\cdot(n+2)\not\in\Hz$.
\\[2mm]
{\bf Subcase 2.2} $\chiomie(\De)=1$. Then $\domf(\al)=0$, $\domf(\De)=i+1$, $\al[0]=\thti(\De[\ual])$ and 
$\al[n+1]=\thti(\De[\al[n]])$ for $n<\om$. Using part \ref{bracketparestimclaim} of Lemma \ref{bracketsmainlem} we have
$\ual\le\De[\ual]^\stari\le\max\{\De^\stari,\ual\}$ and $\De[\al[n]]^\stari=\al[n]$ for all $n<\om$.
We inspect the cases of definition of $\ual$.

\begin{enumerate} 
\item If $\eta=\etapr+1$ for some $\etapr$, we have $\ual=\thti(\De+\etapr)=\De[\ual]^\stari$, and it follows by induction on $n<\om$ that
the $\Omi$-localization of $\al[n]$ is $(\alnod,\ldots,\almmin,\ual,\al[n])$, noting that the fixed point level of $\al[n]$ is strictly less than
the fixed point level of $\al[n+1]$.

\item If $F_i(\De,\eta)$ holds, using Proposition \ref{fixpcondprop} and part \ref{bracketparestimclaim} of Lemma \ref{bracketsmainlem}
we have $\ual=\eta=\almmin=\De[\ual]^\stari$ with $m>1$, 
and as before induction on $n$ shows that the $\Omi$-localization of $\al[n]$ is $(\alnod,\ldots,\almmin,\al[n])$. 

\item If $\eta=0$, due to part \ref{bracketparestimclaim} of Lemma \ref{bracketsmainlem} we have $\De[0]^\stari=\De^\stari$
unless $\De[0]^\stari=0$ and $\De^\stari=1$ with $i=0$ and $\tau=1$, in which case $\almmin=\Om_0=1$, hence $m=1$.     
Thus $\ual=0$.
If $\De=\Omie$, we have $\al[0]=\Omi$ with $\Omi$-localization $(\Omi)$. Otherwise we have $\De>\Omie$, hence $\De[0]\ge\Omie$.
If $\De^\stari=\almmin$ we immediately see that the $\Omi$-localization of $\al[0]$ is $(\alnod,\ldots,\almmin,\al[0])$. 
Now assume that $\De^\stari>\almmin$.
Let $(\ga_0,\ldots,\ga_l)$ be the $\Omi$-localization of $\De^\stari$, and let $\ga_j=\thtj(\Ga_j+\rho_j)$ for $j=0,\ldots,l$.
Then according to Remark \ref{noteworthyrmk} $(\alnod,\ldots,\almmin)$ is a proper initial sequence of $(\ga_0,\ldots,\ga_l)$ and $\Ga_m,\ldots,\Ga_l<\De$ 
by part 1 of Lemma \ref{localipic} since $\eta=0$ and 
$\ga_m,\ldots,\ga_l\in(\almmin,\al)$. Assume that $\De[0]\le\Ga_j<\De$ for any $j\in\{m,\ldots,l\}$. Then part \ref{sandwichclaim}
of Lemma \ref{bracketsmainlem} yields $\Ga_j^\stari\ge\De[0]^\stari=\De^\stari$, which is impossible. Therefore, 
$\Ga_j<\De[0]$ for $j=m,\ldots,l$, and it follows that the $\Omi$-localization of $\al[0]$ is $(\alnod,\ldots,\almmin,\al[0])$, cf.\ Remark \ref{noteworthyrmk}.
Again, an induction on $n$ now yields that the $\Omi$-localization of $\al[n]$ is $(\alnod,\ldots,\almmin,\al[n])$, since, strictly increasing in $n$,
$\De[\al[n]]<\De$ and $\De[\al[n]]^\stari=\al[n]$.
\end{enumerate}
{\bf Subcase 2.3:} Otherwise. Then $\De>0$, $\chiomie(\De)=0$, $\domf:=\domf(\al)=\domf(\De)\le i$, and 
$\al[\ze]=\thti(\De[\ze]+\ual)$.
We are going to make use of Remark \ref{noteworthyrmk} and part \ref{ualcontrolclaim} of Lemma \ref{bracketsmainlem}, that is: $\ual>0\imp\De[\ze]^\stari<\ual$,  
and consider the cases of definition of $\ual$.
\begin{enumerate}
\item If $\eta=\etapr+1$ for some $\etapr$, then $\ual=\thti(\De+\etapr)$, and $(\De[\ze]+\ual)^\stari=\ual$, from which the $\Omi$-localization
of $\al[\ze]$ is seen to be $(\alnod,\ldots,\almmin,\ual,\al[\ze])$. 
\item If $F_i(\De,\eta)$ we have $\ual=\almmin>\De[\ze]^\stari$, hence the
$\Omi$-localization of $\al[\ze]$ is $(\alnod,\ldots,\almmin,\al[\ze])$. 
\item Assume that $\eta=0$. If in this situation $\ual>0$,
that is, $\ual=\almmin$ with $m>1$, we again immediately see that the $\Omi$-localization of $\al[\ze]$ is 
$(\alnod,\ldots,\almmin,\al[\ze])$.
\item Suppose finally that $\eta=\ual=0$, hence $\al=\thti(\De)$ and $\al[\ze]=\thti(\De[\ze])$.  
By part \ref{bracketlocalizationclaim} of Lemma \ref{bracketsmainlem} we know that $\almmin<\al[\ze]$, 
thus $\almmin\le\De[\ze]^\stari$.
In the case $\De[\ze]^\stari=\almmin$ we immediately see that $(\alnod,\ldots,\almmin,\al[\ze])$ is the $\Omi$-localization of $\al[\ze]$.
Suppose that $\De[\ze]^\stari>\almmin$ and let $(\ga_0,\ldots,\ga_l)$ be the $\Omi$-localization of $\De[\ze]^\stari$, 
where $\ga_j=\thti(\Ga_j+\rho_j)$ for $j=0,\ldots,l$.
Then, arguing as in the previous subcase,  $(\alnod,\ldots,\almmin)$ is a proper initial sequence of $(\ga_0,\ldots,\ga_l)$, and 
$\Ga_m,\ldots,\Ga_l<\De$ since $\ga_m,\ldots,\ga_l\in(\almmin,\al)$. 
Assuming that $\De[\ze]\le\Ga_j<\De$ for any $j\in\{m,\ldots,l\}$, part \ref{sandwichclaim}
of Lemma \ref{bracketsmainlem} yields $\Ga_j^\stari\ge\De[\ze]^\stari$, which is impossible. Thus
$\Ga_j<\De[\ze]$ for $j=m,\ldots,l$, and it follows that the $\Omi$-localization of $\al[\ze]$ is $(\alnod,\ldots,\almmin,\al[\ze])$.\qed
\end{enumerate}

\section{Bachmann property for \boldmath $\Ttcirc$\unboldmath}\label{BachmannTtsec}
We turn to proving the \emph{Bachmann property}, assuming that the restriction of our system to 
$\tau+1$ has got this property.
We first provide a lemma that comes in handy in order to modularize the proof of Bachmann property
and is useful itself, e.g.\ in proving Lemma \ref{Gnormlem}.

\begin{lem}\label{bachmannsquarebracketslem}
For $\Ga=\thtje(\Si+\si)$ and $\De=\thtje(\Xi+\xi)$ such that $\chiomie(\De)=1$ and $\De[\ze]<\Ga<\De$ we have
\[\De[\ze]\le\Ga[0].\]
\end{lem}
{\bf Proof.} Note that by assumption we have $\domf(\De)=i+1$ and hence $i\le j$. 
The assumptions also imply that $\Omje\le\De[\ze]<\Ga<\De$.
The claim is shown by induction on the build-up of $\De$ and side induction on the build-up of $\Ga$.
Denote the $\Omje$-localizations of $\Ga$ by $\Omje=\Ga_0,\ldots,\Ga_k=\Ga$ (not to be confused with the
enumeration of $\Ga$-numbers), and of $\De$ by 
$\Omje=\De_0,\ldots,\De_m=\De$, and note that by part \ref{bracketlocalizationclaim} of Lemma \ref{bracketsmainlem} 
we have both $\Ga_{k-1}\le\Ga[0]$ and $\De_{m-1}\le\De[\ze]$.
Note that by Lemma \ref{localipic} we must have \[\Si\le\Xi\] under these conditions, 
as a consequence of the fact that $\Ga\in(\De_{m-1},\De)$, which by Lemma \ref{loclexordlem} also implies that 
\[(\De_0,\ldots,\De_{m-1})\subseteq(\Ga_0,\ldots,\Ga_{k-1}),\] hence in particular $\De_{m-1}\le\Ga_{k-1}$.
Note further that the case where $\Ga=\nu\cdot\om$ for some $\nu$ is trivial, so that we may assume that $\Si>0$ whenever
$\si\not\in\Lim$ or $F_{j+1}(\Si,\si)$ holds.
\\[2mm]
{\bf Case 1:} $\xi\in\Lim\cap\sup_{\rho<\xi}\thtje(\Xi+\rho)$. Then we have $\De[\ze]=\thtje(\Xi+\xi[\ze])$ with $\chiomie(\xi)=1$. 
\\[2mm]
{\bf Subcase 1.1:} $\si\in\Lim\cap\sup_{\mu<\si}\thtje(\Si+\mu)$. Then we have $\Ga[0]=\thtje(\Si+\si[0])$.
\\[2mm]
{\bf 1.1.1:} If $\Si=\Xi$, it follows that $\xi[\ze]<\si<\xi$. The i.h.\ yields $\xi[\ze]\le\si[0]$, which implies the claim.
\\[2mm]
{\bf 1.1.2:} Suppose $\Si<\Xi$. Here we have $\eps_{\Omje+1}\le\De[\ze]\le(\Si+\si)^\starje<\Ga$, and since according to Lemma \ref{localipic}
we can not have $\De[\ze]\in(\Ga_{k-1},\Ga)$, it follows that $\De[\ze]\le\Ga_{k-1}\le\Ga[0]$ using again 
\ref{bracketlocalizationclaim} of Lemma \ref{bracketsmainlem}.
\\[2mm]
{\bf Subcase 1.2:} $\si\not\in\Lim$ or $F_{j+1}(\Si,\si)$. As mentioned before, we may assume that $\Si>0$, so that two possibilities
remain, which are treated almost identically, relying on the fact that $\De[\ze]$ has fixed-point level $\Xi$.
We have $\Ga[0]=\thtje(\Si[\uGa])$ if $\chiomjz(\Si)=1$ and $\Ga[0]=\thtje(\Si[0]+\uGa)$ if $\chiomjz(\Si)=0$, 
thus in either case $\uGa<\Ga[0]$.
We compare $\Si$ and $\Xi$.
\\[2mm]
{\bf 1.2.1:} $\Si=\Xi$. Here scenario (a) of part 2 of Lemma \ref{auxlem} applies, we thus obtain $\De[\ze]\le\uGa$, and clearly $\uGa<\Ga[0]$.
\\[2mm]
{\bf 1.2.2:} $\Si<\Xi$. This implies $\Xi>\Omjz$, and we obtain $\eps_{\Omje+1}<\De[\ze]\le\Ga_{k-1}\le\Ga[0]$ as we can not
have $\De[\ze]\in(\Ga_{k-1},\Ga)$ according to Lemma \ref{localipic}.
\\[2mm]
{\bf Case 2:} $\xi=0$, $\xi=\xipr+1$ for some $\xipr$, or $F_{j+1}(\Xi,\xi)$. We then have $\chiomie(\Xi)=1$, hence $\Xi>0$,
$\chiomjz(\Xi)=0$, and thus $\De[\ze]=\thtje(\Xi[\ze]+\uDe)$.
We further have $\De_{m-1}<\De[\ze]<\De$ by part \ref{bracketlocalizationclaim} of Lemma \ref{bracketsmainlem}
and due to the fact that $\Xi[\ze]<\Xi$. 
Note that $\Si=\Xi$ is impossible since by Lemma \ref{auxlem} we would have  $\Ga\le\uDe<\De[\ze]$.
We thus have $\Si<\Xi$. Note that in presence of $\chiomie(\Xi)=1$ this entails $\Xi>\Om_{j+2}$, 
and hence $\eps_{\Omje+1}\le\De[\ze]$. If $\Si<\Xi[\ze]$, we conclude as before in Subcases 1.1.2 and 1.2.2 that 
$\De[\ze]\le\Ga_{k-1}\le\Ga[0]$.
We may therefore assume that \[\Xi[\ze]\le\Si<\Xi\] and consider the definition of $\Ga[0]$. 
Note that according to part \ref{sandwichclaim}
of Lemma \ref{bracketsmainlem} we have $\Xi[\ze]^\stark\le\Si^\stark$ for $k\ge i$, hence in particular for $k=j+1$.
\\[2mm]
{\bf Subcase 2.1:} $\si\in\Lim\cap\sup_{\mu<\si}\thtje(\Si+\mu)$. Then we have 
$\Ga[0]=\thtje(\Si+\si[0])$.
\\[2mm]
{\bf 2.1.1:} $\Xi[\ze]=\Si$. Then we have $\uDe<\si$ and need to verify that even $\uDe\le\si[0]$ holds.
This is trivial in case $\si$ is additively decomposable, so assume that $\si\in\Hz$.
Checking the cases of definition of $\uDe$ we see that either $\uDe=0$, or $\si\in(\uDe,\uDe^+)$, so that according to part 1
of Remark \ref{noteworthyrmk} $\uDe$ is a predecessor of $\si$ in its $\Omje$-localization and thus $\uDe\le\si[0]$ by 
part \ref{bracketlocalizationclaim} of Lemma \ref{bracketsmainlem}.
\\[2mm]
{\bf 2.1.2:} $\Xi[\ze]<\Si$. We need to show that $\Xi[\ze]^\starje,\uDe<\Ga[0]$. Since  
$\Xi[\ze]^\starje\le\Si^\starje$, which in this subcase is strictly below $\Ga[0]$, we only need to consider cases of definition of
$\uDe$, showing that unless $\uDe=0$, we have $\Ga\in(\uDe,\uDe^+)$ and hence $\uDe<\Ga[0]$, using again 
part \ref{bracketlocalizationclaim} of Lemma \ref{bracketsmainlem} and the fact that the fixed point level of $\uDe$ is larger
than that of $\Ga[0]$.
\\[2mm]
{\bf Subcase 2.2:} $\si\not\in\Lim$ or $F_{j+1}(\Si,\si)$. As in Subcase 1.2 we may assume $\Si>0$ and consider $\chiomjz(\Si)$:
If $\chiomjz(\Si)=0$, we have $\Ga[0]=\thtje(\Si[0]+\uGa)$ and set $\Sinod:=\Si[0]$, and if $\chiomjz(\Si)=1$, we have 
$\Ga[0]=\thtje(\Si[\uGa])$ and set $\Sinod:=\Si[\uGa]$.
Both cases are dealt with similarly. We compare $\Xi[\ze]$ and $\Si$.
\\[2mm]
{\bf 2.2.1:} $\Xi[\ze]=\Si$. Noting that then we have $\uDe<\si$, in particluar $\si>0$, and $\Sinod<\Xi[\ze]$, scenario (b) of part 2 of 
Lemma \ref{auxlem} applies, and we obtain $\De[\ze]\le\uGa<\Ga[0]$. 
\\[2mm]
{\bf 2.2.2:} $\Xi[\ze]<\Si$. Then the i.h.\ applies, since we then have $\Xi[\ze]<\Si<\Xi$ with $\chiomie(\Xi)=1$. 
We obtain \[\Xi[\ze]\le\Si[0]\le\Sinod,\]
and since $\Xi[\ze]\le\Sinod<\Xi$ we have \[\Xi[\ze]^\starje\le\Sinod^\starje<\Ga[0]\]
using part \ref{sandwichclaim} of Lemma \ref{bracketsmainlem}, as even $\domf(\Xi)=i+1<j+2$.

Note that in case of $\uDe=0$ we directly obtain $\De[\ze]\le\Ga[0]$, so assume $\uDe>0$. 
By part \ref{ualcontrolclaim} of Lemma \ref{bracketsmainlem}, recalling that $\chiomjz(\Xi)=0$, we have \[\Xi[\ze]^\starje<\uDe.\] 
Since $\Xi>0$, we have $\uDe>\Omje$ in any case of definition of $\uDe$, so that $\Ga\in(\uDe,\uDe^+)$ implies 
\[\uDe\le\Ga_{k-1}\mbox{  with }k>1.\]
{\bf 2.2.2.1:} $\Xi[\ze]=\Si[0]$. If in this situation $\uGa>0$, we have $\Xi[\ze]<\Si[\uGa]$ in case of $\chiomjz(\Si)=1$, 
and noting that then $\Ga_{k-1}\le\uGa$, we obtain $\uDe\le\uGa$ and hence $\De[\ze]\le\Ga[0]$. 
On the other hand, $\uGa=0$ is impossible, as this would entail $\si=0$ and 
$\Ga_{k-1}\le\Si[0]^\starje=\Xi[\ze]^\starje<\uDe\le\Ga_{k-1}$, contradiction.
\\[2mm]
{\bf 2.2.2.2:} $\Xi[\ze]<\Si[0]$. Since $\Sinod<\Si$ we then have $\uDe\le\Ga_{k-1}<\Ga[0]$, 
so that $\Xi[\ze]^\starje,\uDe<\Ga[0]$
and hence $\De[\ze]<\Ga[0]$.
\qed

\begin{theo}\label{bachmanntheo}
Any Buchholz system over some $\tau$ has the Bachmann property and hence is a Cantorian Bachmann system, 
provided that its restriction to $\tau+1$ is a Cantorian system of fundamental sequences satisfying Bachmann property.
\end{theo}
{\bf Proof.} In order to prove the theorem, we prove the following more general
\\[2mm] 
{\bf Claim:} For any $\al\in\Ttcirc$, $\be\in\Tt$, and $n<\om$ such that $\al[n]<\be<\al$ we have $\al[n]\le\be[0]$.
\\[2mm]   
{\bf Proof of Claim.} The proof is by induction on the build-up of $\al$ and subsidiary induction
on the build-up of $\be$. For $\al\le\tau$ the claim holds by assumption, and for $\al=_\NF\xi+\eta$
the claim follows directly from the i.h., since the additive normal form of $\be$ then must begin with $\xi$, 
as $\cdot[\cdot]$ only affects the last summand. Now suppose that $\al$ is of
the form $\thti(\De+\eta)$ for some $i<\om$, where $\De+\eta>0$, and denote the $\Omi$-localization of $\al$ by
$\Omi=\alnod,\ldots,\alm=\al$. If $\al$ is of the form $\nu\cdot\om$, then the claim is trivial, as $\be$
must fit in between a finite multiple of $\nu=\ual$ and $\ual\cdot\om$ and hence be additively decomposable.
We may therefore assume that $\De>0$ whenever $\eta\not\in\Lim$ or $F_i(\De,\eta)$ holds. 
This implies that $\al[n]\in\Hz\setminus\Omi$, so that 
if $\be$ is additively decomposable, the claim again follows immediately.
Assume therefore that $\be=\thti(\Ga+\rho)>\Omi$, as $\al$, $\al[n]$, and $\be$ must belong to the same segment of
ordinals, i.e.\ $[\Omi,\Omie)$. 
We may further assume that $\Ga>0$ whenever $\rho\not\in\Lim$ or $F_i(\Ga,\rho)$ holds, since otherwise $\be=\ube\cdot\om$ 
and $\be[0]=\ube$, trivially implying the claim.
Let $\Omi=\be_0,\ldots,\be_l=\be$ be the $\Omi$-localization of $\be$, so that $l>0$ and $\be_{l-1}\le\be[0]$ according to 
part \ref{bracketlocalizationclaim} of Lemma \ref{bracketsmainlem}, which also implies that $\almmin\le\al[0]\le\al[n]$.
Since $\be\in(\almmin,\al)$, an application of Lemma \ref{loclexordlem}  yields  
\begin{equation}\label{initialsegequality}
(\alnod,\ldots,\almmin)\subseteq(\be_0,\ldots,\be_{l-1}),
\end{equation} 
hence in particluar $\almmin\le\be_{l-1}$ and $0<m\le l$.
Note also that by part 1 of Lemma \ref{localipic} we must have \[\Ga\le\De,\] since the assumption $\Ga>\De$ would entail 
$\be\le\almmin\le\al[n]$ as $\be<\al$, contradicting the assumption $\al[n]<\be$. 
\\[2mm]
{\bf Case 1:} $\eta\in\Lim\cap\sup_{\si<\eta}\thti(\De+\si)$. Then $\domf(\eta)=\domf(\al)=0$, $\al[n]=\thti(\De+\eta[n])$, 
and \[\thti(\De+\eta[n])<\thti(\Ga+\rho)<\thti(\De+\eta).\]

The fact that $\al[n]$ and $\al$ have the same fixed point level, namely $\De$, keeps this case easy to handle.
\\[2mm]
{\bf Case 1.1:} $\Ga<\De$. Due to part 1 of Lemma \ref{localipic} we then have $\al[n]\not\in(\be_{l-1},\be)$, 
and thus $\al[n]\le\be_{l-1}\le\be[0]$.
\\[2mm]
{\bf Case 1.2:} $\Ga=\De$. We distinguish whether $\thti$ is continuous at $\Ga+\rho$ or not:
\\[2mm]
{\bf 1.2.1:} $\rho\in\Lim\cap\sup_{\nu<\rho}\thti(\Ga+\nu)$. Then $\be[0]=\thti(\Ga+\rho[0])$, and 
it follows that $\eta[n]<\rho<\eta$. The i.h.\ yields $\eta[n]\le\rho[0]$, so that we obtain $\al[n]\le\be[0]$, as desired.
\\[2mm]
{\bf 1.2.2:} $\rho\not\in\Lim$ or $F_i(\Ga,\rho)$.
As mentioned earlier, we may then assume that $\Ga>0$, hence $\Ga\ge\Omie$.
Scenario (a) of part 2 of Lemma \ref{auxlem} applies, and we obtain $\al[n]\le\ube<\be[0]$.\\

For Cases 2 and 3 assume that either $\eta\not\in\Lim$ or that $F_i(\De,\eta)$ holds. 
As mentioned earlier, we may further assume that $\De>0$, hence $\De\ge\Omie$.
\\[2mm] 
{\bf Case 2:} $\chiomie(\De)=1$. Then we have $\domf(\De)=i+1$ while $\domf(\al)=0$, $\al[0]=\thti(\De[\ual])$, 
and $\al[n+1]=\thti(\De[\al[n]])$.
Let $\alpr$ be such that $\al[n]=\thti(\De[\alpr])$, that is $\alpr=\ual$ if $n=0$ and $\alpr=\al[n-1]$ if $n>0$.

The situation $\Ga=\De$ can not occur, as part 1 of Lemma \ref{auxlem} would yield $\be\le\ual<\al[0]$.
We can therefore only have $\Ga<\De$.
\\[2mm]
{\bf Case 2.1:} $\Ga<\De[\alpr]$. Then we obtain $\al[n]\le\be_{l-1}\le\be[0]$ again due to part 1 of Lemma \ref{localipic}.
\\[2mm]
{\bf Case 2.2:} $\Ga\ge\De[\alpr]$. Again, we first distinguish whether $\thti$ is continuous at $\Ga+\rho$ or not:
\\[2mm]
{\bf 2.2.1:} $\rho\in\Lim\cap\sup_{\nu<\rho}\thti(\Ga+\nu)$, hence $\be[0]=\thti(\Ga+\rho[0])$. 
According to part \ref{sandwichclaim} of Lemma \ref{bracketsmainlem} we then have 
$\De[\alpr]^\stari\le\Ga^\stari$ and thus also $\al[n]=\thti(\De[\alpr])\le\thti(\Ga+\rho[0])=\be[0]$.
\\[2mm]
{\bf 2.2.2:}  $\rho\not\in\Lim$ or $F_i(\Ga,\rho)$. As mentioned earlier, we may further assume that $\Ga>0$, hence $\Ga\ge\Omie$.
\\[2mm]
{\bf 2.2.2.1:} $\Ga=\De[\alpr]$.
Then we have $\rho>0$ with either $\rho=\rhopr+1$ and $\ube=\thti(\Ga+\rhopr)$ for some $\rhopr$, 
or $F_i(\Ga,\rho)$ and thus $\Ga^\stari<\rho=\ube=\be_{l-1}$ with fixed point level greater than $\Ga$ by Proposition \ref{fixpcondprop}, 
so that $\al[n]=\thti(\Ga)\le\ube<\be[0]$.
\\[2mm]
{\bf 2.2.2.2:} $\Ga>\De[\alpr]$. Note that according to part \ref{sandwichclaim} of Lemma \ref{bracketsmainlem} and by
Lemma \ref{bachmannsquarebracketslem} we have 
\[\De[\alpr]^\stari\le\Ga^\stari \quad\mbox{ and }\quad \De[\alpr]\le\Ga[0],\] 
and clearly $\Ga[0]^\stari\le\Ga[\ube]^\stari$ and $\Ga[0]\le\Ga[\ube]$. Now, if $\De[\alpr]=\Ga[0]$, then $\al[n]\le\be[0]$ is immediate,
and if $\De[\alpr]<\Ga[0]$, then part \ref{sandwichclaim} of Lemma \ref{bracketsmainlem} even yields 
$\De[\alpr]^\stari\le\Ga[0]^\stari$ and hence $\al[n]<\be[0]$.
\\[2mm] 
{\bf Case 3:} $\chiomie(\De)=0$. Then we have $\domf(\De)=\domf(\al)=0$ and $\al[n]=\thti(\De[n]+\ual)$. 

Note that as in Case 2 we can not have $\Ga=\De$, as using Lemma \ref{auxlem} this would entail $\be\le\ual<\al[n]$, 
contradicting the assumption $\al[n]<\be$. We therefore have $\Ga<\De$ throughout this case as well.
\\[2mm]
{\bf Case 3.1:} $\Ga<\De[n]$. As in cases 1.1 and 2.1 we obtain $\al[n]\le\be_{l-1}\le\be[0]$ due to part 1 of Lemma \ref{localipic}.
\\[2mm]
{\bf Case 3.2:} $\Ga\ge\De[n]$. Again, we first distinguish whether $\thti$ is continuous at $\Ga+\rho$ or not:
\\[2mm]
{\bf Case 3.2.1:} $\rho\in\Lim\cap\sup_{\nu<\rho}\thti(\Ga+\nu)$, hence $\be[0]=\thti(\Ga+\rho[0])$. 
\\[2mm]
{\bf 3.2.1.1:} $\Ga=\De[n]$. This implies $\ual<\rho<\al$, and we need to show that $\ual\le\rho[0]$.
Note that the claim follows immediately if $\rho\not\in\Hz$.
Assume that $\ual>0$ and consider the possible cases. Since $\Omi\le\ual<\rho$ and $\ual$ occurs in the $\Omi$-localization 
of $\rho$ by Lemma \ref{localipic} and Remark \ref{noteworthyrmk}, we obtain $\ual\le\rho[0]$ by part \ref{bracketlocalizationclaim} of 
Lemma \ref{bracketsmainlem}.  
\\[2mm]
{\bf 3.2.1.2:} $\Ga>\De[n]$. Here the claim easily follows if $\ual>0$, since then by part \ref{ualcontrolclaim} of 
Lemma \ref{bracketsmainlem} $\De[n]^\stari<\ual$, and, as a consequence of $\al\le\ual^+$ by Lemma \ref{localipic} and $\be\in(\ual,\ual^+)$,
by Remark \ref{noteworthyrmk} also $\ual\le\be_{l-1}<\be[0]$,
as $\be[0]$ has fixed-point level $\Ga\ge\Omie$.

Now assume that $\ual=0$, hence also $\eta=0$. 
Since $\almmin\le\be_{l-1}<\be[0]$ by (\ref{initialsegequality}) and since by Lemma \ref{localizationlem}
the $\Omi$-localization of $\al[n]$ is $(\alnod,\ldots,\almmin,\al[n])$, 
the assumption $\be[0]\in(\almmin,\al[n])$ would according to Lemma \ref{localipic} imply that the fixed point level $\Ga$ of $\be[0]$ 
is less than or equal to the fixed point level $\De[n]$ of $\al[n]$, contradicting our assumption in 3.2.1.2. Thus $\al[n]\le\be[0]$.
\\[2mm]
{\bf Case 3.2.2:} $\rho\not\in\Lim$ or $F_i(\Ga,\rho)$. 
As mentioned earlier, we may further assume that $\Ga>0$, hence $\Ga\ge\Omie$.
\\[2mm] 
{\bf 3.2.2.1:} $\Ga=\De[n]$.  Then scenario (b) of part 2 of Lemma \ref{auxlem} applies, and we obtain $\al[n]\le\ube<\be[0]$.
\\[2mm]
{\bf 3.2.2.2:} $\Ga>\De[n]$. The i.h.\ for $\De$ and $\Ga$ then yields \[\De[n]\le\Ga[0].\]

We observe that in the case $\ual>0$ we have $\ual>\Omi$  since $\De>0$ and $\ual$ is either equal to $\thti(\De+\etapr)$
where $\eta=\etapr+1$, or equal to $\almmin$ with $m>1$. If $\ual>0$ we also have
$\De[n]^\stari<\ual$ by part \ref{ualcontrolclaim} of Lemma \ref{bracketsmainlem}, and furthermore, since $\be\in(\ual,\ual^+)$ 
as $\al\le\ual^+$ by Lemma \ref{localipic}, by Remark \ref{noteworthyrmk} $\ual\le\be_{l-1}<\be[0]$ with $l>1$.  
\\[2mm]
{\bf 3.2.2.2.1:} $\De[n]=\Ga[0]$. If $\ual=0$ we are fine, so assume that $\ual>0$.
In case of $\ube>0$ we obtain $\ual\le\be_{l-1}\le\ube$ and are fine in both cases for $\chiomie(\Ga)$.
Assuming that $\ube=0$, using the observation from above, namely that $\ual$ must be an element of the $\Omi$-localization of
$\be$ according to Remark \ref{noteworthyrmk}, we would have $\ual\in\{\be_1,\ldots,\be_{l-1}\}$, 
so that $\ual$ would be a subterm of $\Ga[0]$ and hence of $\De[n]$. But this contradicts the fact that $\De[n]^\stari<\ual$.
\\[2mm]
{\bf 3.2.2.2.2:} $\De[n]<\Ga[0]$. In the case $\ual>0$ the claim follows immediately, as we already observed that then 
$\De[n]^\stari<\ual<\be[0]$. 
Assume $\ual=0$.  By Lemma \ref{localizationlem}, the $\Omi$-localization of $\al[n]$ then is $(\alnod,\ldots,\almmin,\al[n])$, 
and since by virtue of part 1 of Lemma \ref{localipic} the assumption $\be[0]\in(\almmin,\al[n])$ entails that the fixed point level of $\be[0]$, 
which depending on $\chiomie(\Ga)$ is either $\Ga[\ube]$ or $\Ga[0]$, is less than or equal to the fixed point level $\De[n]$ of $\al[n]$, 
noting that $\almmin\le\be_{l-1}<\be[0]$ by (\ref{initialsegequality}), we can only have $\al[n]\le\be[0]$.
\qed

\section{Norms, regularity, and Hardy hierarchy}\label{normsec}
Assuming that all additively decomposable terms in $\Tt$, including possible parameters $<\tau$ in the case $\tau\in\Ez$, 
are automatically given in additive normal form, every ordinal represented in $\Tt$ has a \emph{unique} term representation 
composed of the constant $0$, constants $\xi\in\Hz\cap\tau$,
the binary function $+$ for ordinal addition, and the unary functions $\thti$ for each $i<\om$,
since the ordinals $\tau=\Om_0,\Om_1,\Om_2,\ldots$ are represented as $\thtnod(0),\thte(0),\thtz(0),\ldots$.
Counting the number of function symbols $0$, $\xi$ for $\xi\in\Hz\cap\tau$, $+$, and $\thti$ for each $i<\om$ that occur in such 
unique term representation of an ordinal $\al$ in $\Tt$ provides us with a measure of term length of the unique representation of $\al$. 

Unless parameters $\xi\in\Hz\cap\tau$ for some $\tau\in\Ez$ are represented using a suitable notation system (say, $\Ts$ for 
some $\si<\tau$, allowing for a finite iteration of this process) and the measure of term length is accordingly modified, 
sets of ordinals $\be\in\Tt\cap\al$ of term length bounded by some $n<\om$ for given $\al\in\Tt$ are in general not finite. 
We therefore restrict ourselves to $\tau=1$ when using the notion of term length.
The measure of term length is then called the \emph{canonical norm} and gives rise to the notion of \emph{regularity} of a 
Bachmann system. Note however that, as just mentioned, there is a straightforward way to extend term length and canonical norm
to relativized systems $\Tt$ after replacing additively indecomposable parameters below $\tau$ by suitable term notations.

For more background on the origins of using norms of ordinals and an abstract treatment of 
fundamental sequences, normed Bachmann systems, and Hardy hierarchies we again refer to \cite{BCW94}.

\begin{defi}[cf.\ \cite{BCW94}]\label{normregularitydefi}  
\begin{enumerate}
\item For any ordinal $\al$ in $\Tt$ where $\tau=1$ the \emph{canonical norm} $\norm{\al}$ of $\al$ is defined as the total number 
of functions $0$, $+$, and $\thti$ for each $i<\om$ that occur in the unique term representation of $\al$.
\item A Buchholz system $(\Ttcirc,\cdot\{\cdot\})$ is called \emph{regular}, if (assuming $\tau=1$) for any ordinal 
$\al\in\Ttcirc\cap\Lim$ and any $\be\in\Tt\cap\al$ we have \[\be\le\al\{\norm{\be}\}.\]
For a given mapping $N:\Tt\to\om$ we call $(\Ttcirc,\cdot\{\cdot\},N)$ \emph{regular with respect to $N$}, if for any ordinal 
$\al\in\Ttcirc\cap\Lim$ and any $\be\in\Tt\cap\al$ we have \[\be\le\al\{N(\be)\}.\]
\item A system $(\Ttcirc,\cdot\{\cdot\},N)$ where $(\Ttcirc,\cdot\{\cdot\})$ is a Buchholz system and $N$ is a mapping 
$N:\Tt\to\om$ is called a \emph{normed Bachmann system} if 
\begin{enumerate}
\item for all $\al\in\Ttcirc$, $\be\in\Tt$, and $n<\om$ such that $\al\{n\}<\be<\al$ we have $N(\al\{n\})<N(\be)$, and
\item for all $\al\in\Ttcirc\cap\Lim$ we have $N(\al)\le N(\al\{0\})+1$.
\end{enumerate}
\item We call a mapping $N:\Tt\to\om$ a \emph{norm on $\Tt$} if for all $\al\in\Tt$ and all $n<\om$ the set
\[\{\be\in\Tt\cap\al\mid N(\be)\le n\}\]
is finite.
\end{enumerate}
We will also continue to use the more convenient notation $\al[n]$ and write $N\be$ instead of $N(\be)$.
\end{defi}

\begin{rmk} 
\begin{enumerate}
\item Note that for any $\al\in\Tt$ where $\tau=1$ and $n<\om$ the set \[\{\be\in\Tt\cap\al\mid\norm{\be}\le n\}\] is finite, since the 
occurrence of any $\thtk$ with large $k$ in a term $\be<\al$ must come with stepwise collapsing eventually down below 
$\al$ and hence must have a correspondingly large term length and thus canonical norm. 
\item In \cite{BCW94}, instead of $\Ttcirc$, countable initial segments of ordinals are considered as domains of systems
of fundamental sequences and Bachmann systems. In our context this would mean to consider $\Ttcirc\cap\Om_1$ instead
of $\Ttcirc$. However, as we have seen, we need the entire set $\Ttcirc$, or even $\Tt$, for inductive proofs to work.

In \cite{BCW94} it is shown that for a system $(\tau,\cdot[\cdot],N)$ where $\tau$ is a countable ordinal with mappings 
$\cdot[\cdot]:\tau\times\om\to\tau$ and $N:\tau\to\om$ such that
\\[1mm]
$(B1)\quad(\forall \al<\tau,\,n<\om)\;\left[0[n]=0\andsp(\al+1)[n]=\al\andsp[\al\in\Lim\imp\al[n]<\al[n+1]<\al]\right]$,
\\[1mm]
$(B3)\quad(\forall \al,\be<\tau,\,n<\om)\;\left[\al[n]<\be<\al\imp N\al[n]< N\be\right]$, and
\\[1mm]
$(B4)\quad(\forall \al\in\Lim\cap\tau)\;\left[N\al\le N\al[0]+1\right]$,
\\[1mm]
which is called a \emph{normed Bachmann system} in \cite{BCW94}, it \emph{follows} that (see Lemma 1 of \cite{BCW94} and
cf.\ Lemma \ref{Gnormlem} below):

$\;\;$(a) the system $(\tau,\cdot[\cdot],N)$ is a Bachmann system, i.e. property
\\[1mm]
$(B2)\quad(\forall \al,\be<\tau,\,n<\om)\;\left[\al[n]<\be<\al\imp\al[n]\le\be[0]\right]$ 
\\[1mm]
holds, since the assumption $\al[n]<\be<\al\andsp\be[0]<\al[n]$ would imply $N\al[n]<N\be\le N\be[0]+1\le N\al[n]$,

$\;\;$(b) $N$ is a norm, and

$\;\;$(c) $(\tau,\cdot[\cdot],N)$ is a regular Bachmann system, i.e.\ we have 
$(\forall \be<\al<\tau)\;\left[\be\le\al[N\be]\right]$.
\\[2mm]
Furthermore, for any $(\tau,\cdot[\cdot])$ satisfying $(B1)$ and $(B2)$, it \emph{follows} that 
$(\forall \al\in\Lim\cap\tau)\;\al=\sup\{\al[n]\mid n<\om\}$, cf.\ Lemma 3 of \cite{BCW94}.
\end{enumerate}
\end{rmk}

Before turning to another norm of interest, we prove regularity of $\Ttcirc$ regarding the canonical norm, where for simplicity 
we restrict ourselves to $\tau=1$.
\begin{theo} The Buchholz system $\Ttcirc$ for $\tau=1$ is regular: for any $\al\in\Ttcirc\cap\Lim$ and $\be\in\Tt\cap\al$ 
we have \[\be\le\al\{\norm{\be}\}.\]
\end{theo}
{\bf Proof.} The proof is by induction on the build-up of $\al$ and side induction on the build-up of $\be$.
Suppose first that $\al$ is additively decomposable, i.e.\ of a form $\al=_\NF\ga+\de$. If $\be\le\ga$ the claim is immediate.
Otherwise, $\be=\ga+\be_0$ for some $\be_0$ so that $\sumend(\ga)\le\mc(\be_0)$, $\mc(\be_0)$ denoting the maximal 
additively indecomposable component of $\be_0$. Then we have $\be_0<\de$, and the i.h.\
yields $\be_0\le\de\{\norm{\be_0}\}$, which implies the claim.
Now assume that $\al$ is of the form $\thti(\De+\eta)>\Omi$. Let $\be=_\ANF\be_1+\ldots+\be_m$ where $m\ge 1$, since the 
case $\be=0$ is trivial. If $\al=\nu\cdot\om$ for some $\nu$, it follows that $\nu=\ual$, $\be_1\le\ual$, and 
$\be\le\ual\cdot m\le\al\{\norm{\be}\}$, and we are done. Otherwise we have $\al\{n\}\in\Hz\setminus\Omi$ for any $n$. 
If $\be$ is additively decomposable,
i.e.\ $m>1$, then the s.i.h.\ yields $\be_1\le\al\{\norm{\be_1}\}<\al\{\norm{\be}\}$, and we obtain $\be\le\al\{\norm{\be}\}$.
Note further that if $\be\le\Omi$ the claim holds trivially since $\al\{0\}\ge\Omi$. Thus, we may assume that $\be$ is of the 
form $\be=\thti(\Ga+\rho)>\Omi$ from now on.
We have $\be<\al$, hence either $\Ga+\rho<\De+\eta$ and $(\Ga+\rho)^\stari<\al$, or $\be\le(\De+\eta)^\stari$.
\\[2mm]
{\bf Case 1:} $\eta\in\Lim$ and $\neg F_i(\De,\eta)$. Then $\al[n]=\thti(\De+\eta[n])$ for any $n$. Let $\eta=\etapr+\etanod$
be such that $\etanod\in\Hz$ and either $\etapr=0$ or $\eta=_\NF\etapr+\etanod$.
\\[2mm]
{\bf Subcase 1.1:} $\be\le(\De+\eta)^\stari$. If $\be\le(\De+\etapr)^\stari$, then $\be<\al[0]\le\al\{\norm{\be}\}$ is immediate.
Otherwise it follows that $\be\le\eta^\stari=\eta\in\Hz\setminus\Omi$. If $\be<\eta$, the i.h.\ for $\eta$ yields 
$\be\le\eta\{\norm{\be}\}<\al\{\norm{\be}\}$. So, we are left with the situation $\be=\eta>\De^\stari$. Since $F_i(\De,\eta)$
does not hold, we must have $\Ga\le\De$. If $\Ga=\De$, it follows that $\rho<\eta$, and the i.h.\ for $\eta$ yields
$\rho\le\eta\{\norm{\rho}\}<\eta\{\norm{\be}\}$, which implies the claim. Assume $\Ga<\De$ now. 
Note that we have $(\Ga+\rho)^\stari<\be\in\Ttcirc\cap\Lim$, so that the i.h.\ applies to $\be$, as $\be=\eta$ is a subterm of $\al$
in the current situation. We thus have 
$(\Ga+\rho)^\stari\le\be\{\norm{(\Ga+\rho)^\stari}\}<\be\{\norm{\be}\}<\thti(\De+\be\{\norm{\be}\})=\al\{\norm{\be}\}$.
Since $\Ga+\rho<\De+\eta\{\norm{\be}\}$, we obtain $\be<\al\{\norm{\be}\}$.
\\[2mm]
{\bf Subcase 1.2:} $\be>(\De+\eta)^\stari$. Then we must have $\Ga+\rho<\De+\eta$ and $(\Ga+\rho)^\stari<\al$, hence $\Ga\le\De$.
If $\Ga=\De$ and hence $\rho<\eta$, the i.h.\ for $\eta$ yields $\rho\le\eta\{\norm{\rho}\}<\eta\{\norm{\be}\}$. Thus
$\be=\thti(\De+\rho)<\thti(\De+\eta\{\norm{\be}\})=\al\{\norm{\be}\}$, as before.
If on the other hand $\Ga<\De$, hence $\Ga+\rho<\De+\eta\{\norm{\be}\}$, we use the s.i.h.\ for $(\Ga+\rho)^\stari<\al$ to obtain
$(\Ga+\rho)^\stari\le\al\{\norm{(\Ga+\rho)^\stari}\}<\al\{\norm{\be}\}$, which implies the claim.
\\[2mm]
{\bf Case 2:} $\eta\not\in\Lim$ or $F_i(\De,\eta)$. If $\De=0$, hence $\eta>0$ and $\al=\ual\cdot\om$, we obtain 
$\be\le\ual\le\ual\cdot(\norm{\be}+1)=\al\{\norm{\be}\}$ and are done. Assume $\De>0$ from now. 
If $\Ga=\De$, hence $\rho<\eta$, it follows that $\be\le\ual<\al[0]\le\al\{\norm{\be}\}$. And if $\Ga>\De$, it follows that
$\be\le\almmin\le\al[0]\le\al\{\norm{\be}\}$, since $\be\not\in(\almmin,\al)$ due to part 1 of Lemma \ref{localipic}. 
We are therefore left with the situation $\Ga<\De$.
Note that since $(\Ga+\rho)^\stari<\be<\al$, the s.i.h.\ yields $(\Ga+\rho)^\stari\le\al\{\norm{(\Ga+\rho)^\stari}\}<\al\{\norm{\be}\}$.
\\[2mm]
{\bf Subcase 2.1:} $\chiomie(\De)=0$. The i.h.\ for $\Ga<\De\in\Ttcirc\cap\Lim$ implies that 
$\Ga\le\De\{\norm{\Ga}\}<\De\{\norm{\be}\}$, hence $\Ga+\rho<\De\{\norm{\be}\}+\ual$, and as  $(\Ga+\rho)^\stari<\al\{\norm{\be}\}$,
we thus obtain $\be<\al\{\norm{\be}\}$.
\\[2mm]
{\bf Subcase 2.2:} $\chiomie(\De)=1$. 
We have $\norm{\be}>0$ by definition and $\al\{\norm{\be}\}=\thti(\De[\alpr])$ where $\alpr=\al\{\norm{\be}-1\}$. 
If $\Ga<\De[\alpr]$, the claim follows immediately. It turns out that this is the only possibility, since assuming otherwise 
we would have $\De[\alpr]\le\Ga<\De$,
and part \ref{sandwichclaim} of Lemma \ref{bracketsmainlem} would yield $\Ga^\stari\ge\De[\alpr]^\stari$, hence 
$\alpr\le\Ga^\stari\le\al\{\norm{\Ga^\stari}\}$
using part \ref{bracketparestimclaim} of Lemma \ref{bracketsmainlem} and the s.i.h.\ for $\Ga^\stari$.
But since $0<\norm{\Ga^\stari}<\norm{\be}-1$ as we at least need the functions $\thtie$ and $\thti$ to build up $\be$ from
$\Ga^\stari$, we would reach the contradiction $\alpr<\alpr$. 
\qed

We now turn to the definition of the norm of iterated application of $\cdot[0]$ to terms of $\Tt$ for systems $\Tt$ as in Definition
\ref{bsystemdefi}.
\begin{defi}\label{Gnormdefi} In the setting of Definition \ref{bsystemdefi}, for $\al\in\Tt$ let $\al[0]^0:=\al$ and 
$\al[0]^{i+1}:=(\al[0]^i)[0]$.
Define $G:\Tt\to\om$ by $G\al:=\min\{i\mid \al[0]^i=0\}$ as in Lemma 2 of \cite{BCW94}.
\end{defi}

\begin{lem}[cf.\ Lemmas 1 and 2 of \cite{BCW94}]\label{normedBachmannlem}  In the setting of Definition \ref{Gnormdefi}, 
the system $(\Ttcirc,\cdot[\cdot],G)$ is a normed and regular Bachmann system.
\end{lem}
{\bf Proof.} By definition of $G$ we have $G\al=G\al[0]+1$ for any nonzero $\al\in\Ttcirc$.
Since the Buchholz system $(\Ttcirc,\cdot[\cdot])$ enjoys Bachmann property, if $\al[n]<\be<\al$ where $\al\in\Ttcirc$
and $\be\in\Tt$, we must have $\al[n]=\be[0]^i$ for some $i>0$ and thus $G\al[n]<G\be$.  

In order to see that regularity follows, let $\al\in\Ttcirc\cap\Lim$ and observe that property (a) of normed Bachmann systems
implies that $n\le G\al[n]$ for every $n<\om$ as the sequence $(G\al[n])_{n<\om}$ must be strictly increasing, 
which in particular yields $G\be\le G\al[G\be]$ for any $\be\in\Tt$. 
Via contraposition we then obtain from property (a) of normed Bachmann systems that $\be\le\al[G\be]$ for all $\al\in\Ttcirc$ 
and all $\be\in\Tt\cap\al$, thus regularity with respect to $G$. Note that this derivation of regularity from ``normedness'' works
for any $N:\Tt\to\om$.
\qed

\begin{lem}\label{Gnormlem} Suppose $\tau=1$. For terms in $\Tt$ we have 
\begin{enumerate}
\item $G\al=G\be+G\ga$ if $\al=_\NF\be+\ga$.
\item For $\al\in\Tt$ such that $\domf:=\domf(\al)>0$ we have $G\al[\ga]\ge G\al[0]+G\ga$ for all $\ga\in\Tt\cap\aleph_\domf$.
\item $G\thti(\De+\eta)\ge G\De+G\eta+1$.
\end{enumerate}
\end{lem}
{\bf Proof.} Using induction on $\al$ the first part of the lemma follows directly from the definitions of $\cdot[\cdot]$ and $G$, 
since $\al[0]=\be+\ga[0]<\al$. 

For the second part, let $\al\in\Tt$ such that $\domf:=\domf(\al)>0$. If $\al=\thtie(0)$, the claim follows immediately. 
If $\al=_\NF\be+\de$, the claim follows by the i.h.\ for $\de$ using part 1. Assume that $\al=\thti(\De+\eta)$ where $\De+\eta>0$.
Since $\domf>0$, we must have $i>0$ and either $\al[\ga]=\thti(\De+\eta[\ga])$ in the case $\eta\in\Lim$ such that 
$\neg F_i(\De,\eta)$, or $\al[\ga]=\thti(\De[\ga]+\ual)$ with $\De>0$ and either $\eta\not \in\Lim$ or $\eta\in\Lim$ such
that $F_i(\De,\eta)$ holds. In case of $\ga=0$ the claim is immediate, so assume that $\ga>0$ and set $k:=G\ga$.
We then have $k>0$ and $0=\ga[0]^k<\ga[0]^{k-1}<\ldots<\ga[0]^0=\ga$, thus 
\[\al[0]<\al[\ga[0]^{k-1}]<\ldots<\al[\ga]<\al.\]
Lemma \ref{bachmannsquarebracketslem} yields \[\al[\ga[0]^{j+1}]\le\al[\ga[0]^j][0]\]
for $j=0,\ldots,k-1$, and hence \[G\al[0]<G\al[\ga[0]^{k-1}]<\ldots<G\al[\ga].\]
This entails \[G\al[\ga]\ge G\al[0]+G\ga\]
as claimed.  

The third part is shown by induction on $\al:=\thti(\De+\eta)$. If $\De+\eta=0$, we have $G\thti(0)=1$.
If $\eta\in\Lim$ such that $\neg F_i(\De,\eta)$, we apply the i.h.\ to $\al[0]=\thti(\De+\eta[0])$ and obtain
$G\al=G\thti(\De+\eta[0])+1\ge G(\De)+G(\eta[0])+2=G\De+G\eta+1$.
Now assume that either $\eta\not\in\Lim$ or $F_i(\De,\eta)$. If $\De=0$, the claim easily follows distinguishing between the two cases
$\eta=\etapr+1$ for some $\etapr$ or $F_i(0,\eta)$. Assume $\De>0$. 
If $\chiomie(\De)=0$, by the i.h.\ for $\al[0]<\al$ we have $G\al[0]=G\thti(\De[0]+\ual)\ge G\De[0]+G\ual+1=G\De+G\ual$,
and the cases $\eta=0$ and $F_i(\De,\eta)$ are immediate. Consider the case $\eta=\etapr+1$ for some $\etapr$, where
$\ual=\thti(\De+\etapr)$. The i.h.\ applied to $\ual$ yields $G\ual\ge G\De+G\etapr+1=G\De+G\eta\ge G\eta$.
If finally $\chiomie(\De)=1$, we apply the i.h.\ to $\al[0]=\thti(\De[\ual])$ and use part 2 to obtain
\[G\al=G\al[0]+1\ge G\De[\ual]+2\ge G\De[0]+G\ual+2=G\De+G\ual+1.\] 
Considering cases for $\eta$, the claim now follows as before. 
\qed

The above lemma allows us to straightforwardly prove the following generous upper bound of $\norm{\al}$ in terms of $G\al$.
\begin{lem}\label{normGnormestimlem}
Suppose $\tau=1$. For any $\al\in\Tt$ we have \[\norm{\al}\le (G\al+1)^2.\]
\end{lem}
{\bf Proof.} The lemma is shown by induction on the build-up of $\al\in\Tt$, frequently invoking the above Lemma \ref{Gnormlem}.

The claim is trivial for $\al=0,\thti(0)$, $i<\om$.
If $\al=_\NF\be+\ga$, we have $G\al=G\be+G\ga$ and 
\[\norm{\al}=\norm{\be}+\norm{\ga}+1\le(G\be+1)^2+(G\ga+1)^2+1\le(G\be+G\ga+1)^2,\] since $2G\be G\ga\ge2$.
Note that the above argument would not work for a linear expression, say, $c\cdot G\al+d$ for some $c,d<\om$, instead of $(G\al+1)^2$.
Now assume that $\al=\thti(\De+\eta)$ where $\De+\eta>0$.
\\[2mm]
{\bf Case 1:} $\eta\in\Lim$ such that $\neg F_i(\De,\eta)$. We have $G\al\ge G\De+G\eta+1$ according to the above lemma.
Hence $(G\al+1)^2\ge(G\De+G\eta+2)^2$, and using the i.h.\ for $\De$ and $\eta$ we have
\[\norm{\al}\le\norm{\De}+\norm{\eta}+2\le(G\De+1)^2+(G\eta+1)^2+2=(G\De)^2+(G\eta)^2+2(G\De+G\eta)+4<(G\al+1)^2.\]
{\bf Case 2:} $\eta\not\in\Lim$ or $F_i(\De,\eta)$.
\\[2mm]
{\bf Subcase 2.1:} $\De=0$. Then we have $\al=\ual\cdot\om$ and $G\al=G\ual+1$. If in this situation $\eta=\etapr+1$, 
then $\ual=\thti(\etapr)$ and $\norm{\al}=\norm{\eta}+1\le(G\eta+1)^2+1=(G\etapr+2)^2+1<(G\etapr+3)^2\le(G\al+1)^2$, 
using that $G\thti(\etapr)\ge G\etapr+1$. 
And if $F_i(\De+\eta)$, we have $\ual=\eta$, so that $\norm{\al}=\norm{\eta}+1\le (G\eta+1)^2+1<(G\eta+2)^2\le(G\al+1)^2$.
\\[2mm]
{\bf Subcase 2.2:} $\chiomie(\De)=0$ where $\De>0$. Then we have $\al[0]=\thti(\De[0]+\ual)$, 
\[G\al=G\al[0]+1\ge G\De[0]+G\ual+2=G\De+G\ual+1\] according to the above lemma, and consider cases for $\eta$.
If $\eta=0$, we have $\norm{\al}=\norm{\De}+1\le(G\De+1)^2+1<(G\al+1)^2$. Otherwise 
\[\norm{\al}=\norm{\De}+\norm{\eta}+2\le(G\De+1)^2+(G\eta+1)^2+2<(G\De+G\eta+2)^2\le(G\De+G\ual+2)^2\le(G\al+1)^2,\] 
since if $\eta=\etapr+1$ for some $\etapr$, we have $G\ual\ge G\De+G\etapr+1\ge G\eta$ while otherwise simply $G\ual=G\eta$.
\\[2mm]
{\bf Subcase 2.3:} $\chiomie(\De)=1$. Then we have $\al[0]=\thti(\De[\ual])$ and by the above lemma  
\[G\al=G\al[0]+1\ge G\De[\ual]+2\ge G\De[0]+G\ual+2=G\De+G\ual+1,\] whence the claim follows in similar fashion as in Subcase 2.2. 
\qed

\begin{cor} $G$ is a norm on $\Tt$ if $\tau=1$.
\end{cor} 
{\bf Proof.} By the above Lemma \ref{normGnormestimlem}, for any $\al\in\Tt$ and $n<\om$ the set
\[\{\be\in\Tt\cap\al\mid G(\be)\le n\}\]
is contained in the set $\{\be\in\Tt\cap\al\mid\norm{\be}\le(n+1)^2\}$, which is finite since there are only finitely many
terms of bounded length below $\al$. 
\qed

\begin{defi}[cf.\ \cite{BCW94}]\label{Hardyhierarchydefi}
For a Buchholz system $(\Ttcirc,\cdot[\cdot])$ as in Definition \ref{bsystemdefi}, set $\tauinf:=\Tt\cap\Om_1$ and define
the \emph{Hardy hierarchy} $(H_\al)_{\al<\tauinf}$ by \[H_0(n):=n\quad\mbox{ and }\quad H_\al(n):=H_{\al[n]}(n+1)\mbox{ for }\al>0.\]
\end{defi}

By Lemma \ref{normedBachmannlem} we know that the system $(\Ttcirc,\cdot[\cdot],G)$ is a normed and regular Bachmann system,
hence also the restriction $(\tauinf,\cdot[\cdot],G)$ is such a system, and setting $\tau=1$ for simplicity, also 
$(\Ttcirc,\cdot[\cdot],\norm{\cdot})$ and $(\oneinf,\cdot[\cdot],\norm{\cdot})$ are regular Bachmann systems, 
where $\oneinf$ is the proof theoretic ordinal of $\pioneonecanod$. 
\cite{BCW94} provides lemmas of basic properties of the Hardy hierarchy  (when based on a regular Bachmann system), 
which we include here for convenience as they illuminate the interplay of the notions involved.

\begin{lem}[Lemma 3 of \cite{BCW94}]\label{Hardyhierarchylem} Let $(\tauinf,\cdot[\cdot],N)$ be a regular Bachmann system.
\begin{enumerate}
\item $H_\al(n)<H_\al(n+1)$.
\item $\be[m]<\al<\be\imp H_{\be[m]}(n)\le H_\al(n)$.
\item $0<\al\andsp m\le n\imp H_{\al[m]}(n+1)\le H_\al(n)$.
\item $(\forall \be<\al)\;[N\be\le n\imp H_\be(n+1)\le H_\al(n)]$.
\item $H_\al(n)=\min\{k\ge n\mid\al[n][n+1]\ldots[k-1]=0\}=\min\{k\mid \al[n:k]=0\}$, where
\[\al[n:k]:=\al+(n\minusp k)\mbox{ for }k\le n,\quad\mbox{ and }\quad\al[n:k+1]:=(\al[n:k])[k]\mbox{ for }k\ge n.\]
\end{enumerate}
\end{lem}
Let $NF(\al,\be)$ abbreviate the expression stating that $\al,\be>0$ with Cantor normal forms 
$\al=_\CNF\om^{\al_0}+\ldots,\om^{\al_m}$ and $\be=_\CNF\om^{\be_0}+\ldots+\om^{\be_n}$ satisfying $\al_m\ge\be_0$.
Since the Buchholz systems considered here satisfy the properties 
\\[1mm]
$(B5^\prime)\quad(\forall \al,\be,n)\;[NF(\al,\be)\imp(\al+\be)[n]=\al+\be[n]]$ and
\\[1mm]
$(B6^\prime)\quad(\forall m,n<\om)\;[\om^{m+1}[n]=\om^m\cdot(n+1)]$,
\\[1mm]
i.e.\ are Cantorian in the sense of Definition \ref{fundseqdefi}, we have the following lemma, also cited from \cite{BCW94}.
\begin{lem}[Lemma 4 of \cite{BCW94}]\label{HardyhierarchyPRFlem}
In the same setting as in the previous lemma, we have
\begin{enumerate}
\item $NF(\al,\be)\imp H_\al(H_\be(n))\le H_{\al+\be}(n)$.
\item $(H_{\om^m})^{(n+1)}(n+1)\le H_{\om^{m+1}}(n).$
\item For each primitive recursive function $f$ there exists $m$ such that 
\[(\forall \vec{x})\;[f(\vec{x})<H_{\om^m}(\max\{\vec{x}\})].\]
\end{enumerate}
\end{lem}

It is well-known and it also follows from \cite{BCW94} that for any function $f$ that is provably total in $\pioneonecanod$ 
there is an $\al<\oneinf$ such that $f(\vec{x})<H_\al(\max\{\vec{x}\})$.

\section{Fundamental sequences for systems \boldmath$\Tqt$\unboldmath}\label{fundseqTtbarsec}

Fundamental sequences of the desired kind for $\Tqt$ are now obtained from the fundamental sequences for $\Tt$ by application of the strictly
increasing domain transformation functions $\itm$ and $\rtm$ for $m<\om$. We are going to characterize the system of fundamental sequences obtained 
in this straightforward way by an explicit definition on $\Tqt$-terms, which in turn is obtained by careful adaptation of definitions in Section \ref{fundseqTtsec}. 
We begin with defining the counterparts of Definitions \ref{chidefi} and \ref{dominddefi} for $\Tqomie$- and $\Tqt$-terms. We use the same names
as before for these auxiliary functions because of their exactly corresponding meaning.

\begin{defi}[cf.\ Definition \ref{chidefi}]\label{chiqdefi}
We define a characteristic function $\chiomie:\Tqomie\to\{0,1\}$, where $i<\om$, by recursion on the build-up of $\Tqomie$:
\begin{enumerate}
\item $\chiomie(\al):=\left\{\begin{array}{cl} 
                           0&\mbox{ if } \al<\Omie\\
                           1&\mbox{ if } \al=\Omie,
           \end{array}\right.$ 
\item $\chiomie(\al):=\chiomie(\eta)$ if $\al=_\NF\xi+\eta$,
\item $\chiomie(\al):=\left\{\begin{array}{cl}
          \chiomie(\De) & \mbox{ if } \eta\not\in\Lim\mbox{ or }F_j(\De,\eta)\\[2mm]
          \chiomie(\eta)& \mbox{ otherwise,}
      \end{array}\right.$\\[2mm]
if $\al=\thtqj(\De+\eta)>\Omie$ and hence $j\ge i+1$.
\end{enumerate}
\end{defi}

\begin{defi}[cf.\ Definition \ref{dominddefi}]\label{domindqdefi}
We define a domain indicator function $\domf:\Tqt\to\om$ recursively in the term build-up.
\begin{enumerate}
\item $\domf(\al):=0$ if $\al<\tau$,
\item $\domf(\al):=\domf(\eta)$ if $\al=_\NF\xi+\eta>\tau$,
\item for $\al=\thtqi(\De+\eta)$,
\begin{enumerate}
\item[3.1.] $\domf(\al):=\domf(\eta)$, if $\eta\in\Lim$ and $F_i(\De,\eta)$ does not hold,
\item[3.2.] if $\eta\not\in\Lim$ or $F_i(\De,\eta)$:
\begin{enumerate}
\item[3.2.1.] $\domf(\al):=\left\{\begin{array}{cl} i & \mbox{ if }\eta=0\\0 & \mbox{ otherwise}
       \end{array}\right\}$ in case of $\De=0$,\\[2mm] 
\item[3.2.2.] $\domf(\al):=0$ in case of $\chiomje(\De)=1$ {\bf\boldmath for some $j\ge i$\unboldmath},
\item[3.2.3.] $\domf(\al):=\domf(\De)$ otherwise.
\end{enumerate}
\end{enumerate}
\end{enumerate}
We define \[\Tqtcirc=\domfinv(0),\] 
which will turn out to characterize the set of terms of countable cofinality.
\end{defi}

The following lemma shows the partitioning of $\Tqt$ into terms of equal cofinality, using the just introduced auxiliary functions 
$\chi$ and $\domf$.
\begin{lem}[cf.\ Lemma \ref{partitioninglem}]\label{partitioningqlem}
$\Tqt$ is \emph{partitioned} into the union of disjoint sets
\[\Tqt=\Tqtcirc\;\dot{\cup}\;\sum_{i<\om}\{\al\in\Tqt\mid\chiomie(\al)=1\}.\]
\end{lem}
{\bf Proof.} 
The lemma is proved in the same way as Lemma \ref{partitioninglem}. First,
a straightforward induction on the build-up of terms shows for $i<k<\om$ that 
\[\chiomie(\al)+\chiomke(\al)<2\]
for all $\al\in\Tqomie$,
with the canonical embedding $\Tqomie\subseteq\Tqomke$. 
Defining \[M_0:=\{\al\in\Tqt\mid\chiomie(\al)=0\mbox{ for all }i<\om\}\quad\mbox{and }\quad 
                M_{i+1}:=\{\al\in\Tqt\mid\chiomie(\al)=1\}\mbox{ for }i<\om,\]
we see that the sets $(M_i)_{i<\om}$ are pairwise disjoint.

In order to prove the lemma, we are going to show the more informative claim that  
\begin{equation}
\domfinv(0)=M_0
\quad\mbox{ and }\quad
\domfinv(i+1)=M_{i+1}\mbox{ for }i<\om.
\end{equation} 
Since $\domf$ is a well-defined function on the entire domain $\Tqt$, we then obtain the desired partitioning result.
We proceed by induction on the build-up of terms in $\Tqt$ along the definition of $\domf$.
If $\al\le\tau$, we have $\domf(\al)=0$ and $\chiomie(\al)=0$ for all $i<\om$, hence $\al\in M_0$.
If $\al=_\NF\xi+\eta>\tau$, we have $\domf(\al)=\domf(\eta)$, and $\chiomie(\al)=\chiomie(\eta)$ for any $i<\om$, so that
the claim follows from the i.h. Suppose $\al=\thtqi(\De+\eta)>\tau$. 
\\[2mm]
{\bf Case 1:} $\eta\in\Lim$ and $F_i(\De,\eta)$ does not hold. By definition we have both $\domf(\al)=\domf(\eta)$ and
$\chiomje(\al)=\chiomje(\eta)$ for all $j<\om$, as is easily checked considering cases $i>j$ and $i\le j$. The i.h.\ thus
applies to $\eta$, and the claim follows. 
\\[2mm]
{\bf Case 2:} $\eta\not\in\Lim$ or $F_i(\De,\eta)$. 
\\[2mm]
{\bf Subcase 2.1:} $\De=0$. 
\\[2mm]
{\bf 2.1.1:} $\eta=0$. We then have $\al=\Omi$, and clearly $\Omi\in\domfinv(i)$. We already dealt with the case $i=0$ in this
context, and for $j$ such that $i=j+1$ we have $\chiomje(\Omi)=1$, that is, $\al\in M_i$.  
\\[2mm]
{\bf 2.1.2:} $\eta>0$. So, $\al=\thtqi(\eta)\in\domfinv(0)=\Tqtcirc$, and for $j<i$ we have $\chiomje(\al)=\chiomje(0)=0$, while 
$\chiomje(\al)=0$ for $i\le j$ since $\al<\Omje$. This shows that $\al\in M_0$.
\\[2mm]
{\bf Subcase 2.2:} $\chiomje(\De)=1$ for some $j\ge i$. Here we again have $\al\in\domfinv(0)=\Tqtcirc$. For any $k\ge i$ we have $\chiomke(\al)=0$
since $\al<\Omke$, and for $k<i$ we have $\chiomke(\al)=\chiomke(\De)=0$ since $k\not= j$ using the disjointness of the $M$-sets.
\\[2mm]
{\bf Subcase 2.3:} Otherwise, i.e.\ $\De>0$ such that $\chiomje(\De)=0$ for all $j\ge i$. Then by definition $\domf(\al)=\domf(\De)$, as well as
$\chiomke(\al)=\chiomke(\De)$ for all $k<\om$, checking cases $i>k$, $i=k$, and $i<k$. Now the i.h.\ applies to $\De$, and
the claim follows.
\qed

\begin{lem}\label{chitransformlem}
Let $\al\in\Tqomie$. We have \[\chiomie(\al)=\chiomie(\gOmie(\al)),\] 
and if $\al\in\domthtqj$ for some $j\ge i$, we also have \[\chiomie(\al)=\chiomie(\rtj(\al)).\]
\end{lem}
{\bf Proof.}  The proof is by induction on the term length of $\al$, and subsidiary induction on $\htomje(\al)-j$ for the second claim. 
Recall that $(\domthtqj)_{j<\om}$ is $\subseteq$-increasing.
\begin{enumerate}
\item
The interesting case regarding the first claim is where $\al$ is of a form $\al=\thtqj(\Xi+\De+\eta)$ where $j\ge i+1$, $\eta<\Omje\mid\De<\Omjz\mid\Xi$, 
and $\Xi+\De\in\domthtqj$. Then $\rtj(\Xi+\De+\eta)=\Ga+\eta$ where $\Ga:=\rtj(\Xi+\De)$, $\gOmie(\al)=\thtj(\Ga+\eta)$,
and Lemma \ref{fixpcondtransformlem} yields \[F_j(\Xi+\De,\eta)\quad\aeq\quad F_j(\Ga,\eta).\]
{\bf Case 1:} $\eta\not\in\Lim$ or $F_j(\Xi+\De,\eta)$. Then \[\chiomie(\al)=\chiomie(\Xi+\De)=\chiomie(\Ga)=\chiomie(\gOmie(\al))\] by the second part of 
the i.h.\ for the subterm $\Xi+\De$ of $\al$ and the definitions of $\chiomie$ for systems $\Tqt$ and $\Tt$, respectively.\\[2mm]
{\bf Case 2:} Otherwise. Then \[\chiomie(\al)=\chiomie(\eta)=\chiomie(\gOmie(\eta))=\chiomie(\gOmie(\al))\] by the first part of the i.h.\ for the term 
$\eta$ and again the respective definitions of $\chiomie$. 
\item We turn to the second claim. Suppose $\al=\Xi+\De+\eta\in\Tqomie\cap\domthtqj$ for some $j\ge i$, where according to our convention 
$\eta<\Omje\mid\De<\Omjz\mid\Xi$. As above we set $\Ga:=\rtj(\Xi+\De)$. By definition we have
\[\Ga=\left\{\begin{array}{ll}
  \Depr&\mbox{ if }\Xi=0\\
  \thtje(\rtje(\Xi))+\Depr&\mbox{ otherwise,}\end{array}\right.\]
with $\Depr$ defined as in Definition \ref{rtmdefi}, that is, $\Depr=0$ if $\De=0$, and for $\De=_\ANF\De_1+\ldots+\De_k$, $\De_l=\thtqje(\xi_l)$ and
$\Deprl:=\thtje(\rtje(\xi_l))$, $l=1,\ldots,k$, and $\Depr:=\Depre+\ldots+\Deprk$.

{\bf Case 1:} $\eta>0$. Then we have \[\chiomie(\al)=\chiomie(\sumend(\eta))=\chiomie(\gOmie(\sumend(\eta)))=\chiomie(\rtj(\al))\] using the first claim for
$\sumend(\eta)$.\\[2mm]
{\bf Case 2:} $\eta=0$ and $\De>0$. Then \[\chiomie(\al)=\chiomie(\De_k)=\chiomie(\gOmie(\De_k))=\chiomie(\Deprk)=\chiomie(\rtj(\al))\] using the
first claim for $\De_k$. \\[2mm]
{\bf Case 3:} $\eta=\De=0$. As the case $\al=0$ is trivial, assume that $\Xi>0$. We then have $\Ga=\thtje(\rtje(\Xi))$ where $\Omjz\mid\rtje(\Xi)$,
and by side i.h.\ for $j+1$ we obtain \[\chiomie(\al)=\chiomie(\Xi)=\chiomie(\rtje(\Xi))=\chiomie(\Ga)=\chiomie(\rtj(\al)).\] 
\end{enumerate}
This concludes the proof of the lemma.
\qed

\begin{lem}\label{dtransformlem} Let $\al\in\Tqt$. We have \[\domf(\al)=\domf(\gt(\al)).\]
\end{lem}
{\bf Proof.} This now follows by induction on the term length of $\al$. 
The interesting case is that $\al$ is of a form $\al=\thtqi(\Xi+\De+\eta)$ where $\eta<\Omie\mid\De<\Omiz\mid\Xi$.
We set $\Ga:=\rti(\Xi+\De)$, so that $\gt(\al)=\thti(\Ga+\eta)$.\\[2mm] 
{\bf Case 1:} $\eta\in\Lim$ and $\neg F_i(\Xi+\De,\eta)$. Since $F_i(\Xi+\De,\eta)$ if and only if $F_i(\Ga,\gt(\eta))$ by Lemma \ref{fixpcondtransformlem},
we have both $\domf(\al)=\domf(\eta)$ and $\domf(\gt(\al))=\domf(\gt(\eta))$, and the claim follows from the i.h.\ for $\eta$.\\[2mm]
{\bf Case 2:} Otherwise, that is, $\eta\not\in\Lim$ or $F_i(\Xi+\De,\eta)$ holds.\\[2mm] 
{\bf Subcase 2.1:} $\Xi+\De=0$. Here the claim follows immediately.\\[2mm]
{\bf Subcase 2.2:} $\chiomje(\Xi+\De)=1$ for some $j\ge i$. In this case we have $\domf(\al)=0$, $\Omje\mid\Xi+\De>0$.
By Lemma \ref{chitransformlem} we have \[\chiomle(\Xi+\De)=\chiomle(\rtl(\Xi+\De))\] for $l=i,\ldots,j$.
\\[2mm]
{\bf 2.2.1:} $j=i$. Then we have $\chiomie(\Ga)=1$ and $\domf(\gt(\al))=0$.\\[2mm]
{\bf 2.2.2:} $j>i$. This implies that $\De=0$ and $\Omje\mid\Xi>0$, and we observe that \[\Ga=\rti(\Xi)=\thtie(\ldots\thtj(\rtj(\Xi))\ldots),\]
and $\chiomje(\Xi)=\chiomje(\rtj(\Xi))=1$, while $\chiomle(\Xi)=\chiomle(\rtl(\Xi))=0$ for $l=i,\ldots,j-1$. 
By definition of $\domf(\gt(\al))$ we have \[\domf(\gt(\al))=\domf(\Ga)=\ldots=\domf(\rtjmin(\Xi))=\domf(\thtj(\rtj(\Xi)))=0,\]
since $\chiomje(\rtj(\Xi))=1$.  
\\[2mm]
{\bf Subcase 2.3:} Otherwise, that is, $\Xi+\De>0$ and $\chiomje(\Xi+\De)=0$ for all $j\ge i$, so that $\domf(\al)=\domf(\Xi+\De)$ and hence also
$\domf(\gt(\al))=\domf(\Ga)$ as $\chiomie(\Ga)=0$.\\[2mm]
{\bf 2.3.1:} $\De>0$. Then $\Xi+\De$ is not a multiple of $\Omiz$, and for $\Depr$ according to Definition \ref{rtmdefi} we have $\Depr=\gt(\De)$. 
The i.h.\ now yields \[\domf(\al)=\domf(\sumend(\De))=\domf(\sumend(\Ga))=\domf(\gt(\al)).\]
{\bf 2.3.2:} $\De=0$. Then $\Omiz\mid\Xi>0$ and $\Ga=\rti(\Xi)=\thtie(\rtie(\Xi))$. Let $j$ be maximal such that $\Omje\mid\Xi$, hence $j>i$.
We have $0=\chiomje(\Xi)=\chiomje(\rtj(\Xi))$, $\rtjmin(\Xi)=\thtj(\rtj(\Xi))$ and $\rtj(\Xi)=\gt(\Xi)$. By definition of $\domf(\Ga)$ we obtain
\[\domf(\gt(\al))=\domf(\Ga)=\ldots=\domf(\thtj(\rtj(\Xi))=\domf(\rtj(\Xi))=\domf(\gt(\Xi))=\domf(\Xi)=\domf(\al),\]
employing the i.h.\ for $\Xi$.
Thus, corresponding terms in $\Tqt$ and $\Tt$ have the same cofinality.
\qed

\begin{cor}\label{domfcor} For $\al\in\domthtqj$ we have \[\domf(\al)=\domf(\rtk(\al))\]
for all $k\ge\max\{j,\domf(\al)\minusp 1\}$.
\end{cor}
{\bf Proof.} This follows from Lemmas \ref{chitransformlem} and \ref{dtransformlem}, as follows.\\[2mm]
{\bf Case 1:} $\domf(\al)=i+1$, i.e.\ $\chiomie(\al)=1$. Then, since $1=\chiomie(\al)=\chiomie(\rtk(\al))$ for all $k\ge\max\{j,i\}$, we have
\[\domf(\rtk(\al))=i+1\] for all $k\ge\max\{j,i\}$.\\[2mm]
{\bf Case 2:} $\domf(\al)=0$, i.e.\ $\chiomie(\al)=0$ for all $i$. Let $k\ge j$, so that $\al\in\domthtqk$.\\[2mm]
{\bf 2.1:} $\Omkz\mid\al>0$. Then we have $\rtk(\al)=\thtke(\rtke(\al))$ where $\chiomkz(\rtke(\al))=\chiomkz(\al)=0$, hence
\[\domf(\rtk(\al))=\domf(\rtke(\al)).\]
{\bf 2.2:} Otherwise, i.e.\ $0<\sumend(\al)<\Omkz$. Then we have \[\domf(\rtk(\al))=\domf(\sumend(\rtk(\al)))=\domf(\sumend(\al))=\domf(\al).\]
Now the claim follows by induction.
\qed

We may now proceed to give the variant of Definition \ref{bsystemdefi} for $\Tqt$-terms. We will need to verify that terms in the definition below
are actually well-defined (regarding domains of $\thtqi$-functions) and indeed coincide with the terms obtained when using the detour via domain
transformation functions $\itm$ and $\rtm$ mentioned above.

\begin{defi}[cf.\ Definition \ref{bsystemdefi}]\label{bsystemqdefi}
For $\tau\in\Ezone\cap\aleph_1$ let $\cdot\{\cdot\}:(\tau+1)\times\N\to\tau$ be a base system according to Definition \ref{basesysdefi}.  
Fix the canonical assignment $\Om_0:=\tau$ and $\Omie:=\aleph_{i+1}$ for $i<\om$.
Let $\al\in\Tqt$. By recursion on the term length of $\al$ we define the function $\al[\cdot]:\aleph_d\to\Tqomd$ where $d:=\domf(\al)$.
Let $\ze$ range over $\aleph_d$.
\begin{enumerate}
\item $\al[\ze]:=\al\{\ze\}$ if $\al\le\tau$.
\item $\al[\ze]:=\xi+\eta[\ze]$ if $\al=_\NF\xi+\eta>\tau$.
\item For $\al=\thtqi(\De+\eta)$ where $i<\om$, $\thtq_0=\thtqt$, note that $d\le i$, and denote the $\Omi$-localization of 
$\al$ by $\Omi=\al_0,\ldots,\alm=\al$. We define a support term $\ual$ by
\[\ual:=\left\{\begin{array}{cl}
            \almmin & \mbox{ if either } F_i(\De,\eta)\mbox{, or: } \eta=0\mbox{ and }\De[0]^{\star_i}<\almmin=\De^{\star_i}
            \mbox{ where }m>1\\[2mm]
            \thtqi(\De+\etapr) & \mbox{ if } \eta=\etapr+1\\[2mm]
            0 & \mbox{ otherwise.}
     \end{array}\right.\] 
For $\al>\tau$ the definition then proceeds as follows.
\begin{enumerate}
\item[3.1.] If $\eta\in\Lim$ and $\neg F_i(\De,\eta)$, that is, $\eta\in\Lim\cap\sup_{\si<\eta}\thtqi(\De+\si)$, we have $d=d(\eta)$ 
and define 
   \[\al[\ze]:=\thtqi(\De+\eta[\ze]).\]
\item[3.2.] If otherwise $\eta\not\in\Lim$ or $F_i(\De,\eta)$, we distinguish between the following 3 subcases.
\begin{enumerate}
\item[3.2.1.] If $\De=0$, define \[\al[\ze]:=\left\{\begin{array}{cl} \ual\cdot(\ze+1) & \mbox{ if }\eta>0\mbox{ (and hence $d=0$)}\\
            \ze & \mbox{ otherwise.}
       \end{array}\right.\]
\item[3.2.2.] $\chiomje(\De)=1$ {\bf\boldmath for some $j\ge i$\unboldmath}. Then $d=0$ and, proceeding by induction on $j-i$, 
$\al[n]$ is defined recursively in $n<\om$ as follows.
\begin{enumerate}
\item[3.2.2.1.] $j=i$. Then we define
\[\al[0]:=\thtqi(\De[\ual])\quad\mbox{ and }\quad\al[n+1]:=\thtqi(\De[\al[n]]).\]
\item[3.2.2.2.] $j>i$. Then we set $\Si:=\thtqj(\De)$, for technical convenience $\Si[-1]:=0$ as well as 
\[\Sipr[n]:=\left\{\begin{array}{cl}
             \Si[n] & \mbox{ if }\De=\Omje\\[2mm]
             \Si[n-1] & \mbox{ otherwise,}
             \end{array}\right.\] 
and define
\[\al[n]:=\thtqi(\De[\Sipr[n]]+\ual).\]
\end{enumerate}
\item[3.2.3.] Otherwise. Then $d=d(\De)$ and \[\al[\ze]:=\thtqi(\De[\ze]+\ual).\]
\end{enumerate}
\end{enumerate} 
\end{enumerate}
We call the system $(\Tqtcirc,\cdot\{\cdot\})$ 
(more sloppily also simply the entire mapping $\cdot\{\cdot\}$),  where the mapping $\cdot\{\cdot\}$ is simply the restriction
of $\cdot[\cdot]$ to $\Tqtcirc$, 
a \emph{Buchholz system over $\tau$} (for $\Tqt$-terms). Note that this system is determined uniquely modulo the choice of 
$\cdot\{\cdot\}:(\tau+1)\times\N\to\tau$, which in turn is trivially determined if $\tau=1$.   
\end{defi}

\begin{theo}\label{Tqtmaintheo} The system of fundamental sequences for $\Tqt$ obtained by setting \[\al[\ze]:=\ft(\gt(\al)[\ze])\] for $\al\in\Tqt$ 
and $\ze<\aleph_d$, where $d:=\domf(\al)$, is characterized by the above Definition \ref{bsystemqdefi}. 
We therefore have \[\al[\ze]=\gt(\al)[\ze].\] 
\end{theo}
{\bf Proof.} The proof is by induction on the term length of $\al\in\Tqt$, simultaneously proving the following
\begin{claim}\label{theoclaim} If $\al\in\domthtqk$ such that $d\le k+1$, or all $\ze<\aleph_d$ we have \[\al[\ze]=\itk(\rtk(\al)[\ze]),\]
hence in particular $\al[\ze]\in\domthtqk$ and thus \[\rtk(\al[\ze])=\rtk(\al)[\ze].\]
\end{claim}
We will be applying the results from Subsection \ref{orderisosubsec}, in particular Lemmas \ref{fixpcondtransformlem} and \ref{localizationtransformlem},
frequently without explicit mention.

Assuming that Claim \ref{theoclaim} holds for terms of $\Tqt$ strictly shorter than $\al$, we first prove the claim of the theorem.
Since other cases then are trivial, assume that $\al$ is of the form $\al=\thtqi(\Xi+\De+\eta)$ where $\eta<\Omie\mid\De<\Omiz\mid\Xi$, so that 
$\gt(\al)=\thti(\Ga+\eta)$ where $\Ga:=\rti(\Xi+\De)$.\\[2mm]
{\bf Case 1:} $\eta\in\Lim$ and $\neg F_i(\Xi+\De,\eta)$. Then we have $\gt(\al)[\ze]=\thti(\Ga+\gt(\eta)[\ze])$ and $\eta[\ze]=\gt(\eta)[\ze]$ by the i.h.\ for 
$\eta$, so that
\[\al[\ze]=\ft(\gt(\al)[\ze])=\thtqi(\Xi+\De+\eta[\ze]),\] in consistence with Definition \ref{bsystemqdefi}.\\[2mm]
{\bf Case 2:} $\eta\not\in\Lim$ or $F_i(\Xi+\De,\eta)$ holds.\\[2mm]
{\bf Subcase 2.1:} $\Xi+\De=0$. We have $\underline{\gt(\al)}=\ual$ inspecting the definitions of $\ual$ and $\underline{\gt(\al)}$, hence 
\[\al[\ze]=\ft(\gt(\al)[\ze])=\left\{\begin{array}{cl} \ual\cdot(\ze+1) & \mbox{ if }\eta>0\\
            \ze & \mbox{ otherwise.}
       \end{array}\right.\]
{\bf Subcase 2.2:} $\chiomje(\Xi+\De)=1$ for some $j\ge i$.
In this case we have $\domf(\al)=0$, $\Omje\mid\Xi+\De>0$.
By Lemma \ref{chitransformlem} we have \[\chiomle(\Xi+\De)=\chiomle(\rtl(\Xi+\De))\] for $l=i,\ldots,j$, and according to Lemma \ref{dtransformlem}
we have $\domf(\gt(\al))=0$. 
\\[2mm]
{\bf 2.2.1:} $j=i$. Then we have $\chiomie(\Xi+\De)=\chiomie(\Ga)=1$, hence $\domf(\Xi+\De)=\domf(\Ga)=i+1$, and  
\[\gt(\al)[0]=\thti(\Ga[\underline{\gt(\al)}])\quad\mbox{ and }\quad\gt(\al)[n+1]=\thti(\Ga[\gt(\al)[n]])\] according to Definition \ref{bsystemdefi}. 
According to part \ref{bracketparestimclaim} of Lemma \ref{bracketsmainlem} we have $\underline{\gt(\al)}=0$ if $\eta=0$.
By Claim \ref{theoclaim} for $\Xi+\De$ we have \[\rti(\Xi+\De[\ze])=\Ga[\ze]\] for all $\ze<\Omie$, and since $\rti$ does not change the 
$\mbox{}^\stari$-value, it follows that if $\eta=0$, we also have $\ual=0$. This shows that $\underline{\gt(\al)}=\ual$, and we obtain recursively in $n<\om$
\[\al[0]=\ft(\gt(\al)[0])=\thtqi(\Xi+\De[\ual])\quad\mbox{ and }\quad\al[n+1]=\ft(\gt(\al)[n+1])=\thtqi(\Xi+\De[\al[n]]).\]
{\bf 2.2.2:} $j>i$. This implies that $\De=0$ and $\Omje\mid\Xi>0$, and we observe that 
\[\Ga=\rti(\Xi)=\thtie(\ldots\thtj(\rtj(\Xi))\ldots)\quad\mbox{ and }\quad\domf(\Ga)=0.\]

Note first that the $\Omie$-localization of $\Ga$ is $(\Omie,\Ga)$ (considering $\Omie$ the $0$-th element of the localization), since 
$\iti(\Ga)=\iti(\rti(\Xi))=\Xi$. According to Lemma \ref{localizationlem} the $\Omie$-localization of $\Ga[n]$ is then $(\Omie,\Ga[n])$,
unless $\Xi=\Omiz$ and $n=0$, in which case it simply is $(\Omie)$.

According to Definition \ref{bsystemdefi} we have 
\[\gt(\al)[n]=\thti(\Ga[n]+\underline{\gt(\al)}),\]
and we need to verify that, setting $\Si:=\thtqj(\Xi)$, $\Si[-1]:=0$, and 
$\Sipr[n]:=\left\{\begin{array}{cl} \Si[n] & \mbox{ if }\Xi=\Omje\\ \Si[n-1] & \mbox{ otherwise,}\end{array}\right.$ we have
\begin{eqnarray}\label{verifytask} \iti(\Ga[n])=\Xi[\Sipr[n]]\quad\mbox{ and }\quad\ual=\underline{\gt(\al)}.
\end{eqnarray}

Note that the fundamental sequence for $\Si$ is dealt with as in Subcase 2.2.1, 
so that by the i.h.\ we have
\[\Si[0]=\thtqj(\Xi[0])=\thtj(\rtj(\Xi)[0])=\gt(\Si)[0]\quad\mbox{ and }\quad\Si[n+1]=\thtqj(\Xi[\Si[n]])=\thtj(\rtj(\Xi)[\gt(\Si)[n]])=\gt(\Si)[n+1],\]
since $\uSi=0=\underline{\gt(\Si)}$ as $\Xi[0]^\starj=\Xi^\starj$ because of $\chiomje(\Xi)=1$, hence by Claim \ref{theoclaim} for $\Xi$
also $(\rtj(\Xi)[0])^\starj=\rtj(\Xi)^\starj$.
In order to show (\ref{verifytask}) we proceed by side induction on $j-i>0$, showing along the way that 
\begin{eqnarray}\label{Gaclaim} \uGa=0\quad\mbox{ and }\quad\Ga[n]=\thtie(\ldots\thtj(\rtj(\Xi[\Si[n-1]]))\ldots).\end{eqnarray}
We see that $\underline{\gt(\al)}=\ual$ follows immediately if either $\eta$ is of a form $\etapr+1$ or $F_i(\Xi,\eta)$ holds. If $\eta=0$, it follows that
$\ual=0$ since $\Xi[0]^\stari=\Xi^\stari$ because of $\chiomje(\Xi)=1$, so that $\Xi[0]$ and $\Xi$ have the same $\mbox{}^\stari$-values.
Using (\ref{Gaclaim}) we see that $\Ga[0]^\stari=\Xi[0]^\stari=\Xi^\stari=\Ga^\stari$ as $\rtj$ does not change $\mbox{}^\stari$-values for $j\ge i$.
Thus also $\underline{\gt(\al)}=0$. We now consider cases regarding $j-i$.
\begin{itemize}
\item If $j=i+1$, as shown above we have $\Ga=\gt(\Si)$, $\uGa=0$, and $\Ga[n]=\thtj(\rtj(\Xi[\Si[n-1]]))$, satisfying (\ref{Gaclaim}).
In order to calculate $\iti(\Ga[n])$, 
we first consider the case $\Xi=\Omiz$ and obtain \[\iti(\Ga[n])=\iti(\thtie^{(n+1)}(0))=\thtqie^{(n+1)}(0)=\Si[n]=\Sipr[n].\]
If $\Xi>\Omiz$, recall that the first element of the $\Omie$-localization of $\Ga[n]$ is $\Ga[n]$. We then either have $\Xi=\Xipr+\Omiz$ for some 
$\Xipr>0$, which entails \[\iti(\Ga[n])=\Xipr+\Si[n-1]=\Xi[\Sipr[n]],\] 
or we have $\Omiz\mid\Xi[\Si[n-1]]>0$ for all $n<\om$, which yields \[\iti(\Ga[n])=\Xi[\Sipr[n]].\]
\item If $j>i+1$, we have $\Ga=\thtie(\Pi)$ where $\Pi:=\thtiz(\ldots\thtj(\rtj(\Xi))\ldots)$ and use the s.i.h.\ for $\Pi$, that is, $i+1$.
We then have $\chiomiz(\Pi)=0$, $\Pi[n]=\thtiz(\ldots\thtj(\rtj(\Xi[\Si[n-1]]))\ldots)$,
hence according to Definition \ref{bsystemdefi} \[\Ga[n]=\thtie(\Pi[n]),\] as $\uGa=0$ since $\Pi[0]^\starie=\Pi^\starie$ because of $\Xi[0]^\starie=\Xi^\starie$.
The s.i.h.\ also immediately yields \[\iti(\Ga[n])=\itie(\Pi[n])=\Xi[\Sipr[n]].\]
\end{itemize}
{\bf Subcase 2.3:} Otherwise, that is, $\Xi+\De>0$ and $\chiomje(\Xi+\De)=0$ for all $j\ge i$, so that $\chiomie(\Ga)=0$ by Lemma \ref{chitransformlem}
and \[d:=\domf(\al)=\domf(\Xi+\De)=\domf(\Ga)=\domf(\gt(\al))\le i\] by Corollary \ref{domfcor}.
Then for all $\ze<\aleph_d$ we have $\gt(\al)[\ze]=\thti(\Ga[\ze]+\underline{\gt(\al)})$ according to Definition \ref{bsystemdefi}, and by Claim \ref{theoclaim}
for $\Xi+\De$ \[\rti(\Xi+\De[\ze])=\Ga[\ze].\] 
The equality $\underline{\gt(\al)}=\ual$ is clear if either $\eta=\etapr+1$ for some $\etapr$ or if $F_i(\Xi+\De,\eta)$ holds. In the case $\eta=0$, again using 
Claim \ref{theoclaim} for $\Xi+\De$, we observe
that \[(\Xi+\De[0])^\stari<\almmin=(\Xi+\De)^\stari\mbox{ with }m>1\quad\aeq\quad\Ga[0]^\stari<\almmin=\Ga^\stari\mbox{ with }m>1,\]
where $m$ is the maximal index of the $\Omi$-localization of $\al$.
Thus \[\al[\ze]=\ft(\gt(\al)[\ze])=\thtqi(\iti(\Ga[\ze]+\underline{\gt(\al)}))=\thtqi(\Xi+\De[\ze]+\ual),\] as claimed.\\[2mm]
{\bf Proof of Claim \ref{theoclaim}.} Suppose that $\al\in\domthtqk$ where $d=\domf(\al)\le k+1$. We need to show that
\[\itk(\rtk(\al)[\ze])=\al[\ze]\]
for $\ze<\aleph_d$. To this end we write $\al$ in the form $\al=\Xi+\De+\eta$ where $\eta<\Omke\mid\De<\Omkz\mid\Xi$, and set $\Ga:=\rtk(\Xi+\De)$,
so that $\rtk(\al)=\Ga+\eta$. As the case $\al=0$ is trivial, we assume that $\al>0$. Note also that the claim is trivial for $\al<\Omke$, cf.\ Case 1
below. We proceed by subsidiary induction on $\kpr\minusp k$ where $\kpr$ is such that $\al\in[\Om_{\kpr},\Om_{\kpr+1})$.\\[2mm]
{\bf Case 1:} $\eta>0$. Then the claim follows immediately from the i.h.\ for $\eta$, as clearly $\eta\in\domthtqk$:
\[\itk(\rtk(\al)[\ze])=\itk(\Ga)+\eta[\ze]=\Xi+\De+\eta[\ze].\]
{\bf Case 2:} $\eta=0$ and $\De>0$. 
We write $\De=_\ANF\De_1+\ldots+\De_j$, $\De_l=\thtqke(\xi_l)$, and $\Deprl:=\gt(\De_l)=\thtke(\rtke(\xi_l))$ for $l=1,\ldots,j$.
We further set $\Xipr:=\rtke(\Xi)$ due to the frequent occurrence of this term in the sequel.\\[2mm]
{\bf Subcase 2.1:} $j>1$. Then we have $\Ga[\ze]=\rtk(\Xi+\De_1+\ldots+\De_{j-1})+\gt(\De_j)[\ze]$ and using the i.h.\ for $\De_j$
\[\itk(\rtk(\al)[\ze])=\itk(\Ga[\ze])=\Xi+\De_1+\ldots+\De_{j-1}+\De_j[\ze]=\al[\ze].\]
{\bf Subcase 2.2:} $j=1$. We then write $\xi_1=\Si+\rho$ such that $\rho<\Omkz\mid\Si$ and set $\Sipr:=\rtke(\Si)$, so that $\Depre=\thtke(\Sipr+\rho)$.
We inspect cases of the definition of $\Depre[\ze]$.
\\[2mm]
{\bf 2.2.1:} $\xi_1=0$, that is, $\De_1=\Omke$. Then $\Ga[\ze]=\rtk(\Xi)+\ze$ and $\itk(\Ga[\ze])=\Xi+\ze=\al[\ze]$.
\\[2mm]
{\bf 2.2.2:} $\rho\in\Lim$ and $\neg F_{k+1}(\Sipr,\rho)$. Then $\Depre[\ze]=\thtke(\Sipr+\rho[\ze])$. We inspect the definition of $\Ga=\rtk(\Xi+\De)$.
\begin{itemize}
\item If in this situation $\Xi=0$, we have $\itk(\Ga)=\De_1$ with $\Ga=\Depre$, and the first element of the $\Omke$-localization of $\Depre$ cannot
be an epsilon number, hence $\Sipr=0$, $\Ga[\ze]=\thtke(\rho[\ze])$, and the $\Omke$-localization of $\Ga[\ze]$ is either $(\Omke,\Ga[\ze])$ or simply
$(\Omke)$. Applying the (main) i.h.\ to $\De_1$, we therefore obtain \[\itk(\Ga[\ze])=\thtqke(\rho[\ze])=\De_1[\ze].\]
\item $\Xi>0$ such that $\thtke(\Xipr)<\Depre$. We then have $\Ga=\Depre$, but $\itk(\Ga)=\Xi+\De_1$. The first element of the $\Omke$-localization of
$\Depre$ is of a form $\thtke(\Xipr+\nu)$ for some $\nu<\Omke$. By Lemma \ref{localizationlem} the first element of $\Ga[\ze]=\thtke(\Sipr+\rho[\ze])$
must have fixed point level $\Xipr$ since $\Sipr\le\Xipr$. Hence using the (main) i.h.\ for $\De_1$ \[\itk(\Ga[\ze])=\Xi+\De_1[\ze].\]
\item $\thtke(\Xipr)\ge\Depre$. We then have $\Ga=\thtke(\Xipr)+\Depre$ and $\Ga[\ze]=\thtke(\Xipr)+\Depre[\ze]$. The first element of the 
$\Omke$-localization of $\thtke(\Xipr)$ must be itself since $\itk(\Ga)=\Xi+\De_1$, hence again $\itk(\Ga[\ze])=\Xi+\De_1[\ze]$ by the i.h.\ for $\De_1$.
\end{itemize}
{\bf 2.2.3:} Otherwise, that is, $\rho\not\in\Lim$ or $F_{k+1}(\Sipr,\rho)$ holds. Here Definition \ref{bsystemdefi} splits into the following cases.
\\[2mm]
{\bf 2.2.3.1:} $\Sipr=0$. We then have $\Depre=\thtke(\rho)$ where $\rho>0$, $\domf(\Depre)=0$, and $\Depre[n]=\uDepre\cdot(n+1)$ as well as
$\uDepre=\gt(\uDee)$ and $\De_1[n]=\uDee\cdot(n+1)$. 
\begin{itemize}
\item If in this situation $\Xi=0$, we have $\Ga=\Depre>\Omke$, $\itk(\Ga)=\De_1$, and the $\Omke$-localization of $\Ga$ is $(\Omke,\Ga)$ as $\Xipr=0$.
The $\Omke$-localization of $\uDepre$ then is $(\Omke)$ in case of $\rho=1$ and $(\Omke,\uDepre)$ otherwise, where $\rho$ is of a form 
$\rho=\rhopr+1$ for some $\rhopr>0$ as $F_{k+1}(0,\rho)$ cannot hold because of $\Xipr=0$. 
We therefore obtain \[\itk(\Ga[n])=\uDee\cdot(n+1)=\al[n].\] 
\item $\Xi>0$ such that $\thtke(\Xipr)<\Depre$. Then we have $\Ga=\Depre>\Omke$, but $\itk(\Ga)=\Xi+\De_1$, and $\Ga[n]=\uDepre\cdot(n+1)$. 
The first element of the $\Omke$-localization of $\Depre$ is of a form $\thtke(\Xipr+\nu)<\Depre$ for some $\nu<\Omkz$, which thus also is the 
first element of the $\Omke$-localization of $\uDepre$. Thus \[\itk(\Ga[n])=\Xi+\uDee\cdot(n+1)=\al[n].\]
\item $\thtke(\Xipr)\ge\Depre$. We then have $\Ga=\thtke(\Xipr)+\Depre$ and $\Ga[n]=\thtke(\Xipr)+\uDepre\cdot(n+1)$. 
The first element of the $\Omke$-localization of $\thtke(\Xipr)$ must be itself since $\itk(\Ga)=\Xi+\De_1$, thus again $\itk(\Ga[n])=\al[n]$.
\end{itemize}
{\bf 2.2.3.2:} $\chiomkz(\Sipr)=1$. We then have $\Depre[0]=\thtke(\Sipr[\uDepre])$ and $\Depre[n+1]=\thtke(\Sipr[\Depre[n]])$, and since $\Sipr>0$ 
we must have $\Xi,\Xipr>0$. Moreover, since $\al\in\domthtqk$, according to Lemma \ref{tqdomthtqmlem} we must have $\Si+\rho<\Xi+\De_1$, 
hence $\Si\le\Xi$ and $\Sipr\le\Xipr$ by monotonicity of $\rtke$, so that since $\itk(\Ga)=\Xi+\De_1$ the first element of the $\Omke$-localization of 
$\Depre$ must be of the form $\thtke(\Xipr+\nu)$ for some $\nu<\Omkz$.

We either have $\uDepre=\thtke(\Sipr+\rhopr)$ where $\rho=\rhopr+1$ or $\uDepre=\rho$ where $\rho$ satisfies 
$F_{k+1}(\Sipr,\rho)$. Hence $\uDepre=\gt(\uDee)$. By the i.h.\ for $\De_1$ we have $\ft(\Depre[n])=\De_1[n]$.
\begin{itemize}
\item $\thtke(\Xipr)<\Depre$. We then have $\Ga[n]=\Depre[n]$ and inspect the $\Omke$-localization of $\Depre[n]$ using Lemma \ref{localizationlem}.
The case where $\Sipr=\Omkz$ and $\rho=n=0$, hence $\Depre=\thtke(\Omkz)$ and $\Depre[0]=\Omke$, does not occur since this would imply 
$\thtke(\Xipr)\ge\Depre$.
In the cases where $\rho=\rhopr+1$ for some $\rhopr$ or $F_{k+1}(\Sipr,\rho)$ holds we immediately see that the first elements of the 
$\Omke$-localizations of $\Depre$ and $\Depre[n]$ have the same fixed point level $\Xipr$. 
In the case where $\rho=0$ the $\Omke$-localization of $\Depre$ cannot be $(\Omke,\Depre)$
since this would again imply $\thtke(\Xipr)=\Depre$. Therefore, again $\itk(\Ga[n])=\Xi+\ft(\Depre[n])=\Xi+\De_1[n]$ by the i.h.\ for $\De_1$.    
\item $\thtke(\Xipr)\ge\Depre$. Again, we have $\Ga=\thtke(\Xipr)+\Depre$ and $\Ga[n]=\thtke(\Xipr)+\Depre[n]$. The first element of the 
$\Omke$-localization of $\thtke(\Xipr)$ is itself since $\itk(\Ga)=\Xi+\De_1$, hence $\itk(\Ga[n])=\Xi+\De_1[n]$ by the i.h.\ for $\De_1$.
\end{itemize}
{\bf 2.2.3.3:} Otherwise, that is, $\Sipr>0$ such that $\chiomkz(\Sipr)=0$. This again implies that $\Xi>0$, and we argue similarly as in Subcase
2.2.3.2. The case where $\thtke(\Xipr)\ge\Depre$ is handled as before, so we may assume that $\thtke(\Xipr)<\Depre$. 
Hence $\Ga[\ze]=\Depre[\ze]=\thtke(\Sipr[\ze]+\uDepre)$. Since $\itk(\Ga)=\Xi+\De_1$,
the first element of the $\Omke$-localization of $\Depre$ is of a form $\thtke(\Xipr+\nu)$ for some $\nu<\Omkz$, and $\Sipr\le\Xipr$ since $\al\in\domthtqk$.
By Lemma \ref{localizationlem} the fixed point level of the first element of the $\Omke$-localization of $\Depre[\ze]$ must be $\Xipr$, hence
$\itk(\Ga[\ze])=\Xi+\De_1[\ze]$.
\\[2mm]
{\bf Case 3:} $\De+\eta=0$. Then we have $\rtk(\al)=\rtk(\Xi)=\Ga=\thtke(\Xipr)$, and the $\Omke$-localization of $\Ga$ is $(\Omke,\Ga)$ as 
$\itk(\Ga)=\Xi$, which also implies that $\uGa=0$. 
Note that we have $\domf(\Xipr)=\domf(\Xi)\le k+1$ by Corollary \ref{domfcor} and the assumption $d\le k+1$, hence $\Xi$ cannot be a 
successor-multiple of $\Om_{k+l}$ for any $l\ge 2$.
Thus \[\Ga[\ze]=\thtke(\Xipr[\ze])\]
where $\Xipr[\ze]$ is a nonzero multiple of $\Omkz$, and by Lemma \ref{localizationlem} the $\Omke$-localization of $\Ga[\ze]$ is $(\Omke,\Ga[\ze])$,
so that \[\itk(\rtk(\al)[\ze])=\itk(\Ga[\ze])=\itke(\Xipr[\ze])=\Xi[\ze]=\al[\ze]\] by the s.i.h.
\qed

As a corollary we now readily obtain Bachmann property for $\Tqt$-systems.
\begin{cor} The Buchholz system $(\Tqtcirc,\cdot\{\cdot\})$ satisfies Bachmann property, provided that its base system does.
\end{cor}

\emph{Norms for $\Tqt$-terms.} Note that the canonical norm $\norm{\cdot}$ as introduced in Definition \ref{normregularitydefi} cannot be directly 
carried over to $\Tqt$ even for $\tau=1$ because each term $\thtqnod(\thtqi(0))$ where $i<\om$ would receive the same norm $3$. 
This can be solved by assigning weight $i+1$ instead of $1$ to occurrences of $\thtqi$-functions or, smoother for our purposes, defining a norm for 
$\Tqt$-terms via \[\al\mapsto\norm{\gt(\al)}\] for $\al\in\Tqt$. Definition \ref{Gnormdefi} clearly carries over to $\Tqt$. 
By virtue of Theorem \ref{Tqtmaintheo} all results of Section \ref{normsec} now smoothly carry over from $\Tt$ and $\Ttcirc$ to $\Tqt$ and 
$\Tqtcirc$, respectively, which includes obtaining the same Hardy hierarchy as introduced in Definition \ref{Hardyhierarchydefi}.

\section{Conclusion}
We have established systems of fundamental sequences in the context of relativized notation systems $\Tt$ used in the analysis
of patterns of resemblance of orders $1$ and $2$ (\cite{W07b,W07c,CWa,W18,W21}), and investigated ways to relate them to the 
uniform approach to fundamental sequences and hierarchies of fast growing number theoretic functions (Hardy hierarchies) as set 
out in \cite{BCW94}. We extended our results to notation systems $\Tqt$, the definition principle of which is generalizable to larger
segments of ordinals.

Work in progress will rely on the framework established here to contribute to the theory of Goodstein sequences and to develop a theory of pattern 
related fundamental sequences, to elaborate their connection to hierarchies of fast-growing functions, and to derive independence phenomena.

\section*{Acknowledgements}
I would like to express my gratitude to the anonymous referees for both the present journal publication and the related CiE 2024 paper for very 
constructive comments and advice.
I also would like to thank Professor Andreas Weiermann for numerous instructive and motivating discussions
on several topics related to this article.

\end{document}